\documentclass[11pt]{article}
\usepackage{amssymb}
\usepackage{amsmath}

\newcommand{\nextpar}[1]{\vspace{10pt}\noindent\stepcounter{counterpar}\fbox{\textbf{\thecounterpar}}~\textbf{#1}\quad}
\newcounter{counterpar}
\newcounter{ctr}
\newtheorem{theorem}[ctr]{Theorem}
\newtheorem{proposition}[ctr]{Proposition}
\newtheorem{lemma}[ctr]{Lemma}

\newcommand{\startproof}{\noindent{\bf Proof}\quad }
\newcommand{\stopproof}{\hfill$\blacksquare$\\}
\newcommand{\stopstep}{\hfill$\Box$\\}

\newcommand{\defn}[1]{\textbf{#1}}
\newcommand{\foremph}[1]{\textit{#1}}

\newcommand{\Om}{\Omega}
\newcommand{\R}{{\mathbb R}}
\newcommand{\pom}{\partial\Om}
\newcommand{\DD}{\mathcal D}
\newcommand{\dvol}{\DD_{\rm vol}}
\newcommand{\rf}[1]{(\ref{#1})}
\newcommand{\be}{\begin{equation}}
\newcommand{\nn}{\nabla}
\newcommand{\ee}{\end{equation}}
\newcommand{\om}{\omega}
\newcommand{\dd}{\Delta}
\newcommand{\OO}{\mathcal O}
\newcommand{\EE}{\mathcal E}
\newcommand{\pul}{{1\over 2}}
\newcommand{\OOo}{\OO_{\om}}
\newcommand{\ao}{A_\om}
\newcommand{\aoi}{A_{\om}^{-1}}
\newcommand{\apsii}{A_{\psi}^{-1}}
\newcommand{\intom}{\int_\Om}

\title{Local structure of the set of steady-state solutions to the 2D incompressible Euler equations}
\author{Antoine Choffrut\thanks{Hausdorff Center for Mathematics - Universit\"at Bonn}\qquad\qquad Vladim\'ir \v Sver\'ak\thanks{University of Minnesota; supported in part by NSF grant DMS 0800908}}
\date{ }
\begin{document}
\bibliographystyle{amsplain}
\numberwithin{equation}{section}

\maketitle

\begin{abstract}
It is well known that the incompressible Euler equations can be formulated 
in a very geometric language.
The geometric structures provide very valuable
insights into the properties of the solutions.
Analogies with the finite-dimensional model of geodesics on a Lie
group with left-invariant metric can be very instructive,
but it is often difficult to prove analogues of 
finite-dimensional results in the infinite-dimensional setting
of Euler's equations.
In this paper we establish a result in this direction in the simple 
case of steady-state solutions in two dimensions,
under some non-degeneracy assumptions.
In particular, we establish, 
in a non-degenerate situation, 
a local one-to-one correspondence between steady-states
and co-adjoint orbits.
\end{abstract}

\tableofcontents

\section{Introduction}\label{section:Introduction}
\subsection{Background}
We consider the 2d Euler equations for inviscid incompressible fluid in a smooth bounded domain $\Om\subset\R^2$:
\begin{eqnarray}
\left.
\begin{array}{rcl}
u_t+ (u\cdot\nabla) u + \nabla p & = & 0 \label{i-1}\\
{\rm div}~u & = & 0
\end{array}
\right\}&  &
\mbox{in $\Om$\,,} \\
\begin{array}{ccl}
u\cdot N & = & 0\,\,\,\,\,
\end{array} & & \, \mbox{at $\pom$\,,} \label{i-2}
\end{eqnarray}
where $N$ denotes a unit normal to $\pom$. 
It is well-known that the equations have a rich geometric structure,
coming from their interpretation as equations for a geodesic flow
in the group $\dvol(\Om)$ of volume preserving diffeomorphisms.
The modern mathematical investigations exploring this geometric structure
were initiated by the 1966 paper of V.~I.~Arnold~\cite{ArnoldFourier1966}.
The geometric point of view has lead to important insights about
Euler's equations, often by analogies with the 
finite-dimensional situation of the geodesic flow on a
Lie group equipped with a left-invariant metric.~\footnote{We refer the reader to monograph~\cite{ArnoldKhesinTMH} for
examples.}
The passage from such finite-dimensional models to the infinite-dimensional
setting of Euler's equations is often impeded by a
common difficulty in infinite dimensions: the assumptions
which are needed for straightforward generalizations of basic results
of finite-dimensional Calculus (such as the Implicit Function Theorem)
are too strong to be satisfied in situations of interest.
For important advances in this direction, see for example
\cite{EbinMarsden, EbinMisiolekPreston, MisiolekPreston, Preston}.

Our goal in this paper is to establish rigorously, in certain cases, a 
geometric picture of the structure of the set of steady-states of Euler's equations~\rf{i-1},~\rf{i-2}
suggested by the finite-dimensional situation. The main theorem will be
for the case when $\Om$ is diffeomorphic to an annulus, but it seems reasonable 
to set up the problem in greater generality. 

Our results are best described
in the vorticity formulation of the equations.  
We assume that $\Om\subset\R^2$ is a bounded smooth domain.  The connected components of 
$\pom$ will be denoted by
$\Gamma_0,\Gamma_1,\dots,\Gamma_l$, with $\Gamma_0$ bordering the unbounded
connected component of $\R^2\setminus\Om$.

Any (smooth) divergence-free velocity field $u$ in $\Om$ satisfying $u\cdot N=0$ at
$\pom$ can be represented by a stream function $\psi$, defined by
\be\label{i-3}
u=\nn^\perp\psi = \left(\begin{array}{r} -\psi_{x_2}\\ \psi_{x_1}\end{array}\right)\,.
\ee
Clearly, $\psi$ is defined uniquely by $u$ up to a constant.
Therefore, without loss of generality we set $\psi|_{\Gamma_0}=0$.
The vorticity $\om$ is defined by 
\be\label{i-4}
\om=u_{2,x_1}-u_{1,x_2}=\dd\psi\,.
\ee
The stream function $\psi$ is determined by $\om$ and suitable boundary conditions.
To identify these boundary conditions, we note first from (\ref{i-2}) that,  
for each fixed time, $\psi$ is constant also on any other boundary component. 
However, the constants can depend on time, i.e. $\psi|_{\Gamma_i}$ may not, in general,
be constant during the evolution for $i=1,\dots,l$.
But by Kelvin's theorem on conservation of circulation of $u$, 
$\gamma_i=\int_{\Gamma_i} {{\partial\psi}\over{\partial N}}$ are constant along the flow.
(Note that, by the Gauss-Green theorem and the divergence-free condition in (\ref{i-1}),
the circulation $\gamma_0$ around $\Gamma_0$ is determined by $\gamma_i$, $1\leq i\leq l$.)
The constants $\gamma_i$, $1\leq i\leq l$ will be considered as fixed parameters throughout the paper.
Therefore, denoting by $\tau$ the unit tangent vector to the
boundary given by rotating the normal $N$ by $\pi/2$, our boundary conditions will be
\begin{eqnarray}
\psi|_{\Gamma_0}=  0,\, & \label{i-5a} \\ 
{{\partial\psi}\over{\partial\tau}}|_{\Gamma_i}= 0,\, & i=1,\dots,l \,, \label{i-5b}\\
\int_{\Gamma_i}{{\partial\psi}\over{\partial N}}= \gamma_i,& \,i=1,\dots,l\,.\label{i-5c}
\end{eqnarray}
We introduce the subspace
\begin{equation}
\mathcal U_{\gamma_i}=\{\psi\in C^\infty_{\overline\Omega}~\text{satisfying}~(\ref{i-5a})-(\ref{i-5c})\},
\end{equation}
of the space of stream functions
\begin{equation}
\mathcal U=\{\psi\in C^\infty_{\overline\Omega}~\text{satisfying}~(\ref{i-5a})-(\ref{i-5b})\}.
\end{equation}
Together with $\om$ and the equation
\be\label{i-6}
\dd\psi=\om
\ee
the boundary conditions (\ref{i-5a})-(\ref{i-5c}) uniquely determine $\psi$
(see \cite{MarchioroPulvirenti}).
Denoting by
\be\label{i-7}
\{f,g\}=f_{x_1}g_{x_2}-f_{x_2}g_{x_1}
\ee
the 2d Poisson bracket,
equation~\rf{i-1} can be rewritten as
\be\label{i-8}
\om_t+\{\psi,\om\}=0\,,
\ee
where $\psi$ is determined by $\om$ through~\rf{i-6} and the boundary conditions~\rf{i-5a}--~\rf{i-5c}.
Equation~\rf{i-8} describes the transport of $\om=\om(t)$ by the group $\dvol(\Om)$:
we have
\be\label{i-8b}
\om(t)=\om(0)\circ\eta^{-1}(t)\,,
\ee
where $\eta(t)\in\dvol(\Om)$ represents the particle trajectories and $\eta^{-1}$ denotes the inverse of $\eta$.
In other words, letting for each smooth function $\om$ on $\Om$
\be\label{i-9}
\OO_{\om}=\{\om\circ\eta^{-1}\,:\,\eta\in\dvol(\Om)\}\,,
\ee
and 
\be\label{i-10}
\om_0=\om(0),
\ee
the solution of~\rf{i-8} always stays on $\OO_{\om_0}$.
Moreover, equation ~\rf{i-8} can be thought of (formally) as a Hamiltonian
system on $\OO_{\om_0}$, with the Hamiltonian given by the usual energy
\be\label{i-11}
\EE(\om)=\intom\pul|\nn\psi|^2=\intom-\pul\om\psi+\sum_{j=1}^l\pul\gamma_j\psi|_{\Gamma_j}\,.
\ee 
The (formal) symplectic structure on $\OO_{\om_0}$ is of course of great
independent interest, but we will not be concerned with it in this work.
We will only study the equilibria, and these are (formally) characterized as
the critical points of the restriction of  $\EE$ to the orbits $\OO_{\om}$.~\footnote{It is 
important to point out that the space of vorticities is formally
a Poisson manifold,  not a symplectic manifold. The orbits $\OOo$ can be considered
as symplectic leaves of this Poisson manifold. 
See e.\ g.\ ~\cite{MarsdenWeinstein} for details.  See also \cite{Kirillov, Ratiu}.
 The quantities $I_f=\intom f(\om)$ are Casimir functions. They
will also be conserved if $\EE$ is replaced by any other Hamiltonian, and they
do not generate any symmetries. Such situation typically arises in the
process of symplectic reduction, and our situation is an example
of this: Euler's equations appear as a result of the reduction of the 
geodesic flow in the co-tangent bundle of $\dvol$ by the group $\dvol$.
The space of vorticities can be identified with the dual of the Lie
algebra of $\dvol$, and the orbits $\OOo$ are the orbits of the co-adjoint
representation of $\dvol$.}

To summarize, we formally have the following situation: the space of vorticities
is foliated by the orbits $\OOo$, and the equilibria are the critical
points of $\EE$ restricted to the orbits. In finite dimension a routine application of the  
Implicit Function Theorem would imply that if $\OO_{\bar\om}$ is smooth near $\bar\om$ and $\bar\om$ is a 
non-degenerate critical point of $\EE$ on $\OO_{\bar\om}$, then, near $\bar\om$, the equilibria
form a manifold transversal to the foliation by the orbits, of dimension equal to the
 co-dimension of the orbits. In other words, in a non-degenerate
situation, the equilibria
 are locally in one-to-one correspondence with the orbits.
 
Our goal is to establish an analogue of this correspondence in the context of~\rf{i-8}. 
Several obstacles have to be overcome:  the orbits 
$\OOo$ are typically not sub-manifolds of the space
of vorticities if we work with the usual Banach spaces used in PDEs; 
certain linearized operators suffer from loss of derivatives;
it is not clear what are good ``coordinates'' in which to do calculations.
Our main result is, roughly speaking, that the difficulties can 
be resolved for steady-states for which the vorticity $\om$ has no
critical points in $\bar\Om$. This assumption is of course 
restrictive. However, it is likely that in some situations the 
critical points of $\om$ can genuinely  complicate the 
picture and lead to some degeneracies,
especially in the case of hyperbolic critical points.
Elliptic critical points seem to be less dangerous. They still
lead to difficulties for our method (namely certain linear maps are no longer ``tame'',
see Section~\ref{section:IntroNashMoserIFT}), 
but these might perhaps be manageable.

We now outline the main points of our approach. We start with the classical
observation that any function $\psi$ satisfying
\be\label{i-12}
\dd\psi=F(\psi)\,
\ee 
with the boundary conditions~\rf{i-5a}--~\rf{i-5c} gives a steady-state.
This is easily seen from~\rf{i-8}. 
Moreover, if $\bar\psi$ and $\bar F$ solve~\rf{i-12},\,\,\rf{i-5a}--~\rf{i-5c},
and if $\bar\om=\bar F(\bar \psi)$ has no critical points, then any nearby steady-state 
can be obtained in this way. 
This is one reason for the restriction on the geometry of $\Omega$:
the condition that $\bar \omega$ has no
critical points imply that $\Omega$ is diffeomorphic to an annulus,
as $\bar\omega$ is constant on $\partial \Omega$.
(For the case where $\omega$ has a single elliptic critical point in a simply connected domain $\Omega$,
a refinement of our method seems necessary.)
The boundary of $\Omega$ consists then of an inner and outer boundary components,
\begin{equation}\partial\Omega=\Gamma_i\cup\Gamma_o.\end{equation}
This is one characterization
of the steady-states in the situation we wish to investigate: we see that
they are, in some sense, locally parametrized by the functions $F$
such that
\begin{equation}F'\neq 0.\end{equation}
The map $F\mapsto\psi$, defined via~\rf{i-12} and~\rf{i-5a}--~\rf{i-5c} is not quite
one-to-one, since changes to $\bar F$ outside the range of $\bar\psi$ do not 
change the solution $\bar\psi$, but this is not a serious problem.

It is worth remarking that steady-states of the form (\ref{i-12}) with $F'\neq 0$
naturally arise in the statistical theories of 2D Euler flows,
see \cite{MillerWeichmanCross, Robert, Turkington}, 
and in Shnirelman's theory of mixing, see \cite{Shnirelman}.

Our plan is to establish the correspondence between the functions $F$
in~\rf{i-14} and the orbits $\OOo$. By contrast with the finite-dimensional
case, there is a simple obstacle showing that we cannot really consider all
orbits $\OOo$: for any steady-state solution satisfying~\rf{i-5a}--~\rf{i-5c}
the vorticity $\om$ must be constant along each boundary component. 
We introduce the space
\begin{equation}
\mathcal F=\{\omega\in C^\infty_{\overline\Omega}~|~\frac{\partial\omega}{\partial\tau}_{|\partial\Omega}=0\}.
\label{eq:FunctionsLocallyConstantOnBoundary}
\end{equation}
Our main result (Theorem~\ref{theorem:MainTheorem}) establishes a local correspondence
between the functions $F$ and the orbits contained in $\mathcal F$, 
i.e. with constant values of $\om$ at the boundary components.

We need to introduce some (local) parametrization of the ``space of orbits''
$\OOo$. The distribution functions $\ao$ defined by
\be\label{i-13}
\ao(\lambda)=|\{x\in\Omega~|~\om(x)<\lambda\}|
\ee     
provide a good option for those orbits $\mathcal O_\omega\subset \mathcal F$ 
with constant values of $\omega$ on the boundary components,
and for which $\omega$ has no critical points, see Proposition~\ref{proposition:AomegaCharacterizesO*(omega)}.
(The assumption that $\omega$ has no critical points surfaces again and  imposes on $\Omega$ to be diffeomorphic to an annulus.)
In fact it will be better to work with the inverses $\aoi$, for several reasons.
One is that their domain, which is the interval $[0,|\Om|]$, does not change,
and another is the identity
\be\label{i-14}
\aoi=F\circ\apsii\,
\ee
satisfied by the solutions of~\rf{i-12}
in the case when $F'>0$ (when $F'<0$ a similar identity holds; 
in the proof we will perform computations assuming $F'>0$ for simplicity,
the case $F'<0$ being completely analogous). 
This identity will be crucial
for the analysis for the following reason: it shows that there is a chance
for establishing some correspondence between $F$ and $\aoi$. 
The right-hand side of (\ref{i-14})  is non-linear
in $F$ (as $\psi$ depends on $F$), but the non-linear part comes in
only through $\psi$ and hence it is regularized by equation~\rf{i-12}. 
In some sense, the leading part of the dependence of $\aoi$ on $F$
behaves as a composition of $F$ with a fixed function, which, for many purposes, is 
almost the same as identity.  
Heuristically, the function $\aoi$ can be thought of as being obtained from $F$ by applying a 
kind of non-linear Fredholm map to $F$. This ``Fredholmness'' of the map $F\mapsto\aoi$ is crucial
for our approach.
On the other hand, the linearization of $A_\omega^{-1}$ suffers from loss of derivatives,
and \rf{i-14} shows it quite clearly:
\be\label{i-15} \delta (F\circ\apsii)=(\delta F)\circ \apsii+(F'\circ\apsii)\cdot\delta(\apsii)\,.
\ee  
The first term on the right-hand side is very good, but the second term 
contains $F'$. This does not seem to be easy to avoid, and it can be overcome
by working with the Nash-Moser Implicit Function Theorem, which will enable
us to establish a good local correspondence between $F$ and $\aoi$ mentioned above,
see Theorem~\ref{theorem:MainTheorem}.
The correspondence between $F$ and $A_\omega^{-1}$ cannot be one-to-one, as already noted above, but only for the
trivial reason that in~\rf{i-12} the behavior of $F$ outside the range of $\psi$
does not affect the solution. 
Theorem~\ref{theorem:MainTheorem} does contain an injective part, though,
which says
that, if two nearby steady-states have same distribution functions, then they are identical.
Also, some natural non-degeneracy assumptions are needed,
in the form of transversality conditions for linearized operators. Some conditions
of this form are needed even in the finite-dimensional situation.   

\subsection{Statement of main result}
We first introduce the two non-degeneracy assumptions of Theorem~\ref{theorem:MainTheorem}.

We will denote $\bar\psi$, $\bar F$, and $\bar\omega=\overline F(\overline \psi)$ the quantities associated with a reference steady-state.
The function $F$ in (\ref{i-12}) gives a good (local) parametrization of steady-states when $F'\neq 0$
(modulo the lack of injectivity mentioned above).
For this reason $\Omega$ is assumed to be diffeomorphic 
to an annulus.
However, a well-defined map $F\mapsto \psi$ returning a solution to
\begin{equation}\Delta \psi=F(\psi),\qquad \psi_{|\Gamma_o}=0,\quad \frac{\partial\psi}{\partial\tau}_{|\Gamma_i}=0,
\quad \int_{\Gamma_i}\frac{\partial\psi}{\partial N}=\gamma_i,\quad 1\leq i\leq l
\label{eq:SteadyStateEquation}
\end{equation}
can be constructed in a neighborhood of $\overline F$, modulo some non-degeneracy condition discussed shortly,
for an arbitrary number $l$ of boundary components.
Thus for the construction of the map $F\mapsto \psi$ we will make no restriction on the topology of $\Omega$.

The map $F\mapsto \psi$ is well-defined provided the linear map
\begin{equation}
\Delta\phi-F'(\psi)\phi=k,\qquad \phi_{|\Gamma_o}=0,\quad \frac{\partial\phi}{\partial\tau}_{|\Gamma_i}=0,\quad \int_{\Gamma_i}\frac{\partial\phi}{\partial N}=0
\label{eq:Delta-F'(psi)}
\end{equation}
is invertible for each $k\in C^\infty_{\overline\Omega}$ and 
\foremph{for each $F$ in a neighborhood of $\overline F$}.
In fact, it is enough to make this assumption \foremph{at the reference steady-state only}:
\begin{equation}
(\text{ND1})\qquad
\left\{\begin{array}{l}
\text{the reference steady-state}~ \overline\omega=\overline F(\overline\psi)~\text{is non-degenerate in the sense that}
\notag\\\\
\Delta\phi-\overline F'(\overline \psi)\phi=0,\qquad \phi_{|\Gamma_o}=0,\quad \frac{\partial\phi}{\partial\tau}_{|\Gamma_i}=0,\quad \int_{\Gamma_i}\frac{\partial\phi}{\partial N}=0, ~1\leq i\leq l
\notag\\\\
\text{has only the trivial solution}~\phi\equiv 0.\notag
\end{array}\right.
\end{equation}
By the Fredholm alternative, the operator $\Delta -\overline F'(\overline\psi)$ is invertible
with the boundary conditions of (ND1).
Since we do not work in a Banach-space setting but with Fr\'echet spaces, 
it is not automatic that (ND1) implies that $\Delta-F'(\psi)$ is invertible in $C^\infty$ for $F$ near $\overline F$.
See \cite{HamiltonIFTNashMoser} for a discussion of this crucial point and counterexamples,
in particular Section~I.5.5 of Part~I.
\bigskip%

The second non-degeneracy condition says that the steady-states and the co-adjoint orbits intersect trivially.
(Since now we are using the parametrization of the ``space of orbits'' via the distribution functions $A_\omega$,
we will need to impose that $\Omega$ is diffeomorphic to an annulus since then it is crucial that $\omega$ has no critical points.)
More precisely, a tangent $\nu$ to the set of steady-states at $\omega$ is a solution
to the linearized steady-state equation $\omega=F(\psi)$, i.e.  $\nu=\Delta\phi$ where $\phi$ solves, for some $f$,
\begin{equation}
\Delta \phi = F'(\psi)\phi+f(\psi),\qquad \phi_{|\Gamma_o}=0,\qquad \frac{\partial\phi}{\partial\tau}_{|\Gamma_i}=0,
\qquad \int_{\Gamma_i}\frac{\partial\phi}{\partial N}=0.
\label{eq:LinearizedSteadyState}
\end{equation}
A function $\nu$ is tangent to the co-adjoint orbit $\mathcal O_\omega$ at $\omega$ if
there exists a stream function $\alpha\in \mathcal U$ such that $\nu=\{\omega,\alpha\}$.
(This is immediate by linearizing $\omega\circ \eta$ at $\eta=\text{Id}$.)
The second non-degeneracy condition is as follows:
\begin{equation}
(\text{ND2})\qquad
\left\{\begin{array}{l}
\text{if}~\phi~\text{satisfies the linearized steady-state equation (\ref{eq:LinearizedSteadyState}) at}~\overline\omega=\overline F(\overline\psi)
\notag\\\\
\text{and if} ~\nu=\Delta\phi~ \text{is tangent to}~\mathcal O(\overline\omega)~\text{at}~\overline\omega,~\text{i.e.}~ \nu=\{\overline\omega,\alpha\}~\text{for some}~\alpha\in\mathcal U,\notag\\\\
\text{then}~\phi=0.\notag
\end{array}\right.
\end{equation}
Again, we emphasize that this non-degeneracy assumption is made 
\foremph{at the reference steady-state only and not in an entire neighborhood}.
\bigskip%

\begin{theorem}\label{theorem:MainTheorem}
Let $\Omega\subset \mathbb R^2$ be diffeomorphic to an annulus with inner and outer boundary components
$\Gamma_i$, $\Gamma_o$ respectively.
Consider a smooth steady-state solution to Euler's equation on $\Omega$
with vorticity $\overline\omega$ and stream function $\overline\psi$ without critical points.
Let $\overline F$ and $\gamma_i$  such that
\begin{equation}\Delta \overline\psi=\overline F(\overline\psi),\qquad {\overline\psi}_{|\Gamma_0}=0,\qquad \frac{\partial\overline\psi}{\partial\tau}_{|\Gamma_i}=0,
\qquad \gamma_i=\int_{\Gamma_i}\frac{\partial \overline\psi}{\partial N}\end{equation}
(in particular $\overline F'\neq 0$).
Assume further that (ND1) and (ND2) are satisfied.
Then, there exists a neighborhood $\mathcal W$ of $\overline\omega$ in $C^\infty_{\overline\Omega}\cap \mathcal F$
such that 
each co-adjoint orbit intersecting $\mathcal W$ contains exactly one smooth steady-state solution there.
\end{theorem}

\noindent\textbf{Remark}\qquad $\mathcal W$ can be taken to be a $\|\cdot\|_{11}$-neighborhood.
See Proposition~\ref{proposition:InjectivePartOfMainTheorem} (proving the injective part of Theorem~\ref{theorem:MainTheorem})
and the Remark after the statement of Theorem~\ref{theorem:SurjectivePartOfNashMoserIFT} in the Appendix
(on the existence part of the Nash-Moser Inverse Function Theorem).
\bigskip%

An important aspect of the proof of Theorem~\ref{theorem:MainTheorem} is that both non-degeneracy
assumptions (ND1) and (ND2), which are made at the reference steady-state, 
are sufficient to imply non-degeneracy for steady-states \foremph{in an entire neighborhood}.
It turns out that (ND1) and (ND2) are of exactly the same type,
and the proof that $DT(F)f$ has a tame right-inverse will parallel the proof that 
$\Delta\phi-F'(\psi)\phi=k$ (with the boundary conditions of (\ref{eq:LinearizedSteadyState})) has a tame inverse.
There is a kind of Fredholmness  at work in both cases.
To understand (ND1), consider the linear elliptic equation
\begin{equation}\Delta\phi+c\phi=k\label{eq:Delta(phi)+cphi=k}\end{equation}
parametrized by $c$ and 
with the same boundary conditions as in (\ref{eq:LinearizedSteadyState}).
The second term $K(c)\phi=c\phi$ is a ``compact'' perturbation of $\Delta\phi$ in the sense that, for each $n\geq 0$,
$c\phi$ is $C^{n+2,\alpha}$ when $\Delta\phi$ is $C^{n,\alpha}$ (assuming $c$ is smooth).
Assume that $\overline c$ is non-degenerate in the sense that $\Delta+\overline c$ has trivial kernel
in $C^\infty$.
Then, $\Delta+c$ is also non-degenerate if $\|c-\overline c\|_{0,\alpha}$ is sufficiently small.
This comes from the estimate
$\|(c-\overline c)\phi\|_{n,\alpha}\leq C\|c-\overline c\|_{0,\alpha}\|\phi\|_{n,\alpha}+C\|c-\overline c\|_{n,\alpha}\|\phi\|_{0,\alpha}$.
In addition, the elliptic estimates are easily converted into tame estimates for the inverse.

The main thrust in using the Nash-Moser theorem (see Theorem~\ref{theorem:SurjectiveNashMoserIFT}) 
is to show that the derivative $DT$ has a tame \foremph{right}-inverse.
For the problem at hand, this boils down to showing that a map of the form
\begin{equation}h=g+K(F)g\end{equation}
has a tame \foremph{inverse} 
(see Proposition~\ref{proposition:TameInverseToh=g+K(F)g} of Section~\ref{section:DT(F)IsSurjectiveWithSmoothTameRightInverse}).
This is possible since, here again, the second term $K(F)g$ is a ``compact perturbation''
in the sense that, for each $n\geq 0$, $K(F)g$ is $C^{n+2,\alpha}$ when $g$ is $C^{n,\alpha}$
(and $F$ is smooth).
On the other hand, estimates on $\|K(F)g-K(\overline F)g\|_{n,\alpha}$ are
considerably more difficult to establish than those on $\|(c-\overline c)\phi\|_{n,\alpha}$.
Again, estimates for $h=g+K(F)g$ yield easily tame estimates for the inverse.
\bigskip%

\noindent\textbf{The special case $F'>0$}\\
When $F'>0$, it is clear that (ND1) is automatically satisfied at $\omega=F(\psi)$
(multiply (\ref{eq:LinearizedSteadyState}) and integrate by parts).
It turns out that (ND2) is also automatically satisfied.
Let $\nu=\{\omega,\alpha\}$ satisfy (\ref{eq:LinearizedSteadyState}).
Then
\begin{equation}
\int_\Omega\frac{\nu^2}{F'(\psi)}=\int_\Omega\phi\nu +\int_\Omega \frac{f(\psi)}{F'(\psi)}\{\omega,\alpha\}.
\end{equation}
The second term of the right-hand side vanishes using (\ref{eq:Intf{g,h}dvol=etc.}).
Integrating by parts,
\begin{equation}\int_\Omega \left(\frac{\nu^2}{F'(\psi)}+|\nabla\phi|^2\right) = 0
\end{equation}
which forces $\phi=0$ when $F'>0$.

Theorem~\ref{theorem:MainTheorem} is a statement about the \textit{local} structure of the set of steady-states.
However, the proof suggests that a \textit{global} statement should hold 
in the case where $\overline F'>0$.
Namely one might speculate that the entire collection of steady-states satisfying $F'>0$ 
is in one-to-one correspondence with their co-adjoint orbits.
In other words, on any orbit containing a steady flow with $F'>0$, this flow might be unique with this property.
Note that this is certainly true for radial flows
since the profile of $\omega$ and $A_\omega$ are the same up to a change of variables.
In this case, (\ref{eq:SteadyStateEquation}) reduces to a second order ODE with two boundary conditions.
\bigskip%
\pagebreak

\subsection{The Nash-Moser Inverse Function Theorem and tame estimates}\label{section:IntroNashMoserIFT}
The Nash-Moser Inverse Function Theorem for tame Fr\'echet spaces (see Theorem~\ref{theorem:SurjectiveNashMoserIFT})
bears certain important differences
with the classical Inverse Function Theorem in Banach spaces.
These differences concern in particular the notion of differentiability (which is in a sense weaker than the usual one for Banach spaces),
the notion of tameness, and the fact that a right-inverse to the first derivative is assumed to exist in an entire neighborhood
of $\overline F$ and not just at $\overline F$.
The aim of this Section is to address these differences by giving precise definitions 
and clarifying certain assumptions.
See Theorem~\ref{theorem:SurjectivePartOfNashMoserIFT} in the Appendix for the existence part of the Nash-Moser Inverse Function Theorem.

Our terminology and definitions mosly follow \cite{HamiltonIFTNashMoser}.
\bigskip%

\noindent\textbf{Smooth maps of Fr\'echet spaces}\\
Let the topology on a  Fr\'echet space $\mathcal X$ be defined by a countable family of semi-norms $|\cdot|_n$, $n=0,1,2,\dots$
We will be exclusively concerned with spaces $C^\infty_{\mathcal K}$ of smooth functions
on a compact manifold $\mathcal K$ (possibly with boundary),
and the semi-norms will be either the norms $\|\cdot\|_n$-norms (i.e. the sup-norms of derivatives up to $n$-th order)
or the $\|\cdot\|_{n,\alpha}$-norms (i.e. the H\"older-norms of derivatives up to $n$-th order).
Two \defn{gradings} on $\mathcal X$ (i.e. two families of semi-norms $|\cdot|_n$ and $|\cdot|_n'$, $n=0,1,2,\dots$)
are \defn{equivalent} if they define the same topology.

Let now $\mathcal X, \mathcal Y$ be Fr\'echet spaces with semi-norms $|\cdot|_n$, $n=0,1,2,\dots$
(for simplicity, we will use the same notation for the gradings of $\mathcal X$ and $\mathcal Y$).
Let $\mathcal B$ be an open set and $P\colon (\mathcal B\subset \mathcal X)\rightarrow \mathcal Y$ a map between these Fr\'
echet spaces.
\defn{Continuity} is defined as usual.
In particular, when $\mathcal X$ and $\mathcal Y$ are of the form $C^\infty_{\mathcal K}$,
then $P$ is continuous on $\mathcal B$ if, for each $n$, there exists $m=m(n)$ such that 
$P\colon (\mathcal B,|\cdot|_{m(n)})\rightarrow (\mathcal Y,|\cdot|_n)$ is continuous.
Also, two gradings are equivalent if the identity maps
\begin{equation}\text{Id}\colon (\mathcal X, \{|\cdot|_n\}_n)\rightarrow (\mathcal X,\{|\cdot|_n'\}_n),
\qquad \text{Id}\colon (\mathcal X, \{|\cdot|_n'\}_n)\rightarrow (\mathcal X,\{|\cdot|_n\}_n)
\label{eq:IdentityMapsOnFrechetSpaces}
\end{equation}
are both continuous.
Clearly the $C^k$- and the $C^{k,\alpha}$-gradings are equivalent on $C^\infty_{\mathcal K}$.
The corresponding topology is called the $C^\infty$-topology.
For the spaces of the form $C^\infty_{\mathcal K}$,
we will prove continuity of maps using whichever grading is more convenient
(the $C^{k,\alpha}$-grading for operators involving elliptic equations, the $C^k$-grading otherwise).

The notion of differentiability, on the other hand, is in a sense weaker than the usual notion for maps of Banach spaces.
The map $P$ is \defn{differentiable at $u\in\mathcal B$} if for each $v\in\mathcal X$ the limit
\begin{equation}DP(u)v:=\lim_{t\rightarrow 0}\frac{P(u+tv)-P(u)}t\end{equation}
exists in the Fr\'echet-topology, that is, there exists an element $DP(u)v\in \mathcal Y$ such that
\begin{equation}
\lim_{t\rightarrow 0}\left|\frac{P(u+tv)-P(u)}t-DP(u)v\right|_n=0\qquad \text{for each}\quad  n=0,1,2,\dots
\end{equation}
In other words, $P$ has G\^ateaux-derivatives at $u$ in all directions.
$P$ is \defn{continuously differentiable in $\mathcal B$} if 
the map 
\begin{equation}
\left\{
\begin{array}{ccc}
(\mathcal B\subset \mathcal X)&\times& \mathcal X\\
u&&v\end{array}
\longrightarrow \begin{array}c\mathcal Y\\DP(u)v\end{array}
\right\}
\end{equation}
is continuous  \foremph{jointly in the two variables $u$ and $v$}.
In the case the spaces are Banach spaces,
this definition of differentiability is weaker than continuous Fr\'echet-differentiability,
which is a usual assumption for the classical Inverse Function Theorem there.
Partial derivatives for maps of several variables are defined in the usual way,
as well as derivatives of higher order, e.g. (when they exist)
\begin{equation}D^2P(u)(v_1,v_2)=\lim_{t\rightarrow 0}\frac{DP(u+tv_2)v_1-DP(u)v_1}t
\end{equation}
the limit again taken in the Fr\'echet-topology.
A map is smooth if derivatives of all orders exist and are continuous.
All maps will turn out to be smooth, but we will only need at most two derivatives
in order to apply the Nash-Moser Inverse Function Theorem, see Theorem~\ref{theorem:SurjectiveNashMoserIFT}.
(Specifically, we will prove that $T$, $DT$, $D^2T$, and a right-inverse $L$ to $DT$ are continuous.)
Thus, we will call these maps smooth even though we only establish that they have
continuous derivatives of order at most two.

A piece of terminology.
If a map $P(u,v)$ is linear in $v$, then we will say that it is a 
\defn{family of linear maps} and write it as $P(u)v$ to emphasize linearity in $v$.
(Similar terminology and notations apply for maps of more than two variables.)

Many rules of the usual calculus apply.
In particular, the first derivative $DP(u)v$ is linear in $v$
(see \cite{HamiltonIFTNashMoser}, Section~3.2 of Part~I)
the chain rule and the Fundamental Theorem of Calculus hold
(see \cite{HamiltonIFTNashMoser}, Section~2 of Part~I, for a definition and properties of Fr\'echet-space valued integrals),
as well as Taylor's formula with integral remainder
(see Theorem~3.5.6, p.~82 in Part~I of \cite{HamiltonIFTNashMoser}):
\begin{equation}P(u+v)=P(u)+DP(u)v+\int_0^1(1-t)D^2P(u+tv)(v,v)dt.
\label{eq:TaylorSecondOrder}
\end{equation}
Also, the Open Mapping Theorem holds: if a continuous linear map of Fr\'echet spaces is invertible,
then it is a linear isomorphism, i.e. its inverse is again a continuous map of Fr\'echet spaces.
On the other hand, if for a smooth family of linear maps, $P(\overline u)v$ is invertible,
then it is not true in general that $P(u)v$ has an inverse for $u$ in a neighborhood of $\overline u$.
This is in contrast to the Banach-space setting, where the set of invertible operators is open
(this is related to the fact that the set of bounded linear maps on Banach spaces is itself a Banach space).
This explains why the invertibility of the first derivative $DT(u)v$ for $u$ in an entire
neighborhood must be assumed in order to apply the Nash-Moser Inverse Function Theorem.
\bigskip%
\pagebreak

\noindent\textbf{The tame Fr\'echet category}\\
For our purposes (see \cite{HamiltonIFTNashMoser}, Section~II.1,  for a more general notion), 
a \defn{tame Fr\'echet space} $\mathcal X$ is a Fr\'echet space,
with semi-norms $|\cdot|_n$, $n=0,1,2,\dots$, 
which comes with
a family of \defn{smoothing operators} $\{S(t)\}_{t>0}$ such that, for all $t>0$ and $u\in\mathcal X$,
\begin{equation}
|S(t)u|_m\leq C t^{m-l}|u|_l,\qquad |u-S(t)u|_l\leq C t^{l-m}|u|_m,\qquad m\geq l,
\label{ineq:IntroEstimatesSmoothingOperators}
\end{equation}
where the constants depend on $m, l$, but not on $t$ nor $u$.
Spaces of the form $C^\infty_{\mathcal K}$ are tame, see \cite{HamiltonIFTNashMoser}, Part~II, Theorem~1.3.6, p.~137 and Corollary~1.3.7, p.~138.
It is interesting, see \cite{Sergeraert}, as well as Corollary~1.4.2, p.~176,  Part~II of \cite{HamiltonIFTNashMoser},
that the estimates (\ref{ineq:IntroEstimatesSmoothingOperators}) 
\foremph{imply} the interpolation inequalities
\begin{equation}
|u|_i \leq C|u|_{m}^\frac{l-i}{l-m}|u|_{l}^\frac{i-m}{l-m},\qquad m\leq i\leq l
\label{ineq:IntroInterpolationInequalities}
\end{equation}
where the constants depend on $i, m, l$.
These inequalities can otherwise be verified directly ``by hand'' in the $C^k$-, $C^{k,\alpha}$-, or $H^k$-gradings,
for example,
see \cite{delaLlaveObaya}, \cite{HamiltonIFTNashMoser} (Theorem~2.2.1, p.~143, Part~II), \cite{KrylovEllipticParabolicHoelderSpaces}.
\bigskip%

A continuous map $P\colon (\mathcal B\subset \mathcal X)\rightarrow \mathcal Y$ of tame Fr\'echet spaces
is \defn{tame} if, for each $u_0\in \mathcal B$ there exist a neighborhood $\mathcal V$ of $u_0$ in $\mathcal B$,
$r\in\mathbb N$ (the \defn{degree}), $b\in\mathbb N$ (the \defn{base}), and constants $C_n$
such that
\begin{equation}
|P(u)|_n\leq C_n (|u|_{n+r}+1), \qquad n\geq b\label{ineq:TameEstimates}
\end{equation}
for any $u\in\mathcal V$.
These are called \defn{tame estimates} for $P$.
We will usually suppress the dependence on $n$ for the constants and simply write $C$
(even though this dependence on $n$ is, of course, crucial).
It can be proven, see Proposition~\ref{proposition:EquivalentFormsOfTameEstimatesForFamilyOfLinearMaps} below,
that, if a \foremph{linear} map $L$ is tame with degree $r$ and base $b$, then tame estimates can be derived in the form
\begin{equation}
|Lu|_n\leq C |u|_{n+r},\qquad n\geq b\label{ineq:TameLinearEstimates}
\end{equation}
for all $u\in\mathcal X$ \foremph{without any restriction}.
A map is \defn{smooth tame} if it is smooth and derivatives of all orders are tame.

We say that a grading $\{|\cdot|_n'\}_n$ is \defn{tame equivalent} to $\{|\cdot|_n\}_n$
if the identity maps (\ref{eq:IdentityMapsOnFrechetSpaces}) 
are both tame.
(Obviously, this defines an equivalence relation.)
In this case, the smoothing operators $S(t)$ satisfy again inequalities of the form (\ref{ineq:IntroEstimatesSmoothingOperators})
with $|\cdot|_n$ replaced by $|\cdot|_n'$.
A map $P\colon (\mathcal B\subset \mathcal X)\rightarrow \mathcal Y$
of tame Fr\'echet spaces remains tame if one replaces the gradings on $\mathcal X$ and $\mathcal Y$ with tame
equivalent gradings.
Note though that the degrees may be different in different gradings, and thus the choice of grading in 
which the tame estimates are derived should be made with care.
We will derive tame estimates for all maps in the same grading,
and our choice will be the $C^{n,\alpha}$-grading in order to take full advantage of elliptic regularity
afforded by the elliptic system (\ref{eq:SteadyStateEquation}).
(Clearly, the $C^n$- and $C^{n,\alpha}$-gradings on $C^\infty_{\mathcal K}$ are tame equivalent.)
\bigskip%

\noindent\textbf{Further remarks on tame estimates}\\
Clearly, a tame map with degree $r$ also has degree $r'\geq r$.
In turn, it is possible to choose the neighborhood $\mathcal V$ in which
(\ref{ineq:TameEstimates}) holds to be a $|\cdot|_{r+b}$-neighborhood
(by making $r$ larger and $\mathcal V$ smaller).
Composition of tame maps is again a tame map.

On the other hand, the Open Mapping Theorem does not hold in the tame Fr\'echet category:
if a tame linear map is invertible, then an inverse exists and is continuous,
but it need not be tame.
See \cite{HamiltonIFTNashMoser}, Section~1.5.5, Part~I, for counterexamples.

For a map of several variables, these are allowed to have different degrees, e.g.
\begin{equation}|P(u_1,u_2)|_n\leq C_n(|u_1|_{n+r_1}+|u_2|_{n+r_2}+1),
\quad |u|_{r_1+b}<\delta_1,\quad |u|_{r_2+b}<\delta_2,\quad n\geq b.
\end{equation}

For a family of linear maps,  it is possible to do away with the restriction on the variables in which the map is linear
(see Lemma~2.1.7, p.~143, Part~II of \cite{HamiltonIFTNashMoser})
\begin{proposition}[Tame estimates for families of linear maps]\label{proposition:EquivalentFormsOfTameEstimatesForFamilyOfLinearMaps}
Let $P(u)v$ be a family of linear maps. Then there exist constants $C_n$ such that 
\begin{equation}
|P(u)v|_n\leq C_n(|u|_{n+r}+|v|_{n+s}+1),\qquad n\geq b
\label{ineq:TameEstimatesForFamilyOfLinearMaps-Equivalent}
\end{equation}
for $u$ in a $|\cdot|_{r+b}$-neighborhood and $v$ in a $|\cdot|_{s+b}$-neighborhood
if and only if there exist constants $C_n'$ such that
\begin{equation}|P(u)v|_n\leq C_n'(|u|_{n+r}|v|_s+|v|_{n+s}),\qquad n\geq b
\label{ineq:TameEstimatesForFamilyOfLinearMaps}
\end{equation}
for $u$ in a $|\cdot|_{r+b}$-neighborhood and \foremph{any} $v$ (without restriction).
\end{proposition}
This generalizes to maps linear in more than one variable.
In our proof, we will systematically derive tame estimates in the format (\ref{ineq:TameEstimatesForFamilyOfLinearMaps})
when relevant.
\bigskip%

\noindent\textbf{The Nash-Moser Inverse Function Theorem}
\begin{theorem}[Surjective part of the Nash-Moser theorem]\label{theorem:SurjectiveNashMoserIFT}
Let $T\colon(\mathcal B\subset \mathcal X)\rightarrow \mathcal Y$ be a map of tame Fr\'echet spaces.
Suppose that $T$ possesses first and second derivatives $DT$ and $D^2T$, that $DT$ has a right-inverse $L$,
and that all these maps are continuous and satisfy tame estimates.
Then, for any $x_0\in\mathcal B$,
there exists a neighborhood $\mathcal V(x_0)$ of $x_0$, a neighborhood $\mathcal V(y_0)$ of $y_0=T(x_0)$,
and a map $R\colon \mathcal V(y_0)\rightarrow \mathcal V(x_0)$ such that
$T(R(y))=y$ for $y\in \mathcal V(y_0)$.
Furthermore, $R$ is continuous and tame, and if $T$ and $L$ are smooth tame, then so is $R$.
\end{theorem}
See the proof of Theorem~\ref{theorem:SurjectivePartOfNashMoserIFT}
for the existence part of the surjective Nash-Moser Inverse Function Theorem.
We refer to \cite{HamiltonIFTNashMoser}, Section~III.1  for proofs
of further properties of $R$ (smoothness and tameness).
\bigskip%

The proof Theorem~\ref{theorem:MainTheorem} will heavily rely on smoothness and tameness of elementary maps of Fr\'echet spaces.
The necessary lemmas are given in the Appendix and will be used countless times,
often without explicit reference.
\bigskip%

\subsection{Examples of the ``orbit calculus''}
A rigorous interpretation of the orbits $\mathcal O_\omega$ as symplectic leaves
would require some care.
Instead, we give in this Section two examples of the ``orbit calculus''.
Both results go back to Arnold \cite{ArnoldFourier1966}, \cite{ArnoldKhesinTMH},
but our calculations here are slightly different and do not resort to Lie-group theoretical interpretations.

Let $\omega = \Delta \psi$ be a critical point of the kinetic energy
\begin{equation}
{\mathcal E}(\omega)=\frac12\int_\Omega|\nabla\psi|^2\text{dvol}=-\frac12\int_\Omega \omega\psi\text{dvol} +\frac12\sum_{i=1}^l\gamma_i\psi_{|\Gamma_i}
\end{equation}
restricted to its co-adjoint orbit
(the $\gamma_i$'s are fixed).
Let $\nu=\{\alpha,\omega\}$ be tangent to the orbit $\mathcal O_\omega$ at $\omega$,
and $\phi\in\mathcal U_0$ solving (\ref{eq:LinearizedSteadyState}).
Then the first derivative of the energy is given by
\begin{equation}
D{\mathcal E}(\omega)\nu=\int_\Omega\nabla \psi\cdot\nabla \psi\text{dvol}
=-\int_\Omega \psi\nu\text{dvol}+\int_{\partial\Omega}\psi\frac{\partial\phi}{\partial N}dl
=-\int_\Omega \psi\{\alpha,\omega\}\text{dvol}.
\end{equation}
The identity
\begin{equation}
\int_\Omega f\{g,h\}\text{dvol}=\int_\Omega g\{h,f\}\text{dvol}-\int_{\partial\Omega}fg\frac{\partial h}{\partial N}dl
\label{eq:Intf{g,h}dvol=etc.}
\end{equation}
gives
\begin{equation}
D{\mathcal E}(\omega)\nu=\int_\Omega \alpha\{\psi,\omega\}\text{dvol}
\end{equation}
(the boundary terms vanish).
Since $\alpha$ is arbitrary, we conclude that $\{\psi,\omega\}=0$.

Next we compute the second variation of ${\mathcal E}$ at a critical point:
\begin{equation}
D^2{\mathcal E}(\omega)(\nu,\nu)=\int_\Omega\left(|\nabla\phi|^2+\frac{\nu^2}{F'(\psi)}\right)\text{dvol}
\label{eq:SecondVariation}
\end{equation}
where $\psi$ solves (\ref{eq:SteadyStateEquation})
and $\phi$ solves (\ref{eq:LinearizedSteadyState}).
\bigskip%

\noindent\textbf{Proof of (\ref{eq:SecondVariation})}\qquad
Let $\omega_\epsilon=\omega\circ \eta_\epsilon$ with $\omega_0=\text{id}$,
and denote $\psi_\epsilon$ the corresponding stream functions.
From the first derivative
\begin{equation}\frac{d}{d \epsilon}{\mathcal E}(\omega_\epsilon)
=\int_\Omega\nabla\psi_\epsilon\cdot \frac{\partial}{\partial \epsilon}\nabla\psi_\epsilon\text{dvol}
\end{equation}
we find that the second derivative at $\epsilon=0$ is,
posing $\dot\psi=\frac{\partial\psi_\epsilon}{\partial \epsilon}_{|\epsilon=0}$
and $\ddot\psi=\frac{\partial^2\psi_\epsilon}{\partial \epsilon^2}_{|\epsilon=0}$,
\begin{equation}
\frac{d^2}{d\epsilon^2}_{|\epsilon=0}{\mathcal E}(\omega_\epsilon)
=\int_\Omega(|\nabla \dot\psi|^2+\nabla\psi\cdot\nabla \ddot\psi)\text{dvol}.
\end{equation}
The first term is $\int_\Omega |\nabla\phi|^2\text{dvol}$.
As for the second, integrating by parts we find
\begin{equation}
\int_\Omega \nabla\psi\cdot\nabla \ddot\psi\text{dvol}=-\int_\Omega\psi\ddot\omega\text{dvol}
\end{equation}
(the $\gamma_i$'s are fixed).
To calculate $\ddot\omega$, we can take $\eta_\epsilon$ as the flow corresponding to some $v=\nabla^\perp\alpha\in\mathcal U$.
Then
$\eta_\epsilon(x)=x+\epsilon v(x)+\frac{\epsilon^2}2\nabla_vv(x)+\dots$
Taking second derivatives at $\epsilon=0$,
\begin{eqnarray}
\frac{d^2\omega_\epsilon(x)}{d\epsilon^2}_{|\epsilon=0}
&=&\frac{d^2}{d\epsilon^2}_{|\epsilon=0}\omega(x+\epsilon v(x)+\frac{\epsilon^2}2\nabla_vv(x)+\dots)\\
&=&\omega_{,k}v_lv_{k,l}+\omega_{,kl}v_kv_l
=(\omega_{,k}v_kv_l)_{,l}\\
&=&{\rm div}(\{\alpha,\omega\}v)\\
&=&{\rm div}(\nu v).
\end{eqnarray}
Integrating by parts, $-\int_\Omega \psi{\rm div}(\nu v)\text{dvol}=\int_\Omega \nu \nabla\psi\cdot v\text{dvol}$.
Taking gradients of $\omega=F(\psi)$,
\begin{equation}
-\int_\Omega \psi\ddot \omega\text{dvol}
=\int_\Omega \nu \frac{\nabla\omega\cdot v}{F'(\psi)}\text{dvol}=\int_\Omega\frac{\nu^2}{F'(\psi)}\text{dvol}.
\end{equation}
This completes the proof.
\stopproof%

\subsection{Notation}
Constants will generally be denoted by the same letter $C$, even in the derivation of tame estimates
where it is important that they depend on the regularity index $n$.
If $\mathcal K$ is the closure of a smooth, bounded region in Euclidean space
(e.g. $\mathcal K=\overline\Omega$ or $[0,|\Omega|]$), then  $C^n_{\mathcal K}$ is the space
of $n$-times continuously differentiable functions on $\mathcal K$.
The norm is given by
\begin{equation}\|f\|_{n}:=\|f\|_{C^n_{\mathcal K}}:=\sup_{0\leq j\leq n} \sup_{\mathcal K}|\nabla^jf|.\end{equation}
$C^{n,\alpha}_{\mathcal K}$ denotes the subspace of functions in $C^n_{\mathcal K}$ whose derivatives
up to order $n$
are H\"older continuous with exponent $\alpha$.  Throughout, $\alpha$ will be  a fixed constant in $(0,1)$.
The norm is denoted
\begin{equation}\|f\|_{n,\alpha}:=\|f\|_{C^{n,\alpha}_{\mathcal K}}:=\|f\|_{n}+\sum_{j=1}^n[\nabla^jf]_\alpha 
\qquad [f]_\alpha:=\sup_{x\neq y\in \mathcal K}\frac{|f(x)-f(y)|}{|x-y|^\alpha}.
\label{eq:Notation[]alpha}
\end{equation}

We introduce the following spaces of functions with various regularity:
\begin{eqnarray}
\mathcal F^n&:=&\{\omega\in C^n_{\overline\Omega}~|~\frac{\partial\omega}{\partial\tau}_{|\partial\Omega}=0\},\\
\mathcal U^n&:=&\{\psi\in \mathcal F^n~|~\psi_{|\Gamma_0}=0\},\\
\mathcal U^n_{\gamma_i}&:=&\{\psi\in \mathcal U^n~|~\int_{\Gamma_i}\frac{\partial\psi}{\partial {N}}=\gamma_i, 1\leq i\leq l\},\\
\mathcal F^{n,\alpha}&:=&\{\omega\in C^{n,\alpha}_{\overline\Omega}~|~\frac{\partial\omega}{\partial\tau}_{|\partial\Omega}=0\},\\
\mathcal U^{n,\alpha}&:=&\{\psi\in \mathcal F^{n,\alpha}~|~\psi_{|\Gamma_0}=0\},\\
\mathcal U^{n,\alpha}_{\gamma_i}&:=&\{\psi\in\mathcal U^{n,\alpha}~|~\int_{\Gamma_i}\frac{\partial\psi}{\partial {N}}=\gamma_i, 1\leq i\leq l\},\\
\mathcal F&:=&\{\omega\in C^\infty_{\overline\Omega}~|~\frac{\partial\omega}{\partial\tau}_{|\partial\Omega}=0\},\\
\mathcal U&:=&\{\psi\in \mathcal F~|~\psi_{|\Gamma_0}=0\},\\
\mathcal U_{\gamma_i}&:=&\{\psi\in\mathcal U~|~\int_{\Gamma_i}\frac{\partial\psi}{\partial {N}}=\gamma_i, 1\leq i\leq l\}.
\end{eqnarray}
When $l=1$, we will also work with the open subsets
\begin{eqnarray}
\mathcal F_+^n&:=&\{\omega\in\mathcal F^n~|~\nabla\omega\neq 0, \omega_{|\Gamma_i}<\omega_{|\Gamma_o}\},\\
\mathcal F_+^{n,\alpha}&:=&\{\omega\in\mathcal F^{n,\alpha}~|~\nabla\omega\neq 0, \omega_{|\Gamma_i}<\omega_{|\Gamma_o}\},\\
\mathcal F_+&:=&\{\omega\in\mathcal F~|~\nabla\omega\neq 0, \omega_{|\Gamma_i}<\omega_{|\Gamma_o}\}.
\end{eqnarray}
For simplicity, we will sometimes simply write $C^n$, $C^{n,\alpha}$, or $C^\infty$ for these spaces.
\section{The solution operator $\psi=S(F)$}\label{section:SolutionOperator}
A solution operator $F\mapsto \psi$ 
returning a uniquely defined steady-state solution can be constructed in a neighborhood of
the reference steady-state provided it satisfies the non-degeneracy condition (ND1).
For this part of the proof no restriction on the geometry of $\Omega$ is necessary, 
and thus we consider (in this Section only) a bounded domain $\Omega\subset \mathbb R^2$ with outer boundary component $\Gamma_0$
and an arbitrary number $l$ of inner boundary components, $\Gamma_i$, $1\leq i\leq l$.

We assume given a reference steady-state:
\begin{equation}\Delta\overline\psi=\overline F(\overline\psi),\qquad \overline\psi_{|\Gamma_o}=0,
\qquad \frac{\partial\overline\psi}{\partial \tau}_{|\Gamma_i}=0,\qquad \int_{\Gamma_i}\frac{\partial\overline\psi}{\partial N}=\gamma_i, \qquad 1\leq i\leq l
\label{eq:Delta(psibar)=Fbar(psibar)+BC}
\end{equation}
where the $\gamma_i$'s are fixed.
It is assumed to satisfy the non-degeneracy condition (ND1).
The goal of this section is to construct a solution operator $\psi=S(F)$ for $F$ in some neighborhood of $\overline F$
(along with the desired estimates)
returning a uniquely defined solution to the steady-state equation (\ref{eq:SteadyStateEquation}),
\begin{equation}\Delta\psi=F(\psi),\qquad \psi_{|\Gamma_o}=0,
\qquad \frac{\partial\psi}{\partial\tau}_{|\Gamma_i}=0,\qquad \int_{\Gamma_i}\frac{\partial\psi}{\partial N}=\gamma_i,\qquad 1\leq i\leq l.
\label{eq:Delta(psi)=F(psi)+BC}
\end{equation}
Note that we do not assume in this section that $F'$ does not vanish.
Recall also that the case $F'>0$ is special in that the corresponding solution automatically satisfies (ND1).

The boundary conditions in (\ref{eq:Delta(psi)=F(psi)+BC}) define the affine space
\begin{equation}
\mathcal U_{\gamma_i}
=\{\psi\in C^\infty_{\overline\Omega}~|~\psi_{|\Gamma_o}=0, \frac{\partial\psi}{\partial\tau}_{|\Gamma_i}=0, \int_{\Gamma_i}\frac{\partial\psi}{\partial N}dl=\gamma_i, 1\leq i\leq l\}.
\label{eq:AffineSpaceU(gamma(i))}
\end{equation}
The $\gamma_i$'s being fixed, its tangent space is 
\begin{equation}\mathcal U_0=\{\phi\in C^\infty_{\overline\Omega}~|~\phi_{|\Gamma_o}=0, \frac{\partial\phi}{\partial\tau}_{|\Gamma_i}=0, \int_{\Gamma_i}\frac{\partial\phi}{\partial N}dl=0, 1\leq i\leq l\}.
\label{eq:LinearSpaceU(0)}
\end{equation}
We will also consider the following linear equation in $\phi\in\mathcal U_0$  parametrized by $c\in C^\infty_{\overline\Omega}$:
\begin{equation}
\Delta \phi+c\phi=k,\qquad \phi_{|\Gamma_o}=0,\quad \frac{\partial\phi}{\partial\tau}_{|\Gamma_i}=0,\quad \int_{\Gamma_i}\frac{\partial\phi}{\partial N}=0,
\quad 1\leq i\leq l.
\end{equation}

\subsection{Linear estimates}
\begin{lemma}[Estimates for linear elliptic equations]\label{lemma:EstimatesForLinearEllipticEquations}
\hfill\smallskip%
\begin{enumerate}
\item Given $\omega\in C^\infty_{\overline\Omega}$, 
there exists a unique $\psi\in C^\infty_{\overline\Omega}$ solving
\begin{equation}\Delta\psi=\omega, \qquad \psi_{|\Gamma_0}=0,\qquad
\frac{\partial\psi}{\partial\tau}_{|\Gamma_i}=0,\qquad \int_{\Gamma_i}\frac{\partial\psi}{\partial {N}}=\gamma_i,\qquad 1\leq i\leq l.\end{equation}
It satisfies the tame estimates
\begin{equation}
\|\psi\|_{n+2,\alpha}\leq C(\|\omega\|_{n,\alpha}+\sum_{i=1}^l|\gamma_i|),
\qquad n\geq 0.
\label{ineq:HoelderEstimatesDelta(psi)=omega}
\end{equation}
\item For $c\in C^\infty_{\overline\Omega}$ in a $\|\cdot\|_{0,\alpha}$-neighborhood and any $\phi\in \mathcal U_0$ satisfying
\begin{equation}
\Delta \phi+c\phi=k,\qquad \phi_{|\Gamma_o}=0,\quad \frac{\partial\phi}{\partial\tau}_{|\Gamma_i}=0,\quad \int_{\Gamma_i}\frac{\partial\phi}{\partial N}=0,
\quad 1\leq i\leq l,
\end{equation}
we have for $n\geq 0$
\begin{equation}
\|\phi\|_{n+2,\alpha}\leq C\left(\|\Delta\phi+c\phi\|_{n,\alpha}+\|c\|_{n,\alpha}\|\phi\|_{0,\alpha}\right)
\label{ineq:EstimatesForDelta(phi)+cphi=k}
\end{equation}
where the constant depends on $n$, but not on $c$ nor $\phi$.
\end{enumerate}
\end{lemma}
\startproof 

\setcounter{counterpar}{0}
\nextpar{H\"older estimates on $\Delta\psi=\omega$}
The construction of $\psi$ from $\omega$ is standard, see \cite{MarchioroPulvirenti}.
In order to handle the boundary conditions of (\ref{eq:Delta(psi)=F(psi)+BC})
one defines
\begin{equation}u=\psi-\sum_{i=1}^l\psi_{|\Gamma_i}g_i\label{eq:u=psi-psi(Gammai)gi}
\end{equation}
where $g_i\in C^\infty_{\overline\Omega}$, $1\leq i\leq l$, are fixed functions with ${g_i}_{|\Gamma_0}=0$, ${g_i}_{|\Gamma_i}=1$,
and ${g_i}_{|\Gamma_j}=0$ for $j\neq i$.
Then $u$ satisfies 
\begin{equation}\Delta u = \omega-\psi_{|\Gamma_i}\Delta g_i,\qquad u_{|\partial \Omega}=0\end{equation}
for which  we have the Schauder estimates: for $n\geq 0$,
\begin{equation}
\|u\|_{n+2,\alpha}\leq C(\|\Delta u\|_{n,\alpha}+\|u\|_{0,\alpha}).
\end{equation}
Then, 
\begin{eqnarray}
\|\psi\|_{n+2,\alpha}
&\leq&C\cdot(\|\omega\|_{n,\alpha}+\|\psi\|_{0,\alpha}+\sum_{i=1}^l|\psi_{|\Gamma_i}|).
\end{eqnarray}
Similarly we have in Sobolev spaces
\begin{equation}\|\psi\|_{H^2}\leq C(\|\omega\|_{L^2}+\sum_{i=1}^l|\psi_{|\Gamma_i}|).\label{ineq:Gaarding}
\end{equation}
Observe that by the trace theorem, for each $i=1,\dots, l$ we have
\begin{equation}
|\psi_{|\Gamma_i}|\leq C\|\psi\|_{L^2(\partial\Omega)}\leq C\|\nabla\psi\|_{L^2(\Omega)}
\label{ineq:TraceTheorem}
\end{equation}
while an integration by parts gives
\begin{eqnarray}
\int_\Omega |\nabla\psi|^2&=&-\int_\Omega \psi\Delta\psi+\sum_{i=1}^l\psi_{|\Gamma_i}\gamma_i\\
&\leq&\epsilon^2\|\psi\|_{L^2}^2+\frac1{4\epsilon^2}\|\omega\|_{L^2}^2+\epsilon^2\sum_{i=1}^l\psi_{|\Gamma_i}^2+\frac1{4\epsilon^2}\sum_{i=1}^l\gamma_i^2
\end{eqnarray}
and therefore
\begin{equation}
\|\nabla\psi\|_{L^2}\leq \epsilon\|\psi\|_{L^2}+\frac C\epsilon\|\omega\|_{L^2}+\epsilon\sum_{i=1}^l|\psi_{|\Gamma_i}|+\frac C\epsilon\sum_{i=1}|\gamma_i|.
\end{equation}
Taking $\epsilon$ sufficiently small, we can achieve simultaneously
\begin{equation}
\sum_{i=1}^l|\psi_{|\Gamma_i}|\leq C\left(\epsilon\|\psi\|_{L^2}+\|\omega\|_{L^2}+\sum_{i=1}^l|\gamma_i|\right)
\end{equation}
and
\begin{equation}\|\psi\|_{H^2}\leq C\left(\|\omega\|_{L^2}+\sum_{i=1}^l|\gamma_i|\right).
\end{equation}

Finally, by Sobolev's embedding (dimension is $2$),
\begin{equation}\|\psi\|_{0,\alpha}\leq C\|\psi\|_{H^2}\leq C(\|\omega\|_{L^2}+\sum_{i=1}^l|\gamma_i|)
\leq C(\|\omega\|_{0,\alpha}+\sum_{i=1}^l|\gamma_i|).
\end{equation}
\stopstep

\nextpar{Estimates for $\Delta\phi+c\phi=k$}
For $c\in C^\infty_{\overline\Omega}$ and $\phi\in \mathcal U_0$ let $\Delta\phi+c\phi=k$.
Writing $\Delta\phi=k-c\phi$, from (\ref{ineq:HoelderEstimatesDelta(psi)=omega}),
and paying attention that the boundary conditions for $\phi\in \mathcal U_0$ are those of (\ref{eq:LinearSpaceU(0)}),
we deduce by (\ref{ineq:TameEstimatesForProductOfFunctions}) that
\begin{eqnarray}
\|\phi\|_{n+2,\alpha}&\leq&C\left(\|k\|_{n,\alpha}+\|c\|_{0,\alpha}\|\phi\|_{n,\alpha}+\|c\|_{n,\alpha}\|\phi\|_{0,\alpha}\right).
\end{eqnarray}
Now  restricting $c$ to a $\|\cdot\|_{0,\alpha}$-neighborhood
and using the interpolation $\|\phi\|_{n,\alpha}\leq \epsilon \|\phi\|_{n+2,\alpha}+C(\epsilon,n)\|\phi\|_{0,\alpha}$,
we can choose $\epsilon$ sufficiently small (and independent of $c$ in the $\|\cdot\|_{0,\alpha}$-neighborhood)
to get (\ref{ineq:EstimatesForDelta(phi)+cphi=k}).
\stopproof %

We say that $c\in C^\infty_{\overline\Omega}$ is \defn{non-degenerate} if
\begin{equation}
E(c)\phi=\Delta \phi+c\phi=0,\qquad \phi_{|\Gamma_o}=0,\quad \frac{\partial\phi}{\partial\tau}_{|\Gamma_i}=0,\quad \int_{\Gamma_i}\frac{\partial\phi}{\partial N}=0, \quad 1\leq i\leq l
\label{eq:Delta(phi)+cphi=k+BC}
\end{equation}
has only the trivial solution $\phi\equiv 0$.

\begin{proposition}\label{proposition:VE(c)kIsSmoothTame}
Suppose that $\overline c\in C^\infty_{\overline\Omega}$ is non-degenerate.
Then, there exists a $\|\cdot\|_{0,\alpha}$-neighborhood of $\overline c$,
\begin{equation}
\mathcal V_E(\overline c)
=\{c\in C^\infty_{\overline\Omega}~|~\|c-\overline c\|_{0,\alpha}<\epsilon_E\}
\label{eq:VE(overlinec)}
\end{equation}
such that
$E\colon \mathcal V_E(\overline c)\times \mathcal U_0\rightarrow C^\infty_{\overline\Omega}$
has a smooth tame family of inverses
\begin{eqnarray}
VE\colon
\left\{\begin{array}{ccc}\mathcal V_E(\overline c)&\times& C^\infty_{\overline\Omega}\\
c&&k
\end{array}
\longrightarrow
\begin{array}{c}\mathcal U_0\\\phi\end{array}\right\}.
\end{eqnarray}
For $n\geq 0$, $c\in\mathcal V_E(\overline c)$, and $k\in C^\infty_{\overline\Omega}$, 
\begin{equation}
\|\phi\|_{n+2,\alpha}\leq C \left(\|k\|_{n,\alpha}+\|c\|_{n,\alpha}\|k\|_{0,\alpha}\right).
\label{ineq:CnalphaTameEstimatesOnVE(c)k}
\end{equation}
The first derivative with respect to $c$ is given by
\begin{equation}DVE(c)\cdot(k,\chi)=VE(c)\cdot(-\chi\phi).
\label{eq:DVE(c)(k,chi}
\end{equation}
\end{proposition}
\startproof%

\setcounter{counterpar}{0}
\nextpar{Invertibility for $c$ in a $\|\cdot\|_{0,\alpha}$-neighborhood of $\overline c$}
Assume $\overline c$ is non-degenerate:
$\Delta+\overline c$, as an operator of Fr\'echet spaces with the boundary conditions (\ref{eq:Delta(phi)+cphi=k+BC}),
has trivial kernel.
Let $\phi\in C^{2,\alpha}_{\overline\Omega}$ solve 
\begin{equation}
\Delta \phi+\overline c\phi=0,
\qquad \phi_{|\Gamma_o}=0,\quad \frac{\partial\phi}{\partial\tau}_{|\Gamma_i}=0,\quad \int_{\Gamma_i}\frac{\partial\phi}{\partial N}=0, \quad 1\leq i\leq l.
\end{equation}
Then, $\phi\in C^{4,\alpha}$ and, repeating, $\phi\in C^\infty_{\overline\Omega}$.
Thus, $\phi$ is in the kernel of $\Delta+\overline c$ as an operator on Fr\'echet spaces,
and by non-degeneracy assumption $\phi\equiv 0$.
That is, $\Delta+\overline c$ has trivial kernel as an operator of Banach spaces $C^{2,\alpha}\rightarrow C^{0,\alpha}$.
The Fredholm alternative then implies that it is in fact an isomorphism of Banach spaces.
In particular, we have the estimate
\begin{equation}\|\phi\|_{2,\alpha}\leq C\|\Delta\phi+\overline c\phi\|_{0,\alpha}\end{equation}
for all $\phi\in C^{2,\alpha}_{\overline\Omega}$ satisfying the boundary conditions (\ref{eq:Delta(phi)+cphi=k+BC}),
and in particular when $\phi\in C^\infty_{\overline\Omega}$.

Let now $c\in C^\infty_{\overline\Omega}$:
\begin{eqnarray}
\|\phi\|_{0,\alpha}\leq \|\phi\|_{2,\alpha}&\leq&C\|\Delta\phi+\overline c\phi\|_{0,\alpha}\\
&\leq&C\left(\|\Delta\phi+c\phi\|_{0,\alpha}+\|(c-\overline c)\phi\|_{0,\alpha}\right)\\
&\leq&C\left(\|\Delta \phi+c\phi\|_{0,\alpha}+\|c-\overline c\|_{0,\alpha}\|\phi\|_{0,\alpha}\right).
\end{eqnarray}
Take then $\epsilon_E$ sufficiently small in (\ref{eq:VE(overlinec)}) so that the last term 
can be incorporated to the left-hand side
for all $\phi\in\mathcal U_0$ and any $c\in\mathcal V_E(\overline c)$:
\begin{equation}\|\phi\|_{0,\alpha}\leq C\|\Delta\phi+c\phi\|_{0,\alpha}.\end{equation}
But from (\ref{ineq:EstimatesForDelta(phi)+cphi=k}) we have
\begin{eqnarray}
\|\phi\|_{n+2,\alpha}
&\leq& C\left(\|\Delta\phi+c\phi\|_{n,\alpha}+\|c\|_{n,\alpha}\|\phi\|_{0,\alpha}\right)\\
&\leq&C\left(\|\Delta\phi+c\phi\|_{n,\alpha}+\|c\|_{n,\alpha}\|\Delta\phi+c\phi\|_{0,\alpha}\right)
\end{eqnarray}
as desired.
\stopstep%

\nextpar{Continuity in $c$ and $k$}
Let $c,\tilde c\in\mathcal V_E(\overline c)$ and $k,\tilde k\in C^\infty_{\overline\Omega}$,
and let $\Delta\phi+c\phi=k$, $\Delta\tilde \phi+\tilde c\tilde \phi=\tilde k$.
Then, $\phi-\tilde \phi$ solves
\begin{equation}(\Delta +\tilde c)\cdot (\phi-\tilde \phi)=(k-\tilde k)-(c-\tilde c)\phi\end{equation}
and the estimates from the previous paragraph and (\ref{ineq:TameEstimatesForProductOfFunctions}) give
\begin{eqnarray}
&&\|\phi-\tilde \phi\|_{n+2,\alpha}\\
&\leq& C\cdot \Big(\|k-\tilde k\|_{n,\alpha}+ \|(c-\tilde c)\phi\|_{n,\alpha}
+\|\tilde c\|_{n,\alpha}(\|k-\tilde k\|_{0,\alpha}+\|(c-\tilde c)\phi\|_{0,\alpha}\Big)\\
&\leq& C\left(\|k-\tilde k\|_{n,\alpha}
+\|c-\tilde c\|_{n,\alpha}\|\phi\|_{n,\alpha}\right)
\end{eqnarray}
where the constant depends on $\|\tilde c\|_{n,\alpha}$ only.
With $\tilde k$ and $\tilde c$ fixed, $\|\phi\|_{n,\alpha}$ remains bounded 
if $\|k\|_{n,\alpha}$ and $\|c\|_{n,\alpha}$ remain bounded.
Then, $\|\phi-\tilde\phi\|_{n+2,\alpha}$ can be made arbitrarily small provided $\|k-\tilde k\|_{n,\alpha}$
and $\|c-\tilde c\|_{n,\alpha}$ are taken sufficiently small.
That is,
\begin{equation}
\left\{\begin{array}{ccc}
\mathcal V_E(\overline c)^{n,\alpha}&\times&C^{n,\alpha}_{\overline\Omega}\\
c&&k
\end{array}
\longrightarrow
\begin{array}{c}C^{n+2,\alpha}_{\overline\Omega}\\\phi\end{array}
\right\}\qquad (n\geq 0)
\label{eq:ContinuityOfVE(c)kInC^{n,alpha}-Grading}
\end{equation}
is continuous as a map of Banach spaces.
Remark that this also implies that
\begin{equation}
\left\{\begin{array}{ccc}
\mathcal V_E(\overline c)^n&\times&C^n_{\overline\Omega}\\
c&&k
\end{array}
\longrightarrow
\begin{array}{c}C^{n+1}_{\overline\Omega}\\\phi\end{array}
\right\}\qquad (n\geq 0)
\label{eq:ContinuityOfVE(c)kInC^n-Grading}
\end{equation}
is continuous as a map of Banach spaces.
\stopstep%

\nextpar{First derivative of $VE$}
Since $VE(c)\cdot k$ is linear in $k$, its derivative in $k$ exists and is simply $VE(c)\cdot k$.
For the derivative in $c$, fix then $\chi\in C^\infty_{\overline\Omega}$ and let 
$c_t=c+t\chi$. 
Denote the solutions 
\begin{equation}\Delta \phi_t+c_t\phi_t=k,\qquad \Delta \phi+c\phi=k.\end{equation}
Then, 
\begin{equation}
\Delta\frac{\phi_t-\phi}t+c\frac{\phi_t-\phi}t=-\chi\phi_t,\qquad 
\text{or}\qquad
\frac{\phi_t-\phi}t=VE(c)\cdot (-\chi\phi_t).\end{equation}
By continuity of $VE$, the limit of $\frac{\phi_t-\phi}t$ exists in the $C^\infty$-topology.
Furthermore, it is given by \begin{equation}
DVE(c)\cdot (k,\chi)
=VE(c)\cdot (-\chi\phi).
\end{equation}
which is a continuous function of $c,k, \chi$ as a map of Fr\'echet spaces.
This shows that $VE$ is continuously differentiable as a map of Fr\'echet spaces.
Finally, it is clearly a tame map of $c,k, \chi$ since $VE$ is tame.
\stopstep%

\nextpar{$VE$ is smooth tame}
$E(c)\cdot\phi=\Delta\phi+c\phi$ is a smooth tame map as the sum of a linear differential operator
with constant coefficients
and multiplication of functions.
Since its inverse $VE(c)\cdot k$ is tame and continuously differentiable
with tame first derivative,
Theorem~5.3.1, p.~102, Part~I, and Theorem~3.1.1, p.~150, Part~II of \cite{HamiltonIFTNashMoser} 
imply that $VE(c)\cdot k$ is a smooth tame map.
\stopproof%

\subsection{The solution operator}\label{section:ConstructSolutionOperator}
We recall that $\Omega$ is not assumed to be diffeomorphic to an annulus, 
and that the reference steady-state $\overline F$ is assumed to satisfy (ND1),
i.e. $\overline c=-\overline F'(\overline\psi)$ is non-degenerate in the sense
that (\ref{eq:Delta(phi)+cphi=k+BC}) has only the trivial solution.

Let $I\subset \mathbb R$ be a closed, bounded interval strictly larger than $\text{range}(\overline\psi)$.

\begin{proposition}\label{proposition:SolutionMappsi=S(F)}
Let $\overline\omega=\overline F(\overline\psi)$ satisfy (ND1).
Then there exists a smooth tame solution operator to (\ref{eq:Delta(psi)=F(psi)+BC}):
\begin{equation}
S\colon \left\{\begin{array}{c}(\mathcal V_S(\overline F)\subset C^\infty_{I})\\F\end{array}
\longrightarrow
\begin{array}{c}C^\infty_{\overline\Omega}\\\psi
\end{array}\right\}
\end{equation}
where $\mathcal V_S(\overline F)$ is a $\|\cdot\|_{2,\alpha}$-neighborhood.
For $F\in\mathcal V_S(\overline F)$,
\begin{equation}
\|\psi\|_{n+2,\alpha}\leq C(\|F\|_{n,\alpha}+1), \qquad n\geq 0.
\label{ineq:CnalphaTameEstimatesForSolutionOperator}
\end{equation}
Letting $\mathcal V_S^n(\overline F)$ be the completion of $\mathcal V_S(\overline F)$ in $\|\cdot\|_n$,
it is continuous as a map
\begin{equation}
S\colon\left\{\begin{array}{c}
\mathcal V_S^n(\overline F)\\F\end{array}
\longrightarrow
\begin{array}{c}C^{n+1}_{\overline\Omega}\\
\psi
\end{array}\right\}\qquad (n\geq 1).
\end{equation}
The first derivative, given by
\begin{equation}\phi=DS(F)\cdot f=VE(-F'(\psi))\cdot (f\circ\psi)\label{eq:DVE(c)(k,chi)}\end{equation}
is continuous as a map
\begin{equation}
DS\colon\left\{\begin{array}{ccc}
\mathcal V_S^n(\overline F)&\times&C^{n-1}_I\\F&&f\end{array}
\longrightarrow
\begin{array}{c}C^{n}_{\overline\Omega}\\\phi
\end{array}\right\}\qquad (n\geq 1).
\end{equation}
For $F\in\mathcal V_S(\overline F)$, and any $f\in C^\infty_I$,
\begin{eqnarray}
\|\phi\|_{n+2,\alpha}&\leq& C\left(\|f\|_{n,\alpha}+\|F\|_{n+1,\alpha}\|f\|_{1,\alpha}\right),\qquad n\geq 1
\label{ineq:TameEstimatesOnphi=phi(F,f)Forn>=1}\\
\|\phi\|_{2,\alpha}&\leq& C\|f\|_{0,\alpha}.\label{ineq:BoundOn||phi||0,alpha}
\end{eqnarray}
The second derivative $\phi_{12}=D^2S(F)(f_1,f_2)$ is continuous as a map
\begin{equation}
D^2S\colon \left\{\begin{array}{ccccc}
\mathcal V_S^n(\overline F)&\times&C^{n-1}_I&\times&C^{n-1}_I\\F&&f_1&&f_2\end{array}
\longrightarrow
\begin{array}{c}C^{n-1}_{\overline\Omega}\\\phi_{12}
\end{array}\right\}\qquad (n\geq 1).
\end{equation}
It satisfies the tame estimates, for $F\in\mathcal V_S(\overline F)$ and any $f_1,f_2\in C^\infty_I$,
\begin{eqnarray}
\|\phi_{12}\|_{n+2,\alpha}
\leq 
C\left(\|f_1\|_{n+1,\alpha}\|f_2\|_{1,\alpha}+\|f_1\|_{1,\alpha}\|f_2\|_{n+1,\alpha}
+\|F\|_{n+2,\alpha}\|f_1\|_{2,\alpha}\|f_2\|_{2,\alpha}\right)
\end{eqnarray}
for $n\geq 1$ and
\begin{equation}
  \|\phi_{12}\|_{2,\alpha}\leq  C\|f_1\|_{1,\alpha}\|f_2\|_{1,\alpha}.
\end{equation}
\end{proposition}
\startproof%

\setcounter{counterpar}{0}
\nextpar{Construction of the solution operator}
The solution operator will be constructed by first applying the Implicit Function Theorem to suitable Banach spaces.
Since ultimately we want a solution operator in the smooth category, several choices for these spaces are possible.
We make the following choices.
Let 
$V_0(\overline\psi)\subset \{\psi\in H^2_\Omega~|~\psi_{|\Gamma_0}=0,\frac{\partial\psi}{\partial\tau}_{|\Gamma_i}=0, \int_{\Gamma_i}\frac{\partial \psi}{\partial N}dl=\gamma_i, 1\leq i\leq l\}$
a neighborhood of $\overline\psi$ in which $\text{range}(\psi)\subset I$.
(Note that this makes sense since $H^2_\Omega$ continuously embeds into $C^0_{\overline\Omega}$,
and that $\int_{\Gamma_i}\frac{\partial\psi}{\partial N}dl$ makes sense by the trace lemma.)
We may then define
\begin{equation}
H\colon\left\{\begin{array}{ccc}C^1_{I}&\times& V_0(\overline\psi)\\F&&\psi\end{array}
\longrightarrow
\begin{array}{c}L^2_\Omega\\\Delta\psi-F(\psi).
\end{array}\right\}.
\end{equation}
This operator is continuously Fr\'echet-differentiable (for emphasis we work in the Banach-category).
The assumption that $\overline\omega=\overline F(\overline \psi)$ satisfies (ND1) implies that
\begin{equation}D_\psi H(\overline F,\overline \psi)\phi=\Delta\phi-\overline F'(\overline\psi)\phi\end{equation}
has trivial kernel and hence is invertible as an operator
\begin{equation}\left\{\psi\in H^2_\Omega~|~\psi_{|\Gamma_0}=0,\frac{\partial\psi}{\partial\tau}_{|\Gamma_i}=0, \int_{\Gamma_i}\frac{\partial \psi}{\partial N}dl=\gamma_i\right\}\longrightarrow L^2_\Omega.
\end{equation}
The classical Implicit Function Theorem guarantees the existence of neighborhoods $V_S(\overline F)\subset C^1_{I}$
of $\overline F$ 
and $V_S(\overline \psi)\subset V_0(\overline\psi)$ of $\overline \psi$,
and a continuously Fr\'echet-differentiable map of Banach spaces 
$F\in V_S(\overline F)\mapsto \psi\in V_S(\overline \psi)$ such that $\Delta \psi=F(\psi)$.
(We will just say that this is a solution operator $C^1_{I}\rightarrow H^2_\Omega$.)
By Sobolev's embedding theorem, this is in fact a solution operator $C^1_{I}\rightarrow C^{0,\beta}_{\overline\Omega}$ for any $0<\beta<1$.
Now composition $(F,\psi)\mapsto F\circ \psi$ is continuous as an operator $C^1\times C^{0,\beta}\rightarrow C^{0,\alpha}$ for any $0<\alpha<\beta$
(see Lemma~\ref{lemma:CompositionIsSmoothTame}).
Thus, by elliptic regularity the solution operator is continuous $C^1_{I}\rightarrow C^{2,\alpha}_{\overline\Omega}$.
This along with a simple induction implies that the solution operator is continuous
\begin{equation}
\left\{\begin{array}{c}C^l_{I}\\F\end{array}
\longrightarrow
\begin{array}{c}C^{l+1}_{\overline\Omega}\\\psi\end{array}\right\}\qquad (l\geq 1).
\label{eq:ContinuityOfSolutionOperatorInCn-Grading}
\end{equation}
The above implies as well that the solution operator is continuous
\begin{equation}
\left\{\begin{array}{c}C^{l,\alpha}_{I}\\F\end{array}
\longrightarrow
\begin{array}{c}C^{l+1,\alpha}_{\overline \Omega}\\\psi\end{array}\right\}
\qquad (l\geq 1).
\end{equation}

A first requirement on $\mathcal V_S(\overline F)$ is that it be a $\|\cdot\|_{1,\alpha}$-neighborhood such that
\begin{equation}\mathcal V_S(\overline F)\quad\subset\quad \Big(V_S(\overline F)\cap C^\infty_{\overline\Omega}\Big).\end{equation}
\stopstep%

\nextpar{Tame estimates on $\psi=S(F)$ in the $C^{n,\alpha}$-grading}
From Lemma~\ref{lemma:EstimatesForLinearEllipticEquations}, 
$F\in\mathcal V_S(\overline F)$, by (\ref{ineq:TameEstimatesForCompositionOperator}) we have for $n\geq 1$
\begin{eqnarray}
\|\psi\|_{n+2,\alpha}
\leq C\cdot (\|F\circ\psi\|_{n,\alpha}+1)
\leq C\cdot (\|F\|_{n,\alpha}+\|\psi\|_{n,\alpha}\|F\|_{1,\alpha}),
\end{eqnarray}
and as long as $\psi$ remains in a $\|\cdot\|_{1,\alpha}$-neighborhood.
This is the case, e.g. by continuity $F\in C^{1,\alpha}\mapsto \psi\in C^{2,\alpha}$ and since $F\in\mathcal V_S(\overline F)$ remains
in a $\|\cdot\|_{1,\alpha}$-neighborhood.
Recall the inequality
$\|\psi\|_{n,\alpha}\leq \epsilon_1\|\psi\|_{n+1,\alpha}+C(\epsilon_1,n)\|\psi\|_{0,\alpha}$.
Choosing $\epsilon_1$ sufficiently small (depending on $n$ but independent of $F\in\mathcal V_S(\overline F)$)
we find 
\begin{eqnarray}
\|\psi\|_{n+2,\alpha}&\leq&C\cdot (\|F\|_{n,\alpha}+\|\psi\|_{0,\alpha}).
\end{eqnarray}
Since $\|\psi\|_{0,\alpha}$ is bounded (because $\|F\|_{1,\alpha}$ is) we obtain the desired tame estimates
(\ref{ineq:CnalphaTameEstimatesForSolutionOperator}) for $n\geq 1$.
Finally, since $\|F\|_{1,\alpha}$ remains bounded, so does $\|\psi\|_{0,\alpha}$ and the estimates in fact holds for $n\geq 0$
by increasing the constant if necessary.
\stopstep%

\nextpar{First derivative}
Let $F\in \mathcal V_S(\overline F)$ and $f\in C^\infty_I$,
and denote $\psi_t=S(F+tf)$ and $\psi=S(F)$ the corresponding solutions.
Then,
\begin{equation}
\Delta \frac{\psi_t-\psi}t-F'(\psi)\frac{\psi_t-\psi}t=\left(\frac{F(\psi_t)-F(\psi)}t-F'(\psi)\frac{\psi_t-\psi}t\right)+f(\psi_t).
\label{eq:Delta(psi(t)-psi)/t-F'(psi)(psi(t)-psi)/t}
\end{equation}
Now recall that the solution operator obtained from Step~1 is obtained via the classical Implicit Function Theorem,
and therefore is a Fr\'echet-differentiable map of Banach spaces $F\in C^1_I\mapsto \psi\in H^2_\Omega$ .
Thus, the limit $\phi$ of $\frac{\psi_t-\psi}t$ exists in $H^2$, and it satisfies
\begin{equation}\Delta\phi-F'(\psi)\phi=f(\psi)\qquad \text{in}\qquad L^2.\end{equation}
By Sobolev's embedding, $\frac{\psi_t-\psi}t\rightarrow_t \phi$ also in $C^{0,\alpha}$, and in turn
the right-hand side of (\ref{eq:Delta(psi(t)-psi)/t-F'(psi)(psi(t)-psi)/t}) converges in $C^{0,\alpha}$ as well.
This implies that $\frac{\psi_t-\psi}t\rightarrow_t \phi$ in $C^{2,\alpha}$ by continuity of $VE$.
Repeating, one can show that $\frac{\psi_t-\psi}t\rightarrow_t\phi$ in $C^{k,\alpha}$ for each $k$,
in other words that it converges in the $C^\infty$-topology.
This proves that the derivative $\phi$ of $\psi=S(F)$ at $F$ in the direction $f$ exists.
Furthermore, it is given by $\phi=VE(-F'\circ \psi)(f\circ\psi)$ which is clearly continuous and tame.
More precisely, from (\ref{eq:ContinuityOfVE(c)kInC^n-Grading}), we find that
\begin{equation}
\left\{\begin{array}{ccc}
\mathcal V_S^n(\overline F)&\times&C^{n-1}_I\\F&&f\end{array}
\longrightarrow
\begin{array}{c}C^n_{\overline\Omega}\\\phi=DS(F)f
\end{array}
\right\}\qquad (n\geq 1)
\label{eq:ContinuityOf(F,f)->phi=DS(F)f}
\end{equation}
is continuous as a map of Banach spaces.
(Taking $f\in C^n_I$ does not improve $\phi$.)
\stopstep%

\nextpar{Tame estimates on $\phi=DS(F)f$ in the $C^{n,\alpha}$-grading}
Since we will invoke tame estimates from Lemma~\ref{lemma:EstimatesForLinearEllipticEquations} 
to derive tame estimates on $\phi=DS(F)f$,
we will take $\epsilon_S$ sufficiently small so that 
\begin{equation}\mathcal V_S(\overline F)=\{F\colon \|F-\overline F\|_{2,\alpha}<\epsilon_S\}\subset (V_S(\overline F)\cap C^\infty_{\overline\Omega})
\end{equation}
and such that 
$F'(\psi)\in \mathcal V_E(\overline c)$ whenever $F\in\mathcal \mathcal V_S(\overline F)$, 
where $\overline c=-\overline F'(\overline \psi)$ and $\mathcal V_E(\overline c)$ 
is as in Lemma~\ref{proposition:VE(c)kIsSmoothTame}. 
By (\ref{ineq:CnalphaTameEstimatesOnVE(c)k}) and (\ref{ineq:TameEstimatesForCompositionOperator}), 
we have for  $F\in\mathcal V_S(\overline F)$ and $n\geq 1$
\begin{eqnarray}
\|\phi\|_{n+2,\alpha}
&\leq& C\left(\|f(\psi)\|_{n,\alpha}+\|f(\psi)\|_{0,\alpha}\|F'(\psi)\|_{n,\alpha}\right).
\end{eqnarray}
Now 
\begin{eqnarray}
\|f(\psi)\|_{n,\alpha}&\leq&C\left(\|f\|_{n,\alpha}+\|\psi\|_{n,\alpha}\|f\|_{1,\alpha}\right)\\
&\leq&C\left(\|f\|_{n,\alpha}+\|\psi\|_{n+3,\alpha}\|f\|_{1,\alpha}\right)\\
&\leq&C\left(\|f\|_{n,\alpha}+\|F\|_{n+1,\alpha}\|f\|_{1,\alpha}\right)
\end{eqnarray}
for $n\geq 1$ while
\begin{equation}\|f(\psi)\|_{0,\alpha}\leq C\|f\|_{0,\alpha}
\end{equation}
since $\|\psi\|_{1,\alpha}$ remains bounded (because $\|F\|_{1,\alpha}$ does).
Next, for $n\geq 1$
\begin{eqnarray}
\|F'(\psi)\|_{n,\alpha}&\leq&C\left(\|F'\|_{n,\alpha}+\|\psi\|_{n,\alpha}\|F'\|_{1,\alpha}\right)\\
&\leq&C\left(\|F\|_{n+1,\alpha}+\|\psi\|_{n+3,\alpha}\|F\|_{2,\alpha}\right)\\
&\leq&C\|F\|_{n+1,\alpha}
\end{eqnarray}
and this holds in fact for $n=0$ as well since $\|F'(\psi)\|_{0,\alpha}\leq C\|F\|_{1,\alpha}$
because $\|\psi\|_{1,\alpha}$ remains bounded.
With the above, we have for $n\geq 1$
\begin{eqnarray}
\|\phi\|_{n+2,\alpha}
&\leq&C\left(\|f\|_{n,\alpha}+\|F\|_{n+1,\alpha}\|f\|_{1,\alpha}\right)
\end{eqnarray}
and for $n=0$
\begin{eqnarray}
\|\phi\|_{2,\alpha}&\leq& C\left(\|f(\psi)\|_{0,\alpha}+\|F'(\psi)\|_{0,\alpha}\|f(\psi)\|_{0,\alpha}\right)\\
&\leq& C\left(\|f\|_{0,\alpha}+\|F\|_{1,\alpha}\|f\|_{0,\alpha}\right)\\
&\leq&C\|f\|_{0,\alpha}
\end{eqnarray}
since $\|\psi\|_{1,\alpha}$ and $\|F\|_{1,\alpha}$ remain bounded for $F\in\mathcal V_S(\overline F)$.
\stopstep%

\nextpar{Second derivative}
Since $\phi$ is a differentiable tame map of $F$ and $f$, it is immediate that $S(F)=\psi$
is twice continuously differentiable and tame (and in fact smooth tame).
Implicit differentiation on $\Delta\phi_1=F'(\psi)\phi_1+f_1(\psi)$ shows that the second
derivative $\phi_{12}$ in the directions $f_1, f_2$ satisfies
\begin{equation}
\Delta \phi_{12}=F'(\psi)\phi_{12}+F''(\psi)\phi_1\phi_2+f_2'(\psi)\phi_1+f_1'(\psi)\phi_2
\end{equation}
where $\Delta\phi_2=F'(\psi)\phi_2+f_2(\psi)$.
From (\ref{eq:ContinuityOf(F,f)->phi=DS(F)f}), we conclude that
\begin{equation}
D^2S\colon \left\{\begin{array}{ccccc}
\mathcal V_S^n(\overline F)&\times&C^{n-1}_I&\times&C^{n-1}_I\\F&&f_1&&f_2\end{array}
\longrightarrow
\begin{array}{c}C^{n-1}_{\overline\Omega}\\\phi_{12}
\end{array}\right\}\qquad (n\geq 1)
\end{equation}
is continuous as a map of Banach spaces.
\stopstep%

\nextpar{Tame estimates on the second derivatives}
We need first tame estimates on $F''(\psi)\phi_1\phi_2+f_2'(\psi)\phi_1+f_1'(\psi)\phi_2$.
Since $F\in\mathcal V_S(\overline F)$, we have for $n\geq 1$
\begin{eqnarray}
\|F''(\psi)\|_{n,\alpha}&\leq&C\left(\|F''\|_{n,\alpha}+\|\psi\|_{n,\alpha}\|F''\|_{1,\alpha}\right)\\
&\leq&C\left(\|F\|_{n+2,\alpha}+\|F\|_{3,\alpha}\right)\\
&\leq&C\|F\|_{n+2,\alpha}
\end{eqnarray}
and this holds in fact also for $n=0$ since $\|\psi\|_{1,\alpha}$ remains bounded.
Next with $i=1,2$,
\begin{equation}
\|f_i'(\psi)\|_{n,\alpha}\leq C\left(\|f_i\|_{n+1,\alpha}+\|F\|_{n,\alpha}\|f_i\|_{2,\alpha}\right)
\end{equation}
for $n\geq 1$ while $\|f_i'(\psi)\|_{0,\alpha}\leq C\|f_i\|_{1,\alpha}$
since $\|\psi\|_{1,\alpha}$ remains bounded.
Thus, recalling (\ref{ineq:TameEstimatesOnphi=phi(F,f)Forn>=1}) and (\ref{ineq:BoundOn||phi||0,alpha}),
we have for $n\geq 3$
\begin{eqnarray}
&&\|F''(\psi)\phi_1\phi_2+f_2'(\psi)\phi_1+f_1'(\psi)\phi_2\|_{n,\alpha}\\
&\leq&C\Big(\|F''(\psi)\|_{n,\alpha}\|\phi_1\|_{0,\alpha}\|\phi_2\|_{0,\alpha}
+\|\phi_1\|_{n,\alpha}\|\phi_2\|_{0,\alpha}+\|\phi_1\|_{0,\alpha}\|\phi_2\|_{n,\alpha}\\
&&+\|f_1'(\psi)\|_{n,\alpha}\|\phi_2\|_{0,\alpha}+\|f_1'(\psi)\|_{0,\alpha}\|\phi_2\|_{n,\alpha}\\
&&+\|f_2'(\psi)\|_{n,\alpha}\|\phi_1\|_{0,\alpha}+\|f_2'(\psi)\|_{0,\alpha}\|\phi_1\|_{n,\alpha}\Big)\\
&\leq&C\Big(
\|F\|_{n+2,\alpha}\|f_1\|_{0,\alpha}\|f_2\|_{0,\alpha}\\
&&+(\|f_1\|_{n-2,\alpha}+\|F\|_{n-1,\alpha}\|f_1\|_{1,\alpha})\|f_2\|_{0,\alpha}\\
&&+\|f_2\|_{0,\alpha}(\|f_2\|_{n-2,\alpha}+\|F\|_{n-1,\alpha}\|f_2\|_{1,\alpha})\\
&&+(\|f_1\|_{n+1,\alpha}+\|F\|_{n,\alpha}\|f_1\|_{2,\alpha})\|f_2\|_{0,\alpha}\\
&&+\|f_1\|_{1,\alpha}(\|f_2\|_{n-2,\alpha}+\|F\|_{n-1,\alpha}\|f_2\|_{1,\alpha})\\
&&+(\|f_2\|_{n+1,\alpha}+\|F\|_{n,\alpha}\|f_2\|_{2,\alpha})\|f_1\|_{0,\alpha}\\
&&+\|f_2\|_{1,\alpha}(\|f_1\|_{n-2,\alpha}+\|F\|_{n-1,\alpha}\|f_1\|_{1,\alpha})
\Big)\\
&\leq&C\Big(\|f_1\|_{n+1,\alpha}\|f_2\|_{1,\alpha}+\|f_1\|_{1,\alpha}\|f_2\|_{n+1,\alpha}+\|F\|_{n+2,\alpha}\|f_1\|_{2,\alpha}\|f_2\|_{2,\alpha}
\Big)
\end{eqnarray}
and one can verify that this also holds for $n=1, 2$ since $\|\psi\|_{4,\alpha}$ remains bounded for $F\in\mathcal V_S(\overline F)$.
Finally, for $n=0$, we have
\begin{equation}\|F''(\psi)\phi_1\phi_2+f_2'(\psi)\phi_1+f_1'(\psi)\phi_2\|_{n,\alpha}\leq C\|f_1\|_{1,\alpha}\|f_2\|_{1,\alpha}.\end{equation}

The tame estimates on $\phi_{12}$ are thus, for $n\geq 1$,
\begin{eqnarray}
&&\|\phi_{12}\|_{n+2,\alpha}\\
&\leq& C\Big(\|f_1\|_{n+1,\alpha}\|f_2\|_{1,\alpha}+\|f_1\|_{1,\alpha}\|f_2\|_{n+1,\alpha}
+\|F\|_{n+2,\alpha}\|f_1\|_{2,\alpha}\|f_2\|_{2,\alpha}\\
&&+\|F\|_{n+1,\alpha}\left(\|f_1\|_{2,\alpha}\|f_2\|_{2,\alpha}+\|F\|_{3,\alpha}\|f_1\|_{2,\alpha}\|f_2\|_{2,\alpha}\right)
\Big)\\
&\leq&C\Big(\|f_1\|_{n+1,\alpha}\|f_2\|_{1,\alpha}+\|f_1\|_{1,\alpha}\|f_2\|_{n+1,\alpha}
+\|F\|_{n+2,\alpha}\|f_1\|_{2,\alpha}\|f_2\|_{2,\alpha}\Big)
\end{eqnarray}
where we have used interpolation to get $\|F\|_{n+1,\alpha}\|F\|_{3,\alpha}\leq C\|F\|_{n+2,\alpha}\|F\|_{2,\alpha}\leq C\|F\|_{n+2,\alpha}$.
Finally for $n=0$ we have
\begin{equation}
\|\phi_{12}\|_{2,\alpha}\leq C\|f_1\|_{1,\alpha}\|f_2\|_{1,\alpha}.\end{equation}
\stopproof %
\setcounter{counterpar}{0}

\noindent\textbf{Remarks}\qquad
G\aa rding's inequality fails in the $C^n$-grading.
Yet, $\omega=F(\psi)$ is $C^n$ if $F$ is $C^n$, 
and furthermore (see Lemma~\ref{lemma:CompositionIsSmoothTame}) the map
\begin{equation}
\left\{\begin{array}{c}\mathcal V_S^n(\overline F)\\F\end{array}
\longrightarrow
\begin{array}{c}C^n_{\overline\Omega}\\\omega=F(\psi)\end{array}\right\}\qquad (n\geq 1)
\label{eq:ContinuityOfF->omegaInC^n-Grading}
\end{equation}
is continuous.
Likewise, from $\nu=\Delta \phi$ it would appear that $\nu$ loses two derivatives from $F$.
However, writing $\nu=F'(\psi)\phi+f(\psi)$, one finds that in fact
\begin{equation}
\left\{\begin{array}{ccc}
\mathcal V_S^n(\overline F)&\times&C_I^{n-1}\\F&&f\end{array}
\longrightarrow
\begin{array}{c}C^{n-1}_{\overline\Omega}\\
\nu
\end{array}\right\}\qquad (n\geq 1)
\label{eq:ContinuityOf(F,f)->nuInC^n-Grading}
\end{equation}
is continuous as a map of Banach spaces.
Finally, the second derivative $\nu_{12}$ of $F\mapsto \omega$ in the directions $(f_1, f_2)$ 
is continuous
\begin{equation}
\left\{\begin{array}{ccccc}
\mathcal V_S^n(\overline F)&\times&C^{n-1}_I&\times&C^{n-1}_I\\F&&f_1&&f_2\end{array}
\longrightarrow
\begin{array}{c}C^{n-2}_{\overline\Omega}\\\nu_{12}
\end{array}\right\}\qquad (n\geq 1)
\label{eq:ContinuityOf(F,f1,f2)->nu12InC^n-Grading}
\end{equation}
as a map of Banach spaces.
\stopproof%

\section{Distribution functions and co-adjoint orbits}\label{section:DistributionFunctions}
We are now assuming that $\Omega$ is diffeomorphic to an annulus:
\begin{equation}\partial \Omega =\Gamma_o\cup \Gamma_i.\end{equation}
This Section is concerned with establishing properties and estimates on distribution functions $A_\omega$
for functions $\omega$, which are locally constant on the boundary and have no critical points.
These results will apply equally to stream functions, which further satisfy $\psi_{|\Gamma_o}=0$,
and to steady-state vorticity functions, which are locally constant on $\partial\Omega$ as observed
in the Introduction.
Without loss of generality we will assume that $\omega_{|\Gamma_i}<\omega_{|\Gamma_o}$,
and therefore we will work with spaces $\mathcal F_+^n$ and $\mathcal F_+$ introduced at the end of Section~\ref{section:Introduction}.

We begin with the following useful result.
\begin{lemma}[Global coordinates on $\overline\Omega$ induced by $\omega$]\label{lemma:CoordinatesDependingOnomega}
For each $\omega\in\mathcal F_+$ there exists a global coordinate system for $\overline\Omega$,
$z\colon [0,1]\times\mathbb S^1\rightarrow \overline\Omega$,
such that $\{t\}\times \mathbb S^1$ is mapped onto the level set 
$\{x\in\overline\Omega~\colon~\omega(x)=\min\omega+t(\max\omega-\min\omega)\}$:
\begin{equation}\omega(z(t,s))=\min\omega+t(\max\omega-\min\omega).
\label{eq:omega(z)=t}
\end{equation}
The map $\omega\mapsto z$ is continuous $\mathcal F_+^{n}\rightarrow C^{n-1}_{[0,1]\times \mathbb S^1}$,
$n\geq 2$.
\end{lemma}
\noindent\textbf{Remark}\qquad
The proof of Lemma~\ref{lemma:CoordinatesDependingOnomega} can easily be adapted to achieve a continuous map
$\mathcal F_+^n\rightarrow C^n_{[0,1]\times \mathbb S^1}$.\footnote{The loss of derivative
is in the $s$-direction only.
To regain this derivative, write the inverse mapping as
$z^{-1}=\left(\frac{\omega-\min\omega}{\max\omega-\min\omega}, \theta\right)$,
apply a suitable smoothing operator to the second factor $\theta$, and note that this perserves the property (\ref{eq:omega(z)=t}).}
However, this will make no difference in the rest of the paper
as $z$ will always be used in connection with a factor $\frac1{|\nabla\omega|}$ which is of class $C^{n-1}$
when $\omega\in \mathcal F_+^n$.
\bigskip%

\startproof%
Let $c_i(s)$, $0\leq s\leq 1$, be a smooth parametrization of $\Gamma_i$.
For each $s$, let $z=z(t,s)$ be the solution to
\begin{equation}\dot z=(\max\omega-\min\omega)\frac{\nabla\omega(z)}{|\nabla\omega(z)|^2},\qquad z(0)=c_i(s).
\label{eq:ODEForGlobalCoordinatesOnOmega}
\end{equation}
Since $\frac{d}{dt}\omega(z(t))=\max\omega-\min\omega$, we have $\omega(z(t,s))=\min\omega+t(\max\omega-\min\omega)$
for $t\in[0,1]$ and $s\in\mathbb S^1$.
It is standard that if $\omega\in C^n$, $n\geq 2$, then $z$ is $C^n$ in $t$ and $C^{n-1}$ in $s$.

For $\omega,\omega_1$, denote $f=(\max\omega-\min\omega)\frac{\nabla\omega}{|\nabla\omega|^2}$, 
$f_1=(\max\omega_1-\min\omega_1)\frac{\nabla\omega_1}{|\nabla \omega_1|^2}$,
and $z,z_1$ the corresponding coordinate systems.
For $t\in[0,1]$, $s\in\mathbb S^1$,
\begin{eqnarray}
&&z(t,s)-z_1(t,s)\\
&=&\int_0^{t}\left(f(z(\theta,s))-f(z_1(\theta,s))\right)d\theta+\int_0^t(f-f_1)(z_1(\theta,s))d\theta\\
&=&\int_0^{t}\int_0^1\nabla f(z_1(\theta,s)+\tau(z(\theta,s)-z_1(\theta,s))\cdot (z(\theta,s)-z_1(\theta,s)d\tau d\theta\\
&&+\int_0^t(f-f_1)(z_1(\theta,s))d\theta
\end{eqnarray}
so that
\begin{eqnarray}
|z(t,s)-z_1(t,s)|
&\leq&\sup|\nabla f|\int_0^t|z(\theta,s)-z_1(\theta,s)|d\theta+\sup |f-f_1|.
\end{eqnarray}
This is the same as
\begin{equation}
\frac{d}{dt}\left(e^{-\sup|\nabla f|t}\int_0^t|z-z_1|(\theta,s)d\theta\right)\leq e^{-\sup|\nabla f|t}\sup|f-f_1|
\end{equation}
hence integrating
\begin{equation}
\int_0^t|z-z_1|(\theta,s)d\theta\leq e^{\sup|\nabla f|t}\int_0^te^{-\sup|\nabla f|\theta}\sup|f-f_1|d\theta
\leq \sup |f-f_1|\frac{e^{|\sup\nabla f|}}{|\sup\nabla f|}.
\end{equation}
In turn,
\begin{equation}
|z(t,s)-z_1(t,s)|\leq \sup |f-f_1|\left(e^{\sup|\nabla f|}+1\right).
\end{equation}
Given $f_1\in C^1$, for any $f\in C^1$ in a neighborhood of $f_1$ such that $|\nabla f|\leq M$,
then $\sup|z-z_1|$ can be made arbitrarily small provided $\sup|f-f_1|$ is chosen sufficiently small.
In other words, $f\mapsto z$ is continuous as a map $C^1\rightarrow C^0$.

Estimate now
\begin{eqnarray}
&&|\dot z(t,s)-\dot z_1(t,s)|\\
&\leq&|f(z(t,s))-f(z_1(t,s))|+|(f-f_1)(z(t,s))|\\
&\leq&\int_0^1|\nabla f(z_1(t,s)+\tau(z(t,s)-z_1(t,s)))||z(t,s)-z_1(t,s)|d\tau+\sup|f-f_1|\\
&\leq&\sup|\nabla f||z-z_1|+\sup|f-f_1|
\end{eqnarray}
which shows that $f\mapsto \dot z$ is continuous $C^1\rightarrow C^0$.

To estimate the derivative $z'$ in $s$, differentiate (\ref{eq:ODEForGlobalCoordinatesOnOmega}):
\begin{eqnarray}
\frac{d}{dt}z'(t,s)=\nabla f(z(t,s))\cdot z'(t,s),\qquad 
\frac{d}{dt}z_1'(t,s)=\nabla f_1(z_1(t,s))\cdot z_1'(t,s)
\end{eqnarray}
and observe that $z'(0,s)=z_1'(0,s)$, $s\in\mathbb S^1$.
Integrating, 
\begin{eqnarray}
z'(t,s)-z_1'(t,s)
&=&\int_0^t\nabla f(z_1(\theta,s))\left(z'(\theta,s)-z_1'(\theta,s)\right)d\theta\\
&&+\int_0^t \left(\nabla f(z(\theta,s))-\nabla f_1(z_1(\theta,s))\right)z_1'(\theta,s)d\theta
\end{eqnarray}
and thus
\begin{eqnarray}
|z'(t,s)-z_1'(t,s)|
\leq\sup|\nabla f|\int_0^t|z'(\theta,s)-z_1'(\theta,s)|d\theta +\sup |z_1'|\sup |\nabla (f-f_1)|.
\end{eqnarray}
As above, we find
\begin{equation}
|z'(t,s)-z_1'(t,s)|\leq \sup|z_1'|\sup |\nabla (f-f_1)|\left(\frac{e^{\sup |\nabla f|}}{\sup|\nabla f|}+1\right),
\end{equation}
showing that $f\mapsto z'$ is continuous $C^1\rightarrow C^0$.

Therefore, $f\mapsto z$ is continuous $C^1\rightarrow C^1$.
Differentiating further in $t$ and $s$, one establishes easily by induction that 
$f\mapsto z$ is continuous $C^k\rightarrow C^k$ for each $k\geq 1$.
\stopproof%

\subsection{Properties of the distribution function $A_\omega$}\label{section:PropertiesOfAomega}
Linearizing $\omega\circ\eta$ at $\eta=\text{id}$,
a tangent $\nu$ to $\mathcal O_\omega$ at $\omega$ is of the form
\begin{equation}\nu=\{\omega,\alpha\}\qquad \text{for some stream function}~\alpha\in\mathcal U.
\label{eq:nu={omega,alpha}}
\end{equation}
This is a first order PDE in $\alpha$ hence can be locally integrated along the characteristics,
which are here closed curves.
The compatibility conditions (to be able to integrate round these closed curves)
are
\begin{equation}\int\limits_{\{x\colon\omega(x)=\lambda\}}\frac{\nu(x)}{|\nabla\omega(x)|}dl(x)=0,\qquad \min\omega\leq\lambda \leq\max\omega.
\label{eq:CharacterizationOfnu={alpha,omega}}\end{equation}
On the other hand, if $\omega\in\mathcal F_+$,
then the distribution function $A_\omega(\lambda)= |\left\{x\colon\omega(x)<\lambda\right\}|$ can be expressed as
\begin{equation}
A_\omega(\lambda) = \int_{\min\omega}^\lambda\int_{x:\omega(x)=\lambda}\frac1{|\nabla\omega(x)|}dl(x)d\lambda,
\qquad \lambda\in [\min\omega,\max\omega].
\label{eq:Apsi(lambda)}
\end{equation}
This is easily seen using the coarea formula (see \S~3.2 in \cite{Federer})
\begin{equation}
\int_\Omega u(x)|\nabla\omega(x)|\zeta(\omega(x))~dx 
= \int_{\min\omega}^{\max\omega} \zeta(\lambda')\int_{x:\omega(x)=\lambda'}u(x)dl(x)~d\lambda'
\label{eq:BasicCoareaFormula}
\end{equation}
with
$u(x)=\frac1{|\nabla\omega(x)|}$ and $\zeta(\lambda')$ the characteristic function
over the interval $[\min\omega,\lambda]$.
Much of the work will be devoted to the detailed study of
\begin{equation}
J_\omega u(\lambda)
=\int_{x:\omega(x)=\lambda}u(x)dl(x),
\qquad \lambda \in [\min \omega,\max\omega].
\label{eq:Jpsiu(lambda)}
\end{equation}
With this notation, the coarea formula reads
\begin{equation}
\int_\Omega u(x) |\nabla\omega(x)|\zeta(\omega(x))dx = \int_{\min\omega}^{\max\omega} \zeta(\lambda) J_\omega u(\lambda)~d\lambda.
\label{eq:CoareaFormula}
\end{equation}
The derivative of $A_\omega$ can then be written as
\begin{equation}A_\omega'=J_\omega\frac1{|\nabla\omega|}\label{eq:DerivativeOfApsi}
\end{equation}
and, differentiating $A_\omega^{-1}(A_\omega(\lambda))=\lambda$,  that of $A_\omega^{-1}$ is
\begin{equation}(A_\omega^{-1})'=\frac1{J_\omega\frac1{|\nabla\omega|}\circ A_\omega^{-1}}.
\label{eq:DerivativeOfApsiInverse}
\end{equation}

\begin{proposition}[$A_\omega$ characterizes $\mathcal O(\omega)$ locally]\label{proposition:AomegaCharacterizesO*(omega)}
\hfill\smallskip%
\begin{enumerate}
\item 
For $\omega,\omega_1\in\mathcal F_+$, let $z,z_1$ be the global coordinate transformations from Lemma~\ref{lemma:CoordinatesDependingOnomega}.
If $\omega,\omega_1\in\mathcal F_+$ and $A_{\omega_1}=A_\omega$, 
then there exists $\eta\in\mathcal D_\text{vol}$ such that $\omega=\omega_1\circ\eta$.
\item Let $\omega\in \mathcal F_+$ and let $\nu\in\mathcal F$ such that $\nu_{|\partial\Omega}=0$.
Then, 
there exists a stream function $\alpha\in\mathcal U$ such that 
\begin{equation}\nu=\{\omega,\alpha\}\label{eq:nu={alpha,omega}SecondTime}
\end{equation}
if and only  if
\begin{equation}\int_{\omega=\lambda}\frac\nu{|\nabla\omega|}dl=0,\qquad \lambda\in\text{range}(\omega).
\label{eq:FirstOrderPDECompatibility}
\end{equation}
\item Let $\omega_\epsilon\in\mathcal F_+$, $\omega_0=\omega$.
Then, 
$\frac{dA_{\omega_\epsilon}}{d\epsilon}_{|\epsilon=0}=0$,
if and only if $\nu=\frac{d\omega_\epsilon}{d\epsilon}_{|\epsilon=0}$ is tangent to the orbit $\mathcal O(\omega)$ at $\omega$,
i.e. $\nu=\{\omega,\alpha\}$ for some $\alpha\in\mathcal U$.
\end{enumerate}
\end{proposition}
\startproof%

\setcounter{counterpar}{0}
\nextpar{$A_\omega$ characterizes $\mathcal O(\omega)$ locally}
If $A_{\omega_1}=A_\omega$, then setting $\phi=z_1\circ z^{-1}$, we have $\omega=\omega_1\circ\phi$.
Next we construct a diffeomorphism $\psi$ which moves points along the level sets of $\omega_1$ 
in such a way that $\eta=\psi\circ \phi\in\mathcal D_\text{vol}$.
We write it in the form
\begin{equation}\psi(x)=z_1(t, a(t,s)),\qquad \text{where}\qquad x=z_1(t,s)
\end{equation}
or $\psi=z_1\circ \alpha$, where $\alpha(t,s)=(t,a(t,s))$.
The condition $\text{det}\left(\frac{\partial (\psi\circ\phi)}{\partial(x,y)}\right)=1$
can be written in the form $(Z\circ \alpha)\partial_s a=F$ where $Z=\text{det}\left(\frac{\partial z_1}{\partial(t,s)}\right)$
and $F=\frac1{\text{det}\left(\frac{\partial\phi}{\partial (x,y)}\right)\circ\phi^{-1}}$.
This is a collection of ODEs in $s$ parametrized by $t$.
Let then $a(t,s)$, $(t,s)\in [0,1]\times [0,1]$, denote the solution such that $a(t,0)=0$, $0\leq t\leq 1$.
Setting $\eta=\psi\circ \phi$ defines a \foremph{local} diffeomorphism $\overline\Omega\rightarrow \overline\Omega$
such that $\text{det}\left(\frac{\partial\eta}{\partial (x,y)}\right)=1$.
One needs to check that $\eta$ is a global diffeomorphism of $\overline \Omega$.

Fix $\lambda\in\text{range}(\omega)$ so that $z\left(\frac{\lambda-\min\omega}{\max\omega-\min\omega},0\right)$ is a point
on the level set $\{\omega =\lambda\}$, and denote $q=q(t)$ the solution to
\begin{equation}\dot q=\nabla^\perp \omega(q),\qquad q(0)=z\left(\frac{\lambda-\min\omega}{\max\omega-\min\omega},0\right).\end{equation}
Clearly, $q(t)$ travels around the level set $\{\omega=\lambda\}$.
This is a Hamiltonian flow with respect to the symplectic form $\text{dvol}=dx\wedge dy$.
Now $\eta$ is a local symplectomorphism, so that $q_1=\eta\circ q$ solves
\begin{equation}
\dot q_1 =\nabla^\perp\omega_1(q_1),\qquad q_1(0)=\eta\left(z\left(\frac{\lambda-\min\omega}{\max\omega-\min\omega},0\right)\right).
\end{equation}
Clearly, $q_1(t)$ travels around the level set $\{\omega_1=\lambda\}$.
But $A_\omega=A_{\omega_1}$, so that the travel time of $q(t)$ and $q_1(t)$ around 
$\{\omega=\lambda\}$ and $\{\omega_1=\lambda\}$ respectively are the same:
\begin{equation}
\int_{\omega_1=\lambda}\frac{dl}{|\nabla\omega_1|}=\frac{d}{d\lambda}A_{\omega_1}(\lambda)=\frac{d}{d\lambda}A_\omega(\lambda)
=\int_{\omega=\lambda}\frac{dl}{|\nabla\omega|}.\end{equation}
Finally, $q_1=\eta\circ q$ so $\eta$ takes $\{\omega=\lambda\}$ onto $\{\omega_1=\lambda\}$ for each $\lambda$ and the claim is proved.
\stopstep%

\nextpar{Characterizing $\nu\parallel \mathcal O(\omega)$}
Equation (\ref{eq:nu={alpha,omega}SecondTime}) is a first order PDE in $\alpha$,
and the characteristics are the level sets of $\omega$.
Thus, it is locally solvable and (\ref{eq:FirstOrderPDECompatibility}) are precisely the compatibility conditions
that ensure that $\alpha$ is globally defined
(the characteristics are simple closed curves).
\stopstep%

\nextpar{Characterizing $\nu\parallel \mathcal O(\omega)$ in terms of distribution functions}
We will show that (\ref{eq:FirstOrderPDECompatibility}) holds.  This is immediate once the identity
\begin{equation}
\frac{\partial}{\partial\epsilon}_{|\epsilon=0}A_{\omega_\epsilon}(\lambda)
=-J_\omega\frac\nu{|\nabla\omega|}(\lambda)=-\int_{\omega=\lambda}\frac\nu{|\nabla\omega|}dl
\label{eq:ddepsilonAomegaepsilon}
\end{equation}
is proved.
Let $\lambda$ be in the interior of $\text{range}(\omega)$ and let $f$ be a $C^1$ function of $\mathbb R$
such that $f'$ has support contained in the interior of $\text{range}(\omega)$.
Then, by the coarea formula (\ref{eq:CoareaFormula})
\begin{equation}
\frac{\partial}{\partial}_{|\epsilon=0}\int_\Omega f(\omega_\epsilon(x))dx
=\int_\Omega f'(\omega(x))\nu(x)dx=\int_{\mathbb R}f'(\lambda)J_\omega\frac\nu{|\nabla\omega|}(\lambda)d\lambda.
\end{equation}
By approximation, this holds for $f$ a continuous, piecewise linear function.
Fix then $\lambda$ and for $\delta>0$ let $f^\delta(\lambda')$ have value $1$ for $\lambda'<\lambda$,
$0$ for $\lambda'>\lambda+\delta$, and linear in between.
Then
$\frac{\partial}{\partial\epsilon}_{|\epsilon=0}\int_\Omega f^\delta(\omega_\epsilon(x))dx\rightarrow_\delta
\frac{\partial}{|\partial\epsilon}_{|\epsilon=0}A_{\omega_\epsilon}(\lambda)$
while 
$\int_{\mathbb R}(f^\delta)'(\lambda')J_\omega\frac\nu{|\nabla\omega|}(\lambda')d\lambda'\rightarrow_\delta -(J_\omega u)(\lambda).$
\stopproof%

\begin{lemma}\label{lemma:ddlambdaJpsiu(lambda)}
Setting ${N}=\frac{\nabla \omega}{|\nabla\omega|}$,
\begin{equation}
\frac{d}{d\lambda}(J_\omega u)(\lambda)
=J_\omega \left(\frac{{\rm div}(u{N})}{|\nabla\omega|}\right)(\lambda),\qquad \lambda\in[\min\omega,\max\omega]
\label{eq:ddlambdaJpsiu}
\end{equation}
\end{lemma}
\startproof%
Let $\zeta$ be an arbitrary smooth function with compact support in $(\min\omega,\max\omega)$.
Integrating by parts the coarea formula (\ref{eq:CoareaFormula}),
using the identities $\nabla (\zeta\circ \omega)=\zeta'(\omega)\nabla \omega$ and
${\rm div}((\zeta\circ \omega)u{N})=(\zeta\circ \omega){\rm div}(u{N})+\langle \nabla(\zeta\circ\omega),u{N}\rangle$,
and the coarea formula again, we have
\begin{eqnarray}
\int_{\min\omega}^{\max\omega} \zeta(\lambda)(J_\omega u)'(\lambda)~d\lambda 
&=& -\int_\Omega u(x)|\nabla\omega(x)|\zeta'(\omega(x))~dx\\
&=& \int_\Omega {\rm div}(u{N})(\zeta\circ\omega)~dx\\
&=&\int_{\min\omega}^{\max\omega} \zeta(\lambda)J_\omega \left(\frac{{\rm div}~(u{N})}{|\nabla\omega|}\right)(\lambda)~d\lambda.
\end{eqnarray}
\stopproof%

\begin{lemma}\label{lemma:DerivativeOfJpsiepsilon(lambda)Inepsilon}
Let $\omega_\epsilon\in\mathcal F$ and $u\in C_{\overline\Omega}^\infty$.
Then we have the pointwise derivative
\begin{equation}\frac\partial{\partial \epsilon}_{|\epsilon=0}J_{\omega_\epsilon}u(\lambda)
=-J_\omega\left\{\frac{\nu~{\rm div}(u{N})}{|\nabla\omega|}\right\}(\lambda),
\qquad \text{where}\qquad \nu:=\frac{\partial}{\partial \epsilon}_{|\epsilon=0}\omega_\epsilon.
\label{eq:ddepsilonJpsiepsilonu(lambda)}
\end{equation}
\end{lemma}
\startproof
For an arbitrary function $\zeta$ with compact support in $(\min\omega,\max\omega)$,
differentiating the co-area formula (\ref{eq:CoareaFormula})
\begin{equation}
\int_\Omega u(x)|\nabla\omega_\epsilon(x)|\zeta(\omega_\epsilon(x))dx
=\int_{\min\omega}^{\max\omega} \zeta(\lambda)J_{\omega_\epsilon}u(\lambda)~d\lambda
\end{equation}
at $\epsilon=0$ we obtain
\begin{equation}
\int_\Omega u(x)\left[\frac{\langle \nabla \omega,\nabla\nu\rangle}{|\nabla\omega|}\zeta(\omega(x))+|\nabla\omega|\zeta'(\omega(x))\nu(x)\right]~dx
=\left[\int_{\min\omega}^{\max\omega} \zeta(\lambda)\left(\frac{\partial}{\partial\epsilon}_{|\epsilon=0}J_{\omega_\epsilon}u\right)(\lambda)d\lambda\right].
\end{equation}
On the one hand, the coarea formula (\ref{eq:CoareaFormula}) 
and an integration by parts give
\begin{equation}
\int_\Omega u(x)\nu(x)|\nabla\omega|\zeta'(\omega(x))dx 
=\int_{\min\omega}^{\max\omega} \zeta'(\lambda)(J_\omega (u\nu))(\lambda)~d\lambda
=-\int_{\min\omega}^{\max\omega} \zeta(\lambda)(J_\omega (u\nu))'(\lambda)~d\lambda
\end{equation}
and on the other,
\begin{eqnarray}
\int_\Omega u(x)\frac{\langle \nabla \omega,\nabla\nu\rangle}{|\nabla \omega|}\zeta(\omega(x))~dx
&=&\int_{\min\omega}^{\max\omega} \zeta(\lambda)J_\omega \left(u \frac{\langle \nabla \omega,\nabla \nu\rangle}{|\nabla \omega|^2}\right)(\lambda)~d\lambda
\end{eqnarray}
Since $\zeta$ is arbitrary,
and using (\ref{eq:ddlambdaJpsiu}), we conclude with 
\begin{equation}
\frac{\partial}{\partial\epsilon}_{|\epsilon=0}\left(J_{\omega_\epsilon}u(\lambda)\right)
=-(J_\omega(u\nu))'(\lambda)+J_\omega\left(u\frac{\langle \nabla\omega,\nabla \nu\rangle}{|\nabla\omega|^2}\right)(\lambda)
=-J_\omega\left\{\frac{\nu~{\rm div}(u{N})}{|\nabla\omega|}\right\}(\lambda).
\end{equation}
\stopproof%

\subsection{$A_\omega^{-1}$ is a smooth tame map of $\omega$}\label{section:AomegaInverseIsSmoothTame}
We will heavily rely on continuity and tameness of operators introduced in the Appendix. 
Also, we recall the Fa\`a  di Bruno formula: if $f,g$ are two functions of one variable, then the $n$-th derivative ($n\geq 0$) of $f\circ g$ is
\begin{equation}
(f\circ g)^{(n)}
=\sum_{k=1}^n \sum_{\begin{subarray}{c}j_1+\cdots+j_k=n\\j_1,\dots,j_k\geq 1\end{subarray}}c_{k;j_1,\dots,j_k}(f^{(k)}\circ g)g^{(j_1)}\cdots g^{(j_k)},\qquad n\geq 0
\label{eq:nThDerivativeOf(fog)}\end{equation}
where $c_{k;j_1,\dots,j_k}$ are constants.
If $f=g^{-1}$, then $f'=\frac1{g'\circ f}$ and  $(f\circ g)^{(n)}=0$ for $n\geq 2$, hence 
\begin{equation}
f^{(n)}=-(f')^n\sum_{k=1}^{n-1} \sum_{\begin{subarray}{c}j_1+\cdots+j_k=n\\j_1,\dots,j_k\geq 1\end{subarray}}c_{k;j_1,\dots,j_k}f^{(k)}(g^{(j_1)}\circ f)\cdots (g^{(j_k)}\circ f),\qquad n\geq 2.
\label{eq:nThDerivativeOff}
\end{equation}

\begin{proposition}[$Q(\omega)=A_\omega^{-1}$ is smooth tame]\label{proposition:AomegaInverseIsSmoothTame}
\hfill\smallskip%
\begin{enumerate}
\item
The operator
\begin{equation}
Q\colon\left\{
\begin{array}{c}
\mathcal F_+\\\omega\end{array}
\longrightarrow
\begin{array}{c}C^\infty_{[0,|\Omega|]}\\A_\omega^{-1}
\end{array}\right\}
\end{equation}
is a smooth tame map of Fr\'echet spaces with first derivative
\begin{equation}
DQ(\omega)\cdot \nu
=\frac{J_\omega\frac\nu{|\nabla\omega|}\circ A_\omega^{-1}}{J_\omega\frac1{|\nabla\omega|}\circ A_\omega^{-1}}.
\label{eq:DQ(psi)phi}
\end{equation}
For $n\geq 0$, $\omega\in\mathcal F_+$ in a $\|\cdot\|_{2,\alpha}$-neighborhood and any $\nu, \nu_1, \nu_2\in \mathcal F$,
\begin{eqnarray}
\|A_\omega^{-1}\|_{n,\alpha}&\leq& C (\|\omega\|_{n,\alpha}+1),\\
\|DQ(\omega)\nu\|_{n,\alpha}&\leq&C \left(\|\nu\|_{n,\alpha}+\|\omega\|_{n+1,\alpha}\|\nu\|_{1,\alpha}\right),\\
\|D^2Q(\omega)(\nu_1,\nu_2)\|_{n,\alpha}
&\leq& C\Bigg\{\|\nu_1\|_{n+1,\alpha}\|\nu_2\|_{1,\alpha}+\|\nu_1\|_{1,\alpha}\|\nu_2\|_{n+1,\alpha}\\
&&+\|\omega\|_{n+2,\alpha}\|\nu_1\|_{2,\alpha}\|\nu_2\|_{2,\alpha}\Bigg\}.
\end{eqnarray}
\item The operator
\begin{equation}
\left\{
\begin{array}{ccc}
\mathcal F_+&\times& C^\infty_{\overline\Omega}\\\omega&&u\end{array}
\longrightarrow
\begin{array}{c}C^\infty_{[0,|\Omega|]}\\J_\omega u\circ A_\omega^{-1}
\end{array}
\right\}
\end{equation}
is smooth tame.
For $n\geq 0$, $\omega$ in a $\|\cdot\|_{2,\alpha}$-neighborhood and any $u$ (without restriction),
\begin{equation}
\left\|J_\omega\frac{u}{|\nabla\omega|}\circ A_\omega^{-1}\right\|_{n,\alpha}
\leq C\cdot (\|u\|_{n,\alpha}+\|\omega\|_{n+1,\alpha}\|u\|_{1,\alpha}).
\label{ineq:TameEstimatesOnJpsiuOfApsiInverse}
\end{equation}
\end{enumerate}
\end{proposition}
What prevents one from working directly with $A_\omega$ is that the range of $\omega$ is not fixed
and that $A_\omega'$ is necessarily discontinuous at the endpoints of $\text{range}(\omega)$.
This problem is resolved by working instead with the inverse distribution function $A_\omega^{-1}$,
at the cost of a fair amount of technical complications.

Proposition~\ref{proposition:AomegaInverseIsSmoothTame} 
is split into Lemmas~\ref{lemma:QIsContinuous}, \ref{lemma:QIsTame}, and \ref{lemma:QIsSmoothTame}.
For continuity, the $C^n$-grading  will be more convenient.
However, tame estimates will still be derived in the $C^{n,\alpha}$-grading.

\begin{lemma}[$Q(\omega)=A_\omega^{-1}$ is continuous]\label{lemma:QIsContinuous}
$Q(\omega)=A_\omega^{-1}$ is continuous as a map of Fr\'echet spaces.
More precisely, the following are continuous as maps of Banach spaces:
\begin{equation}
Q\colon\left\{\begin{array}{c}\mathcal F_+^{m}\\\omega\end{array}
\longrightarrow\begin{array}{c}C^m_{[0,|\Omega|]}\\A_\omega^{-1}
\end{array}\right\}\qquad (m\geq 2),
\label{eq:ContinuityOfomega->AomegaInverse}
\end{equation}
\begin{equation}
\left\{\begin{array}{ccc}\mathcal F_+^{m}&\times&C^{m-1}_{\overline\Omega}\\\omega&&u\end{array}
\longrightarrow\begin{array}{c}C^{m-1}_{[0,|\Omega|]}\\
(J_\omega u)\circ A_\omega^{-1}
\end{array}\right\}\qquad (m\geq 2),
\end{equation}
\end{lemma}
\startproof%

\setcounter{counterpar}{0}
\nextpar{$\omega\mapsto A_\omega^{-1}$ is continuous $C^m\rightarrow C^{m}$, $m\geq 2$}
In order to alleviate some of the complications of working with $A_\omega^{-1}$ (rather than $A_\omega$ directly),
we will use the following device.
Fix $\omega_1\in \mathcal F_+$.
For $\omega\in \mathcal F_+$, let $\omega_2\in \mathcal F_+$ and $a$ a monotone increasing function such that
$\text{range}(\omega_2)=\text{range}(\omega_1)$ and $\omega=a\circ \omega_2$.  
We will take $a$ to be affine:
\begin{equation}\omega=a\circ\omega_2=\min\omega+\frac{\max\omega-\min\omega}{\max\omega_2-\min\omega_2}(\omega_2-\min\omega_2).
\end{equation}
Fix $k\geq 0$. 
Then
\begin{eqnarray}
\|A_\omega^{-1}-A_{\omega_1}^{-1}\|_k&\leq& \|A_{a\circ \omega_2}^{-1}-A_{\omega_2}^{-1}\|_k+\|A_{\omega_2}^{-1}-A_{\omega_1}^{-1}\|_k\\
&=& \|(a - \text{id})\circ A_{\omega_2}^{-1}\|_k+\|A_{\omega_2}^{-1}-A_{\omega_1}^{-1}\|_k\\
&=&I_k+II_k.
\end{eqnarray}
\stopstep%

\nextpar{Estimating $I_k$}
In view of (\ref{eq:nThDerivativeOf(fog)}),
$I_k$ is arbitrarily small provided $\|a-\text{id}\|_k$ is sufficiently small
while $\|A_{\omega_2}^{-1}\|_k$ remains bounded.
Since $a$ is affine, the former is small provided $\|\omega-\omega_1\|_0$ is small.
Next, we show that $\|A_{\omega_2}^{-1}\|_k$ remains bounded
provided $\|\omega_2\|_k$ is bounded and $\|\omega_2-\omega_1\|_2$ is sufficiently small.

Observe first that $\|A_{\omega_2}^{-1}\|_0=\max (|\max \omega_2|,|\min\omega_2|)=\max(|\max\omega|,|\min\omega|)$ which is fixed.
For $\lambda\in\text{range}(\omega_2)$,
\begin{equation}
A_{\omega_2}'(\lambda)=\int_{\omega_2=\lambda}\frac{dl}{|\nabla\omega_2|}
=\int_{s\in\mathbb S^1}
\frac{\left|\frac{\partial {z_2}}{\partial s}\right|(\frac{\lambda-\min{\omega_2}}{\max{\omega_2}-\min{\omega_2}},s)}{|\nabla{\omega_2}({z_2}(\frac{\lambda-\min{\omega_2}}{\max{\omega_2}-\min{\omega_2}},s))|}
ds.
\end{equation}
Fix then $L_1$ and a neighborhood 
\begin{equation}\|{\omega_2}-\omega_1\|_2<\epsilon_0\label{ineq:NbhdWhereCertainTechnicalBoundsHold}
\end{equation}
where
\begin{equation}\frac1{L_1}<\int_{\omega_2=\lambda}dl<L_1\qquad \text{and}\qquad \frac1{L_1}<A_{\omega_2}'(\lambda)<L_1\end{equation}
for $\lambda\in\text{range}(\omega_2)$.
Now
\begin{equation}\left(A_{\omega_2}^{-1}\right)'(\mu)=\frac1{A_{\omega_2}'\circ A_{\omega_2}^{-1}(\mu)},\qquad \mu\in [0,|\Omega|]\end{equation}
so that $\|A_{\omega_2}^{-1}\|_1$ remains bounded when $\|{\omega_2}-\omega_1\|_2<\epsilon_0$.

For the higher order derivatives, we use (\ref{eq:nThDerivativeOff}):
$\frac{d^nA_{{\omega_2}}^{-1}}{d\mu^n}$, $n\geq 2$, is given by
\begin{eqnarray}
-\left(\frac{dA_{\omega_2}^{-1}}{d\mu}\right)^n\sum_{k=1}^{n-1}\sum_{\begin{subarray}{c}j_1+\dots+j_k=n\\j_1,\dots,j_k\geq 1\end{subarray}}c_{k;j_i}
\left(\frac{d^kA_{\omega_2}^{-1}}{d\mu^k}\right)
\left(\frac{d^{j_1}A_{\omega_2}}{d\lambda^{j_1}}\circ A_{\omega_2}^{-1}\right)
\cdots \left(\frac{d^{j_k}A_{\omega_2}}{d\lambda^{j_k}}\circ A_{\omega_2}^{-1}\right).
\label{eq:nThDerivativeOfApsiInverse}
\end{eqnarray}
An induction shows that $\frac{d^nA_{\omega_2}^{-1}}{d\mu^n}$ is bounded provided $\frac{dA_{\omega_2}}{d\mu}$,\dots, $\frac{d^nA_{\omega_2}}{d\mu^n}$.
We show that this holds when $\|\omega_2\|_n$ is bounded.
Setting
\begin{equation}v_0=\frac1{|\nabla{\omega_2}|},\qquad v_m=\frac{{\rm div}(v_{m-1}N_2)}{|\nabla{\omega_2}|}, \qquad m\geq 1
\end{equation}
where $N_2=\frac{\nabla{\omega_2}}{|\nabla{\omega_2}|}$,
Lemma~\ref{lemma:ddlambdaJpsiu(lambda)} yields
\begin{eqnarray}
&&\frac{d^jA_{\omega_2}(\lambda)}{d\lambda^j}\\
&=&(J_{\omega_2} v_{j-1})(\lambda)=\int_{{\omega_2}=\lambda}v_{j-1}dl\\
&=&\int_{s\in\mathbb S^1}v_{j-1}\left(z\left((\frac{\lambda-\min\omega_2}{\max\omega_2-\min\omega_2},s\right)\right)
\left|\frac{\partial z}{\partial s}\right|\left(\frac{\lambda-\min\omega_2}{\max\omega_2-\min\omega_2},s\right)ds.
\label{eq:djAomega2(lambda)dlambdaj}
\end{eqnarray}
But $v_{j-1}$ is a smooth expression of the derivatives of ${\omega_2}$ up to order $j$.
Therefore, with $\|{\omega_2}-\omega_1\|_2<\epsilon_0$ and $\|{\omega_2}\|_j$ bounded,
$\|A_{\omega_2}\|_{C^j_{\text{range}({\omega_2})}}$ is bounded.
\stopstep%

\nextpar{Estimating $II_k$}
Since $\omega_2$ and $\omega_1$ have same range, 
we may invoke general results on the inversion operator, see Lemma~\ref{lemma:InversionIsSmoothTame} in the Appendix:
for $k\geq 1$, $II_k$ is arbitrarily small
provided $\|A_{\omega_2}-A_{\omega_1}\|_{C^k_{\text{range}(\omega_1)}}$ is taken sufficiently small.
Restricting to functions such that $\text{range}(\omega_2)=\text{range}(\omega_1)$,
it remains to show that $\omega_2\mapsto A_{\omega_2}$ is continuous $C^k\rightarrow C^{k}$, $k\geq 2$.
But this is immediate in view of (\ref{eq:djAomega2(lambda)dlambdaj})
and the fact that $\omega\mapsto \frac{\partial z}{\partial s}$ is continuous $C^j\rightarrow C^0$ for $j\geq 2$.
\stopstep%

\nextpar{$(\omega,u)\mapsto (J_\omega u)\circ A_\omega^{-1}$ is continuous $C^m\times C^{m-1}\rightarrow C^{m-1}$, $m\geq 2$}
Let $u,u_1$, and $\omega, \omega_1$ with corresponding coordinate systems $z, z_1$.
For $\mu\in [0,|\Omega|]$, write$\lambda=A_\omega^{-1}(\mu)$, $\lambda_1=A_{\omega_1}^{-1}(\mu)$,
and $t=\frac{\lambda-\min\omega}{\max\omega-\min\omega}$, $t_1=\frac{\lambda_1-\min\omega_1}{\max\omega_1-\min\omega_1}$.
Then
\begin{eqnarray}
&&J_\omega u\circ A_\omega^{-1}(\mu)-J_{\omega_1}u_1\circ A_{\omega_1}^{-1}(\mu)\\
&=&\int_{\omega=\lambda}udl - \int_{\omega_1=\lambda_1}u_1dl\\
&=&\int_{s\in\mathbb S^1}
\left(u(z(t,s))\left|\frac{\partial z}{\partial s}\right|(t,s)-u_1(z_1(t_1,s))\left|\frac{\partial z_1}{\partial s}\right|(t_1,s)\right)
ds
\end{eqnarray}
and the integrand is small (uniformly in $s$ and $\mu$) provided $\|u-u_1\|_0$ and $\|\omega-\omega_1\|_2$ are taken sufficiently small.
This shows that $(\omega,u)\mapsto J_\omega u\circ A_\omega^{-1}$ is continuous $C^2\times C^0\rightarrow C^0$.

Write
$\frac{d(J_\omega u\circ A_\omega^{-1})}{d\mu}=\left(J_\omega\frac{{\rm div}(uN)}{|\nabla\omega|}\circ A_\omega^{-1}\right)\times \left(\frac{dA_\omega^{-1}}{d\mu}\right)$
where $N=\frac{\nabla\omega}{|\nabla\omega|}$.
Since $(\omega,u)\mapsto \frac{{\rm div}(uN)}{|\nabla\omega|}$ is continuous $C^2\times C^1\rightarrow C^0$
and $\omega\mapsto A_\omega^{-1}$ is continuous $C^2\rightarrow C^1$, we conclude thanks to the previous paragraph
that $(\omega,u)\mapsto J_\omega u\circ A_\omega^{-1}$ is continuous $C^2\times  C^1\rightarrow C^1$.

As for higher order derivatives, 
set $v_0=u$, $v_m=\frac{{\rm div}(v_{m-1}N)}{|\nabla\omega|}$, $m\geq 1$.
A simple induction shows that 
$(\omega, u)\mapsto v_m$ is continuous $\mathcal F_+^{m+1}\times C^m_{\overline\Omega}\rightarrow C^0_{\overline\Omega}$.
Now from (\ref{eq:nThDerivativeOf(fog)}) we find 
\begin{eqnarray}
&&\frac{d^n(J_\omega u\circ A_\omega^{-1})}{d\mu^n}\\
&=&\sum_{k=1}^n\sum_{\begin{subarray}{c}j_1+\cdots+j_k=n\\j_1, \dots, j_k\geq 1\end{subarray}}
c_{k;j_i}\left(\frac{d^kJ_\omega u}{dc^k}\circ A_\omega^{-1}\right)
\left(\frac{d^{j_1}A_\omega^{-1}}{d\mu^{j_1}}\right)
\dots\left(\frac{d^{j_k}A_\omega^{-1}}{d\mu^{j_k}}\right)\label{eq:nThDerivativeOfJpsiuOfApsiInverse}\\
&=&\sum_{k=1}^n\sum_{\begin{subarray}{c}j_1+\cdots+j_k=n\\j_1, \dots, j_k\geq 1\end{subarray}}
c_{k;j_i}\left(J_\omega u_k\circ A_\omega^{-1}\right)
\left(\frac{d^{j_1}A_\omega^{-1}}{d\mu^{j_1}}\right)
\dots\left(\frac{d^{j_k}A_\omega^{-1}}{d\mu^{j_k}}\right).
\end{eqnarray}
Again by induction $(\omega,u)\mapsto J_\omega u\circ A_\omega^{-1}$ is continuous $C^n\times C^{n-1}\rightarrow C^{n-1}$,
$n\geq 2$.
\stopproof%

\begin{lemma}[$Q(\omega)=A_\omega^{-1}$ is tame]\label{lemma:QIsTame}
\hfill\smallskip%
\begin{enumerate}
\item For $n\geq 0$ and $\omega$ in $\|\cdot\|_{2,\alpha}$-neighborhood,
\begin{equation}
\|A_\omega^{-1}\|_{n,\alpha}\leq C\cdot (\|\omega\|_{n,\alpha}+1).
\label{ineq:TameEstimatesOnQ}
\end{equation}
\item
For $n\geq 0$, $\omega$ in a $\|\cdot\|_{2,\alpha}$-neighborhood and any  $u$ (without restriction),
\begin{equation}
\|J_\omega u\circ A_\omega^{-1}\|_{n,\alpha}\leq C\cdot (\|u\|_{n,\alpha}+\|\omega\|_{n+1,\alpha}\|u\|_{1,\alpha}).
\label{ineq:TameEstimateOnJpsiuMpsiInverse}
\end{equation}
\end{enumerate}
\end{lemma}
\pagebreak
\startproof%

\setcounter{counterpar}{0}
\nextpar{Estimate on $\|J_\omega u\circ A_\omega^{-1}\|_{0,\alpha}$}
Recall that $\epsilon_0$ is defined in (\ref{ineq:NbhdWhereCertainTechnicalBoundsHold}).
Clearly,
\begin{equation}\|J_\omega u\circ A_\omega^{-1}\|_0\leq L_1\|u\|_0.\label{ineq:||JomegauOfAomegaInverse||0}
\end{equation}
To estimate the H\"older-constant of $J_\omega u\circ A_\omega^{-1}$,
fix $\mu,\mu'\in [0,|\Omega|]$ and set
$\lambda=A_\omega^{-1}(\mu)$, $\lambda'=A_\omega^{-1}(\mu')$, $t=\frac{\lambda-\min\omega}{\max\omega-\min\omega}$, $t'=\frac{\lambda'-\min\omega}{\max\omega-\min\omega}$.
Then,
\begin{eqnarray}
&&J_\omega u(\lambda')-J_\omega u(\lambda)\\
&=&\int_{s\in\mathbb S^1}\left(u(z(t',s))\left|\frac{\partial z}{\partial s}\right|(t',s)
-u(z(t,s))\left|\frac{\partial z}{\partial s}\right|(t,s)\right)ds\\
&=&\int_{s\in\mathbb S^1}u(z(t',s))\left(\left|\frac{\partial z}{\partial s}\right|(t',s)-\left|\frac{\partial z}{\partial s}\right|(t,s)\right)\\
&&+\left(u(z(t',s))-u(z(t,s))\right)\left|\frac{\partial z}{\partial s}\right|(t,s)ds.
\end{eqnarray}
Since $\omega$ is restricted to a $\|\cdot\|_{2,\alpha}$-neighborhood, 
Lemma~\ref{lemma:CoordinatesDependingOnomega} gives that
the first term is bounded by $C\|u\|_0|t'-t|$ 
while the second is bounded by $C[u]_\alpha|z(t',s)-z(t,s)|^\alpha\leq C [u]_\alpha|t'-t|^\alpha$.
Now from the proof of Lemma~\ref{lemma:QIsContinuous}, $(A_\omega^{-1})'$ is bounded, 
so that $|t'-t|=C|\lambda'-\lambda|\leq C|\mu'-\mu|$.
In conclusion, for $\omega$ in a $\|\cdot\|_{2,\alpha}$-neighborhood and any $u$,
\begin{equation}\|J_\omega u\circ A_\omega^{-1}\|_{0,\alpha}\leq C\|u\|_{0,\alpha}
\label{ineq:||JomegauOfAomegaInverse||0,alpha}\end{equation}
for $\omega$ in a $\|\cdot\|_{2,\alpha}$-neighborhood and any $u$ (without restriction).
This implies (\ref{ineq:TameEstimateOnJpsiuMpsiInverse}) for $m=0$.
\stopstep%

\nextpar{Estimate on $\|A_\omega^{-1}\|_{1,\alpha}$}
With $\left(A_\omega^{-1}\right)'=\frac1{J_\omega \frac1{|\nabla\omega|}\circ A_\omega^{-1}}$,
with $\|\omega-\overline\omega\|_1$ sufficiently small $\frac1{|\nabla\omega|}$ remains bounded,
and thus making $\epsilon_0$ smaller if necessary, $\left(A_\omega^{-1}\right)'$ remains bounded in $\|\cdot\|_{0,\alpha}$
for $\|\omega-\overline\omega\|_{2,\alpha}<\epsilon_0$.
\stopstep%

\nextpar{Estimate on $\|u_m\|_{0,\alpha}$}
Suppose first (for simplicity) that $u, a$ are smooth functions of one variable and set
$u_0:=u$, $u_{m}:=(u_{m-1}a)'a$, $m\geq 1$.
By induction, one verifies that $u_m$ is then of the form
\begin{equation}
u_m= \sum_{\begin{subarray}{c}j_0+\cdots+j_m=m\\j_0,\dots,j_m\geq 1\end{subarray}}c_{m;j_1,\dots,j_m}u^{(j_0)}a^{(j_1)}\cdots a^{(j_m)}a^m.
\end{equation}
Using (\ref{ineq:TameEstimatesForProductOfFunctions}), 
and interpolation inequalities (\ref{ineq:IntroInterpolationInequalities}) on all factors
(between their $\|\cdot\|_{0,\alpha}$ and $\|\cdot\|_{m-j_0,\alpha}$-norms),
\begin{eqnarray}
\|u_m\|_{0,\alpha}
&\leq&C\sum_{\begin{subarray}{c}j_0+\cdots+j_m=m\\j_0,\dots,j_m\geq 1\end{subarray}}\|u\|_{j_0,\alpha}\|a\|_{j_1,\alpha}\cdots\|a\|_{j_m,\alpha}\|a^m\|_{0,\alpha}\\
&\leq&C\|a\|_{0,\alpha}^{l_m}\sum_{\begin{subarray}{c}j_0+\cdots+j_m=m\\j_0,\dots, j_m\geq 0\end{subarray}}
\|u\|_{j_0,\alpha}\|a\|_{m-j_0,\alpha}^\frac{j_1}{m-j_0}\cdots
\|a\|_{m-j_0,\alpha}^\frac{j_m}{m-j_0}\\
&\leq &C\|a\|_{0,\alpha}^{l_m}\|u\|_{j_0,\alpha}\|a\|_{m-j_0,\alpha}\\
&\leq &C\|a\|_{0,\alpha}^{l_m}(\|u\|_{0,\alpha}\|a\|_{m,\alpha}+\|u\|_{m,\alpha}\|a\|_{0,\alpha})
\end{eqnarray}
where $l_m$ is some positive integer depending on $m$.
In turn,
\begin{equation}
\|u_m\|_{0,\alpha}\leq C(\|u\|_{m,\alpha}+\|a\|_{m,\alpha}\|u\|_{0,\alpha}),\qquad m\geq 1
\end{equation} for all $u$ without restriction, and $a$ in a $\|\cdot\|_{0,\alpha}$-neighborhood.

The situation with $u_m$ as defined by $u_m=\frac{{\rm div}(u_{m-1}{N})}{|\nabla\omega|}$
can be dealt with in a similar fashion, only the details are more tedious.
Here, $a$ plays the r\^ole of $\frac1{|\nabla\omega|}$ or $N=\frac{\nabla\omega}{|\nabla\omega|}$.
In conclusion,
\begin{equation}\|u_m\|_{0,\alpha}\leq C (\|u\|_{m,\alpha}+\|\omega\|_{m+1,\alpha}\|u\|_{0,\alpha}),\qquad m\geq 1
\label{ineq:EstimateOn||um||0,alpha}
\end{equation}
for all $\omega$ in a $\|\cdot\|_{1,\alpha}$-neighborhood and all $u$ without restriction.
\stopstep %

\nextpar{Estimates on $\|J_\omega u_m\circ A_\omega^{-1}\|_{0,\alpha}$}
We easily conclude from (\ref{ineq:||JomegauOfAomegaInverse||0,alpha}) 
and (\ref{ineq:EstimateOn||um||0,alpha}) that
\begin{equation}
\left\|\frac{d^mJ_\omega u}{d\lambda^m}\circ A_\omega^{-1}\right\|_{0,\alpha}=
\left\|J_\omega u_m\circ A_\omega^{-1}\right\|_{0,\alpha}
\leq C\cdot (\|u\|_{m,\alpha}+\|\omega\|_{m+1,\alpha}\|u\|_{0,\alpha}),\qquad m\geq 0
\label{ineq:||dmdlambdamJomegauOfAomegaInverse||0,alpha}
\end{equation}
for $\omega$ in a $\|\cdot\|_{2,\alpha}$-neighborhood.
\stopstep

\nextpar{Estimate (\ref{ineq:TameEstimatesOnQ}) on $\|A_\omega^{-1}\|_{n,\alpha}$}
In this paragraph only, we write
\begin{equation}f=A_\omega^{-1},\qquad g=A_\omega,\qquad \text{and}\qquad J=J_\omega \frac1{|\nabla\omega|}
\end{equation}
with derivatives $f'$, $g'$, $J'$.
Since $g'=J$,
\begin{eqnarray}
f^{(n)}
&=&-(f')^n\sum_{k=1}^{n-1}\sum_{\begin{subarray}{c}j_1+\dots+j_k=n\\j_1,\dots,j_k\geq 1\end{subarray}} c_{k;j_1,\dots, j_k}f^{(k)}(J^{(j_1-1)}\circ f)\cdots (J^{(j_k-1)}\circ f)
\end{eqnarray}
hence by (\ref{ineq:TameEstimatesForProductOfFunctions})
\begin{eqnarray}
\|f^{(n)}\|_{0,\alpha}
\leq C\|f'\|_{0,\alpha}^n\sum_{k=1}^{n-1}\sum_{\begin{subarray}{c}j_1+\dots+j_k=n\\j_1,\dots,j_k\geq 1\end{subarray}}\|f^{(k)}\|_{0,\alpha}
\|J^{(j_1-1)}\circ f\|_{0,\alpha}\cdots \|J^{(j_k-1)}\circ f\|_{0,\alpha}.
\end{eqnarray}
We prove by induction the estimate (\ref{ineq:TameEstimatesOnQ}).
We have already seen that $\|f\|_{1,\alpha}$ is bounded for $\omega$ in a $\|\cdot\|_{2,\alpha}$-neighborhood.
Suppose the estimate (\ref{ineq:TameEstimatesOnQ}) proved up to some $n-1\geq 1$.
From (\ref{ineq:||dmdlambdamJomegauOfAomegaInverse||0,alpha})
$\|J^{(j-1)}\circ f\|_{0,\alpha}\leq C\cdot (\|\omega\|_{j,\alpha}+1)$, $j\geq 1$,
hence by the induction hypothesis
\begin{eqnarray}
\|f^{(n)}\|_{0,\alpha}
\leq C\|f'\|_{0,\alpha}^n\sum_{k=1}^{n-1}\sum_{\begin{subarray}{c}j_1+\dots+j_k=n\\j_1,\dots,j_k\geq 1\end{subarray}}(\|\omega\|_{k,\alpha}+1)
(\|\omega\|_{j_1,\alpha}+1)\cdots (\|\omega\|_{j_k,\alpha}+1).
\end{eqnarray}
The double sum is the sum of $1$, $\|\omega\|_{k,\alpha}\|\omega\|_{j_1,\alpha}\cdots \|\omega\|_{j_k,\alpha}$,
and products of fewer factors.
For simplicity, consider only the  term $\|\omega\|_{k,\alpha}\|\omega\|_{j_1,\alpha}\cdots \|\omega\|_{j_k,\alpha}$
(the other terms are in fact easier).
We interpolate each factor between its $\|\cdot\|_{1,\alpha}$- and $\|\cdot\|_{n,\alpha}$-norms
using (\ref{ineq:IntroInterpolationInequalities})
\begin{eqnarray}
&&\|\omega\|_{k,\alpha}\|\omega\|_{j_1,\alpha}\cdots \|\omega\|_{j_k,\alpha}\\
&\leq& C\cdot \|\omega\|_{1,\alpha}^\frac{n-k}{n-1}\|\omega\|_{n,\alpha}^\frac{k-1}{n-1}\cdot
\|\omega\|_{1,\alpha}^\frac{n-j_1}{n-1}\|\omega\|_{n,\alpha}^\frac{j_1-1}{n-1}\cdots
\|\omega\|_{1,\alpha}^\frac{n-j_k}{n-1}\|\omega\|_{n,\alpha}^\frac{j_k-1}{n-1}\\
&\leq&C\cdot \|\omega\|_{n,\alpha}\label{ineq:ArriveAt||psi||n,alpha}
\end{eqnarray}
since $\|\omega\|_{1,\alpha}$ remains bounded and $j_1+\cdots+j_k=n$.
This establishes the tame estimate (\ref{ineq:TameEstimatesOnQ}) on $A_\omega^{-1}$ for 
$\omega$ in a $\|\cdot\|_{2,\alpha}$-neighborhood.
\stopstep

\nextpar{$J_\omega u\circ A_\omega^{-1}$ is tame}
Setting $f=J_\omega u$ and $g=A_\omega^{-1}$ in (\ref{eq:nThDerivativeOf(fog)}) and with
\begin{equation}u_0=u,\qquad u_m=\frac{{\rm div}(u_{m-1}N)}{|\nabla\omega|},\qquad m\geq 1,\label{eq:um}
\end{equation}
we have
\begin{eqnarray}
\frac{d^n(J_\omega u\circ A_\omega^{-1})}{d\mu^n}
=\sum_{k=1}^n\sum_{\begin{subarray}{c}j_1+\cdots+j_k=n\\j_1,\dots,j_k\geq 1\end{subarray}} c_{k;j_i}\left(J_\omega u_k\circ A_\omega^{-1}\right)
\left(\frac{d^{j_1}A_\omega^{-1}}{d\mu^{j_1}}\right)\cdots \left(\frac{d^{j_k}A_\omega^{-1}}{d\mu^{j_k}}\right).
\end{eqnarray}
From (\ref{ineq:||dmdlambdamJomegauOfAomegaInverse||0,alpha}), (\ref{ineq:TameEstimatesOnQ}), and (\ref{ineq:EstimateOn||um||0,alpha}),
and by (\ref{ineq:TameEstimatesForProductOfFunctions}),
we have for $\omega$ in a $\|\cdot\|_{2,\alpha}$-neighborhood and $u$ in a $\|\cdot\|_{0,\alpha}$-neighborhood,
\begin{eqnarray}
&&\left\|\left(J_\omega u_k\circ A_\omega^{-1}\right)
\left(\frac{d^{j_1}A_\omega^{-1}}{d\mu^{j_1}}\right)\cdots \left(\frac{d^{j_k}A_\omega^{-1}}{d\mu^{j_k}}\right)\right\|_{0,\alpha}\\
&\leq& C\cdot (\|u\|_{k,\alpha}+\|\omega\|_{k+1,\alpha})(\|\omega\|_{j_1,\alpha}+1)\cdots (\|\omega\|_{j_k,\alpha}+1).
\end{eqnarray}
Products of $\|\omega\|_{j,\alpha}$'s can be estimated by $\|\omega\|_{n,\alpha}$ in a way similar to that 
which led to (\ref{ineq:ArriveAt||psi||n,alpha}).
Consider now $\|u\|_{k,\alpha}\|\omega\|_{j_1,\alpha}\cdots\|\omega\|_{j_k,\alpha}$
(the remaining terms are easier and their estimation will be omitted).
Interpolating each factor between its $\|\cdot\|_{1,\alpha}$ and $\|\cdot\|_{n,\alpha}$-norms
with (\ref{ineq:IntroInterpolationInequalities})
and using the inequality $x^\delta y^{1-\delta}\leq x+y$ ($0\leq\delta\leq 1$), 
we find
\begin{eqnarray}
\|u\|_{k,\alpha}\|\omega\|_{j_1,\alpha}\cdots \|\omega\|_{j_k,\alpha}
&\leq& C\|u\|_{1,\alpha}^\frac{n-k}{n-1}\|u\|_{n,\alpha}^\frac{k-1}{n-1}\|\omega\|_{1,\alpha}^\frac{(k-1)n}{n-1}\|\omega\|_{n,\alpha}^\frac{n-k}{n-1}\\
&\leq& C\cdot\|\omega\|_{1,\alpha}^{k-1}\Big[\|u\|_{1,\alpha}\|\omega\|_{n,\alpha}+\|u\|_{n,\alpha}\|\omega\|_{1,\alpha}\Big]\\
&\leq& C\cdot (\|u\|_{n,\alpha}+\|\omega\|_{n,\alpha})
\end{eqnarray}
for $\omega$ in a $\|\cdot\|_{2,\alpha}$-neighborhood and $u$ in a $\|\cdot\|_{1,\alpha}$-neighborhood.
Replacing $u$ by $\nabla\omega$ one immediately obtains
\begin{equation}
\|\omega\|_{k+1,\alpha}\|\omega\|_{j_1,\alpha}\cdots \|\omega\|_{j_k,\alpha}\leq C\cdot \|\omega\|_{n+1,\alpha},\qquad n\geq 1\end{equation}
for $\omega$ in a $\|\cdot\|_{2,\alpha}$-neighborhood.
Putting these together, we arrive at 
\begin{eqnarray}
&&\left\|\left(J_\omega u_k\circ A_\omega^{-1}\right)
\left(\frac{d^{j_1}A_\omega^{-1}}{d\mu^{j_1}}\right)\cdots \left(\frac{d^{j_k}A_\omega^{-1}}{d\mu^{j_k}}\right)\right\|_{0,\alpha}\\
&\leq& C\cdot (\|\omega\|_{n+1,\alpha}+\|u\|_{n,\alpha})
\end{eqnarray}
and summing, 
\begin{equation}
\|J_\omega u\circ A_\omega^{-1}\|_{n,\alpha}\leq C\cdot (\|\omega\|_{n+1,\alpha}+\|u\|_{n,\alpha}+1),
\qquad n\geq 1
\end{equation}
for $\omega$ in a $\|\cdot\|_{2,\alpha}$-neighborhood and all $u$ in a $\|\cdot\|_{1,\alpha}$-neighborhood.
Since $J_\omega u\circ A_\omega^{-1}$ is linear in $u$, see Proposition~\ref{proposition:EquivalentFormsOfTameEstimatesForFamilyOfLinearMaps},
\begin{equation}\|J_\omega u\circ A_\omega^{-1}\|_{n,\alpha}\leq C\cdot (\|u\|_{n,\alpha}+\|\omega\|_{n+1,\alpha}\|u\|_{1,\alpha}),
\qquad n\geq 1
\end{equation}
for $\omega$ in a $\|\cdot\|_{2,\alpha}$-neighborhood and all $u$ without restriction.
This is (\ref{ineq:TameEstimateOnJpsiuMpsiInverse}) for $m\geq 1$.
\stopproof %

\pagebreak

Even though $Q$ is in fact smooth, two derivatives are sufficient for the Moser iteration.
\begin{lemma}[$Q(\omega)=A_\omega^{-1}$ is smooth tame]\label{lemma:QIsSmoothTame}
$Q\colon\mathcal F_+\rightarrow C^\infty_{[0,|\Omega|]}$ is twice continuously differentiable
as a map of Fr\'echet spaces.
Its first derivative is given by
\begin{equation}
DQ(\omega)\cdot \nu
=\frac{J_\omega\frac\nu{|\nabla\omega|}\circ A_\omega^{-1}}{J_\omega\frac1{|\nabla\omega|}\circ A_\omega^{-1}}.
\label{eq:DQ(omega)nu}
\end{equation}
For $n\geq 0$, $\omega$ in a $\|\cdot\|_{2,\alpha}$-neighborhood and any $\nu$ without restriction, 
we have 
\begin{eqnarray}
\left\|J_\omega\frac{\nu}{|\nabla\omega|}\circ A_\omega^{-1}\right\|_{n,\alpha} \quad\text{and}\quad
\|DQ(\omega)\nu\|_{n,\alpha} \leq C\left\{\|\nu\|_{n,\alpha}+\|\omega\|_{n+1,\alpha}\|\nu\|_{1,\alpha}\right\}.
\end{eqnarray}
The first derivative
\begin{equation}
\left\{
\begin{array}{ccc}
\mathcal F_+^m&\times&C^{m-1}_{\overline\Omega}\\\omega&&\nu\end{array}
\longrightarrow\begin{array}{c}C^{m-1}_{\overline\Omega}\\DQ(\omega)\nu
\end{array}\right\}\qquad (m\geq 2)
\end{equation}
is continuous as a map of Banach spaces.

More generally, 
$(\omega,u)\in\mathcal F_+\times C^\infty_{0,|\Omega|]}\mapsto J_\omega u\circ A_\omega^{-1}\in C^\infty_{[0,|\Omega|]}$
is smooth tame.
For $\omega$ in a $\|\cdot\|_{2,\alpha}$-neighborhood, $u\in C^\infty_{\overline\Omega}$, $\nu\in\mathcal F$,
the first derivative in $\omega$ satisfies 
\begin{equation}
\left\|D_\omega(J_\omega u\circ A_\omega^{-1})\nu\right\|_{n,\alpha}
\leq C\bigg\{\|u\|_{n+1,\alpha}\|\nu\|_{1,\alpha}+\|u\|_{2,\alpha}\|\nu\|_{n,\alpha}+\|\omega\|_{n+2,\alpha}\|u\|_{2,\alpha}\|\nu\|_{1,\alpha}\bigg\}.
\label{ineq:TameEstimatesForDerivativeOfJ(omega)uOfA(omega)InverseInomega}
\end{equation}
For $\omega\in\mathcal F_+$ in a $\|\cdot\|_{2,\alpha}$-neighborhood, $\nu_1,\nu_2\in\mathcal F$, and $n\geq 0$,
\begin{equation}
\|D^2Q(\omega)(\nu_1,\nu_2)\|_{n,\alpha}
\leq C\Bigg\{\|\nu_1\|_{n+1,\alpha}\|\nu_2\|_{1,\alpha}+\|\nu_1\|_{1,\alpha}\|\nu_2\|_{n+1,\alpha}
+\|\omega\|_{n+2,\alpha}\|\nu_1\|_{2,\alpha}\|\nu_2\|_{2,\alpha}\Bigg\}.
\end{equation}
\end{lemma}
\startproof\\
\setcounter{counterpar}{0}
\nextpar{Preliminary remark}
Let $f(\mu,\epsilon)$ be a smooth function on $[0,|\Omega|]\times[0,\epsilon_0]$.
Then, 
\begin{equation}
\frac{f(\mu,\epsilon)-f(\mu,0)}\epsilon \rightarrow_\epsilon \left(\frac{\partial f}{\partial \epsilon}\right)(\mu,0)
\end{equation}
uniformly in $\mu\in[0,|\Omega|]$, and thus 
\begin{equation}
\left\|\frac{f(\cdot,\epsilon)-f(\cdot,0)}\epsilon - \left(\frac{\partial f}{\partial \epsilon}\right)(\cdot,0)\right\|_0\rightarrow_\epsilon 0.
\end{equation}
Since $f$ is smooth, the same holds for all derivatives $\frac{\partial^nf}{\partial \mu^n}$, i.e.
\begin{equation}
\left\|\frac{f(\cdot,\epsilon)-f(\cdot,0)}\epsilon - \left(\frac{\partial f}{\partial \epsilon}\right)(\cdot,0)\right\|_n\rightarrow_\epsilon 0.
\end{equation}
In other words, $\frac{f(\cdot,\epsilon)-f(\cdot,0)}\epsilon\rightarrow_\epsilon \left(\frac{\partial f}{\partial \epsilon}\right)(\cdot,0)$
in the $C^\infty$-topology.
\stopstep%

\nextpar{Differentiability of $Q\colon \mathcal F_+\rightarrow C^\infty_{[0,|\Omega|]}$}
Let $\omega\in \mathcal F_+$, $\nu\in \mathcal F$,
and set $\omega_\epsilon=\omega+\epsilon\nu$.
From the previous paragraph, it is enough to show that $A_{\omega_\epsilon}^{-1}(\mu)$ is a smooth function of $\mu$ and $\epsilon$.
But by definition, we have $A_{\omega_\epsilon}^{-1}(A_{\omega_\epsilon}(\lambda))=\lambda$.
The classical Implicit Function Theorem with parameter $\epsilon$ shows that $A_{\omega_\epsilon}^{-1}(\mu)$ is 
smooth in $\mu$ and $\epsilon$ provided $A_{\omega_\epsilon}(\lambda)$ is smooth in $\lambda$ and $\epsilon$.
But observe that, using the change of coordinates $z^\epsilon$ corresponding to $\omega_\epsilon$, 
see Lemma~\ref{lemma:CoordinatesDependingOnomega}, 
\begin{eqnarray}
\frac{d}{d\lambda}A_{\omega_\epsilon}(\lambda)
&=&J_{\omega_\epsilon}\frac1{|\nabla\omega_\epsilon|}(\lambda)\\
&=&\int_{\omega_\epsilon=\lambda}\frac1{|\nabla\omega_\epsilon|}dl\\
&=&\int_{s\in\mathbb S^1}\frac1{|\nabla\omega_\epsilon\left(z^\epsilon(t^\epsilon,s)\right)|}\left|\frac{\partial z^\epsilon}{\partial s}(t^\epsilon,s)\right|ds
\end{eqnarray}
where $t^\epsilon=\frac{\lambda - \min\omega_\epsilon}{\max \omega_\epsilon-\min\omega_\epsilon}$
is obviously a smooth function of $\lambda$ and $\epsilon$.
This shows that $\frac{d}{d\lambda}A_{\omega_\epsilon}(\lambda)$ is smooth in $\lambda$ and $\epsilon$,
hence that $A_{\omega_\epsilon}(\lambda)$ is, as desired.
\stopstep %

\nextpar{First derivative $DQ(\omega)\cdot\nu$}
Differentiating $A_{\omega_\epsilon}^{-1}(A_{\omega_\epsilon}(\lambda))=\lambda$ at $\epsilon=0$, we find 
\begin{equation}
\frac{\partial}{\partial \epsilon}_{|\epsilon=0}A_{\omega_\epsilon}^{-1}(\mu)
+\frac{dA_\omega^{-1}}{d\mu}(\mu)\frac{\partial}{\partial\epsilon}_{|\epsilon=0}A_{\omega_\epsilon}(\lambda)=0,
\qquad\text{where}\qquad  \lambda=A_\omega^{-1}(\mu)
\end{equation}
so that, thanks to (\ref{eq:ddepsilonAomegaepsilon}),
\begin{equation}
\frac{\partial}{\partial \epsilon}_{|\epsilon=0}A_{\omega_\epsilon}^{-1}(\mu)
=\frac{dA_\omega^{-1}}{d\mu}(\mu)~\left(J_\omega\frac\nu{|\nabla\omega|}\circ A_\omega^{-1}\right)(\mu).\end{equation}
i.e.
\begin{eqnarray}
DQ(\omega)\cdot \nu
&=& \frac{J_\omega\frac\nu{|\nabla\omega|}\circ A_\omega^{-1}}{J_\omega\frac1{|\nabla\omega|}\circ A_\omega^{-1}}
\end{eqnarray}
which is a rational function of continuous tame maps of $\omega$ and $\nu$.
In particular, $DQ(\omega)\cdot \nu$ is continuous as a map of Fr\'echet spaces and tame in $\omega$ and $\nu$ by Lemma~\ref{lemma:QIsContinuous}.
More precisely, Lemma~\ref{lemma:QIsContinuous} implies that
\begin{equation}
\left\{
\begin{array}{ccc}\mathcal F_+^m&\times&C^{m-1}_{\overline\Omega}\\\omega&&\nu\end{array}
\longrightarrow\begin{array}{c}C^{m-1}_{[0,|\Omega|]}\\DQ(\omega)\nu
\end{array}\right\}
\end{equation}
is a continuous map of Banach spaces,
and Lemma~\ref{lemma:QIsTame} that it is tame.
\stopstep 

\nextpar{Tame estimates on $J_\omega\frac\nu{|\nabla\omega|}\circ A_\omega^{-1}$}
For $\omega$ in a $\|\cdot\|_{2,\alpha}$-neighborhood and any $\nu$,
we have for $n\geq 0$
\begin{eqnarray}
&&\left\|J_\omega\frac{\nu}{|\nabla\omega|}\circ A_\omega^{-1}\right\|_{n,\alpha}\\
&\leq& C\left(\left\|\frac{\nu}{|\nabla\omega|}\right\|_{n,\alpha}+\|\omega\|_{n+1,\alpha}\left\|\frac{\nu}{|\nabla\omega|}\right\|_{1,\alpha}\right)\\
&\leq& C\Big\{\|\nu\|_{n,\alpha}
+\|\nu\|_{0,\alpha}(\|\omega\|_{n+1,\alpha}+1)
+\|\omega\|_{n+1,\alpha}\|\nu\|_{1,\alpha}\Big\}\\
&\leq& C\left(\|\nu\|_{n,\alpha}+\|\omega\|_{n+1,\alpha}\|\nu\|_{1,\alpha}\right)
\label{ineq:TameEstimatesForJ(omega)FractionOfA(omega)Inverse}
\end{eqnarray}
by (\ref{ineq:TameEstimateOnJpsiuMpsiInverse}), (\ref{ineq:TameEstimatesForProductOfFunctions}),
and (\ref{ineq:TameEstimatesForNemitskiiOperator}).
\stopstep%

\nextpar{Tame estimates on $DQ(\omega)\nu$}
For $\omega$ in a $\|\cdot\|_{2,\alpha}$-neighborhood and any $\nu$,
we have for $n\geq 0$
\begin{eqnarray}
\left\|DQ(\omega)\nu\right\|_{n,\alpha}
&\leq& \left\|J_\omega\frac\nu{|\nabla\omega|}\circ A_\omega^{-1}\right\|_{n,\alpha}\left\|\frac1{J_\omega\frac1{|\nabla\omega|}\circ A_\omega^{-1}}\right\|_{0,\alpha}\\
&&+ \left\|J_\omega\frac\nu{|\nabla\omega|}\circ A_\omega^{-1}\right\|_{0,\alpha}\left\|\frac1{J_\omega\frac1{|\nabla\omega|}\circ A_\omega^{-1}}\right\|_{n,\alpha}\\
&\leq&
C\bigg\{\|\nu\|_{n,\alpha}+\|\omega\|_{n+1,\alpha}\|\nu\|_{1,\alpha}
+\|\nu\|_{1,\alpha}(\|\omega\|_{n+1,\alpha}+1)\bigg\}\\
&\leq& C\left\{\|\nu\|_{n,\alpha}+\|\omega\|_{n+1,\alpha}\|\nu\|_{1,\alpha}\right\}
\end{eqnarray}
by (\ref{ineq:TameEstimatesForJ(omega)FractionOfA(omega)Inverse}) and (\ref{ineq:TameEstimatesForProductOfFunctions}).
\stopstep%

\nextpar{First derivative of $J_\omega u\circ A_\omega^{-1}$}
The operator $J_\omega u\circ A_\omega^{-1}$ is linear in $u$ so we only need
to worry about differentiability in $\omega$.
$J_{\omega_\epsilon}u\circ A_{\omega_\epsilon}^{-1}(\mu)$
is a smooth function of $(\mu,\epsilon)$,
and a similar argument as for $A_\omega^{-1}$ shows that 
$\frac{d}{d\epsilon}_{|\epsilon=0}(J_{\omega_\epsilon}u\circ A_{\omega_\epsilon}^{-1})$ exists
in the $C^\infty$-topology.
For $\mu\in (0,|\Omega|)$, setting $\lambda=A_\omega^{-1}(\mu)$,
and using (\ref{eq:ddepsilonJpsiepsilonu(lambda)}), (\ref{eq:ddlambdaJpsiu}), and (\ref{eq:DQ(psi)phi}),
\begin{eqnarray}
&&\frac{\partial}{\partial\epsilon}_{|\epsilon=0}(J_{\omega_\epsilon}u(A_{\omega_\epsilon}^{-1}(\mu)))
\quad= \quad D_\omega\left(J_\omega u\circ A_\omega^{-1}\right)\nu
\label{eq:ddepsilonJomegaepsilonuOfAomegaepsilonInverse}
\\
&=&
-J_\omega\left(\frac{\nu{\rm div}(u{N})}{|\nabla\omega|}\right)\circ A_\omega^{-1}(\mu)\\
&&+\left(J_\omega \left(\frac{{\rm div}(u{N})}{|\nabla\omega|}\right)\circ A_\omega^{-1}\right)(\mu)
\left(\frac{J_\omega\frac\nu{|\nabla\omega|}\circ A_\omega^{-1}}{J_\omega\frac1{|\nabla\omega|}\circ A_\omega^{-1}}\right)(\mu).
\end{eqnarray}
Lemma~\ref{lemma:QIsContinuous} implies that $(\omega,u)\mapsto J_\omega u\circ A_\omega^{-1}$ is continuously differentiable,
and Lemma~\ref{lemma:QIsTame} that the derivative is tame.
(From this expression, it is not too difficult to see, by an induction argument, that in fact $J_\omega u\circ A_\omega^{-1}$ is 
infinitely differentiable and all derivatives are tame.)
\stopstep %

\nextpar{Tame estimates on $D_\omega\left(J_\omega u\circ A_\omega^{-1}\right)\nu$}
We estimate $\left\|D_\omega\left(J_\omega u\circ A_\omega^{-1}\right)\nu\right\|_{n,\alpha}$ by the sum of three terms:
setting $N=\frac{\nabla\omega}{|\nabla\omega|}$,
(\ref{ineq:TameEstimatesForProductOfFunctions}) gives
\begin{eqnarray}
I+II+III&=&\left\|J_{\omega}\frac{\nu{\rm div}(uN)}{|\nabla\omega|}\circ A_{\omega}^{-1}\right\|_{n,\alpha}\\
&+&\left\|J_{\omega}\frac{{\rm div}(uN)}{|\nabla\omega|}\circ A_{\omega}^{-1}\right\|_{n,\alpha}\|DQ(\omega)\nu\|_{0,\alpha}\\
&+&\left\|J_{\omega}\frac{{\rm div}(uN)}{|\nabla\omega|}\circ A_{\omega}^{-1}\right\|_{0,\alpha}\|DQ(\omega)\nu\|_{n,\alpha}.
\end{eqnarray}

\noindent\textbf{Estimates on $I$}\qquad
By (\ref{ineq:TameEstimatesForJ(omega)FractionOfA(omega)Inverse}),
(\ref{ineq:TameEstimatesForProductOfFunctions}), and (\ref{ineq:TameEstimatesForNemitskiiOperator}),
we have for $\omega$ in a $\|\cdot\|_{2,\alpha}$-neighborhood, any $\nu$, $u$, and $n\geq 0$,
\begin{eqnarray}
&&\left\|J_\omega\frac{\nu{\rm div}(uN)}{|\nabla\omega|}\circ A_\omega^{-1}\right\|_{n,\alpha}\\
&\leq&C\bigg\{\|\nu{\rm div}(uN)\|_{n,\alpha}+\|\omega\|_{n+1,\alpha}\|\nu{\rm div}(uN)\|_{1,\alpha}\bigg\}\\
&\leq&C\bigg\{\|\nu\|_{n,\alpha}\|u\|_{1,\alpha}+\|uN\|_{n+1,\alpha}\|\nu\|_{0,\alpha}
+\|\omega\|_{n+1,\alpha}\|\nu\|_{1,\alpha}\|u\|_{2,\alpha}\bigg\}\\
&\leq&C\bigg\{\|\nu\|_{n,\alpha}\|u\|_{1,\alpha}+\|u\|_{n+1,\alpha}\|\nu\|_{0,\alpha}
+\|N\|_{n+1,\alpha}\|u\|_{0,\alpha}\|\nu\|_{1,\alpha}\\
&&+\|\omega\|_{n+1,\alpha}\|\nu\|_{1,\alpha}\|u\|_{2,\alpha}\bigg\}\\
&\leq&C\bigg\{\|\nu\|_{n,\alpha}\|u\|_{1,\alpha}+\|u\|_{n+1,\alpha}\|\nu\|_{0,\alpha}
+(\|\omega\|_{n+2,\alpha}+1)\|u\|_{2,\alpha}\|\nu\|_{1,\alpha}\bigg\}\\
&\leq&C\bigg\{\|\nu\|_{n,\alpha}\|u\|_{1,\alpha}+\|u\|_{n+1,\alpha}\|\nu\|_{0,\alpha}
+\|\omega\|_{n+2,\alpha}\|u\|_{2,\alpha}\|\nu\|_{1,\alpha}\bigg\}\\
&\leq&C\bigg\{\|u\|_{1,\alpha}\|\nu\|_{n,\alpha}+\|u\|_{n+1,\alpha}\|\nu\|_{0,\alpha}+\|\omega\|_{n+2,\alpha}\|u\|_{2,\alpha}\|\nu\|_{1,\alpha}
\bigg\}.
\end{eqnarray}

\noindent\textbf{Estimates on $J_{\omega}\frac{{\rm div}(uN)}{|\nabla\omega|}\circ A_{\omega}^{-1}$}\qquad
By (\ref{ineq:TameEstimatesForJ(omega)FractionOfA(omega)Inverse}),
(\ref{ineq:TameEstimatesForProductOfFunctions}), and (\ref{ineq:TameEstimatesForNemitskiiOperator}),
we have for $n\geq 0$
\begin{eqnarray}
&&\left\|J_{\omega}\frac{{\rm div}(uN)}{|\nabla\omega|}\circ A_{\omega}^{-1}\right\|_{n,\alpha}\\
&\leq& C\left\{\|uN\|_{n+1,\alpha}+\|\omega\|_{n+1,\alpha}\|uN\|_{2,\alpha}\right\}\\
&\leq& C\left\{\|u\|_{n+1,\alpha}+\|u\|_{0,\alpha}\|N\|_{n+1,\alpha}+\|\omega\|_{n+1,\alpha}\|u\|_{2,\alpha}(\|\omega\|_{3,\alpha}+1)\right\}\\
&\leq& C\left\{\|u\|_{n+1,\alpha}+\|u\|_{0,\alpha}(\|\omega\|_{n+2,\alpha}+1)+\|\omega\|_{n+2,\alpha}\|u\|_{2,\alpha}\right\}\\
&\leq& C\left\{\|u\|_{n+1,\alpha}+\|\omega\|_{n+2,\alpha}\|u\|_{2,\alpha}\right\}
\end{eqnarray}
provided $\omega$ remains in a $\|\cdot\|_{2,\alpha}$-neighborhood.
(We have used that, by interpolation inequalities,
$\|\omega\|_{n+1,\alpha}\|\omega\|_{3,\alpha}\leq C\|\omega\|_{n+2,\alpha}\|\omega\|_{2,\alpha}\leq C\|\omega\|_{n+2,\alpha}$
for $\omega$ in a $\|\cdot\|_{2,\alpha}$-neighborhood.)
\bigskip%

\noindent\textbf{Estimates on II+III}\qquad
For $\omega$ in a $\|\cdot\|_{2,\alpha}$-neighborhood, and any $u$, $\nu$,
we conclude from the above and (\ref{ineq:TameEstimatesForProductOfFunctions}) that for $n\geq 0$
\begin{eqnarray}
&&II+III\\
&\leq& C\Bigg\{\bigg(\|u\|_{n+1,\alpha}+\|\omega\|_{n+2,\alpha}\|u\|_{2,\alpha}\bigg)\|\nu\|_{1,\alpha}\\
&&+\|u\|_{2,\alpha}\bigg(\|\nu\|_{n,\alpha}+\|\omega\|_{n+1,\alpha}\|\nu\|_{1,\alpha}
\bigg)\Bigg\}\\
&\leq&C\bigg\{
\|u\|_{n+1,\alpha}\|\nu\|_{1,\alpha}+\|u\|_{2,\alpha}\|\nu\|_{n,\alpha}+\|\omega\|_{n+2,\alpha}\|u\|_{2,\alpha}\|\nu\|_{1,\alpha}
\bigg\}
\end{eqnarray}

\noindent\textbf{Conclusion}\qquad
Putting the above together, 
\begin{eqnarray}
&&\left\|D_\omega(J_\omega u\circ A_\omega^{-1})\nu\right\|_{n,\alpha}\\
&\leq& C\bigg\{\|u\|_{n+1,\alpha}\|\nu\|_{1,\alpha}+\|u\|_{2,\alpha}\|\nu\|_{n,\alpha}+\|\omega\|_{n+2,\alpha}\|u\|_{2,\alpha}\|\nu\|_{1,\alpha}\bigg\}.
\end{eqnarray}
\stopstep%

\nextpar{Second derivatives of $Q$}
From (\ref{eq:DQ(omega)nu}) it is clear that $DQ(\omega)\nu$ is continuously differentiable in $\omega$.
Using (\ref{eq:ddepsilonJomegaepsilonuOfAomegaepsilonInverse}) and after some simplification,
the second derivative $D^2Q(\omega)(\nu_1,\nu_2)$, being the partial derivative of $DQ(\omega)\nu_1$ 
with respect to $\omega$ in the direction $\nu_2$, is then given by
\begin{eqnarray}
&&D^2Q(\omega)(\nu_1,\nu_2)\\
&=&
\frac{\left(J_\omega\frac{\nu_1}{|\nabla\omega|}\circ A_\omega^{-1}\right)\left(J_\omega\frac{{\rm div}\left(\frac{\nu_2N}{|\nabla\omega|}\right)}{|\nabla\omega|}\circ A_\omega^{-1}\right)}{\left(J_\omega\frac1{|\nabla\omega|}\circ A_\omega^{-1}\right)^2}\\
&+&\frac{\left(J_\omega\frac{\nu_2}{|\nabla\omega|}\circ A_\omega^{-1}\right)\left(J_\omega\frac{{\rm div}\left(\frac{\nu_1N}{|\nabla\omega|}\right)}{|\nabla\omega|}\circ A_\omega^{-1}\right)}{\left(J_\omega\frac1{|\nabla\omega|}\circ A_\omega^{-1}\right)^2}\\
&-&\frac{J_\omega \frac{{\rm div}\left(\frac N{|\nabla\omega|}\right)}{|\nabla\omega|}\circ A_\omega^{-1}}{\left(J_\omega\frac1{|\nabla\omega|}\circ A_\omega^{-1}\right)^3}
\left(J_\omega \frac{\nu_1}{|\nabla\omega|}\circ A_\omega^{-1}\right)\left( J_\omega \frac{\nu_2}{|\nabla\omega|}\circ A_\omega^{-1}\right)\\
&-&\frac{J_\omega \frac{{\rm div}\left(\frac{\nu_1\nu_2N}{|\nabla\omega|}\right)}{|\nabla\omega|}\circ A_\omega^{-1}}{J_\omega \frac1{|\nabla\omega|}\circ A_\omega^{-1}}.
\end{eqnarray}
(One verifies that this expression is symmetric in $\nu_1$ and $\nu_2$.)
By Lemma~\ref{lemma:QIsContinuous} and Lemma~\ref{lemma:QIsTame}, 
a moment's concentration shows that this is a continuous map
\begin{equation}
D^2Q\colon\left\{
\begin{array}{ccccc}
\mathcal F_+^m&\times&C^{m-1}_{\overline\Omega}&\times&C^{m-1}_{\overline\Omega}\\\omega&&\nu_1&&\nu_2\end{array}
\longrightarrow\begin{array}{c}C^{m-2}_{[0,|\Omega]}\\D^2Q(\omega)(\nu_1,\nu_2)
\end{array}\right\}
\end{equation}
and that it is tame.
\stopstep%

\nextpar{Tame estimates on $D^2Q(\omega)(\nu_1,\nu_2)$}
We write the above expression as $D^2Q(\omega)(\nu_1,\nu_2)=I_2+II_2+III_2+IV_2$.
Note that all the factors not depending on $\nu_1$ nor $\nu_2$ have their $\|\cdot\|_{n,\alpha}$-norms
bounded by $C(\|\omega\|_{n+1,\alpha}+1)$.

\noindent\textbf{Estimates on $I_2$}\\
Using (\ref{ineq:||JomegauOfAomegaInverse||0,alpha}) for $\|\cdot\|_{0,\alpha}$-estimates,
\begin{eqnarray}
\|I_2\|_{n,\alpha}
&\leq&C\Bigg\{(\|\omega\|_{n+1,\alpha}+1)\left\|\frac{\nu_1}{|\nabla\omega|}\right\|_{0,\alpha}\left\|\frac{{\rm div}\frac{\nu_2N}{|\nabla\omega|}}{|\nabla\omega|}\right\|_{0,\alpha}\\
&&+(\|\nu_1\|_{n,\alpha}+\|\omega\|_{n+1,\alpha}\|\nu_1\|_{1,\alpha})\|\nu_2\|_{1,\alpha}\\
&&+\|\nu_1\|_{0,\alpha}\left(\left\|{\rm div}\left(\frac{\nu_2N}{|\nabla\omega|}\right)\right\|_{n,\alpha}
+\|\omega\|_{n+1,\alpha}\left\|{\rm div}\left(\frac{\nu_2N}{|\nabla\omega|}\right)\right\|_{1,\alpha}\right)\Bigg\}\\
&\leq&C\Bigg\{(\|\omega\|_{n+1,\alpha}+1)\|\nu_1\|_{0,\alpha}\|\nu_2\|_{1,\alpha}\\
&&+(\|\nu_1\|_{n,\alpha}+\|\omega\|_{n+1,\alpha}\|\nu_1\|_{1,\alpha})\|\nu_2\|_{1,\alpha}\\
&&\|\nu_1\|_{_0,\alpha}
\left(\left\|\frac{\nu_2N}{|\nabla\omega|}\right\|_{n+1,\alpha}+\|\omega\|_{n+1,\alpha}\left\|\frac{\nu_2N}{|\nabla\omega|}\right\|_{2,\alpha}
\right)
\Bigg\}\\
&\leq&C\Bigg\{
\|\omega\|_{n+1,\alpha}\|\nu_1\|_{1,\alpha}\|\nu_2\|_{1,\alpha}+\|\nu_1\|_{n,\alpha}\|\nu_2\|_{1,\alpha}\\
&&+\|\nu_1\|_{0,\alpha}
\Big(\|\nu_2\|_{n+1,\alpha}+(\|\omega\|_{n+2,\alpha}+1)\|\nu_2\|_{0,\alpha}\\
&&+\|\omega\|_{n+1,\alpha}\|\nu_2\|_{2,\alpha}(\|\omega\|_{3,\alpha}+1)
\Big)
\Bigg\}\\
&\leq&C\bigg\{\|\nu_1\|_{n,\alpha}\|\nu_2\|_{1,\alpha}+\|\nu_1\|_{0,\alpha}\|\nu_2\|_{n+1,\alpha}
+\|\omega\|_{n+2,\alpha}\|\nu_1\|_{1,\alpha}\|\nu_2\|_{2,\alpha}\bigg\}
\end{eqnarray}
where we have used interpolations to get
$\|\omega\|_{n+1,\alpha}\|\omega\|_{3,\alpha}\leq C\|\omega\|_{n+2,\alpha}\|\omega\|_{2,\alpha}\leq C\|\omega\|_{n+2,\alpha}$.
\stopstep%

\noindent\textbf{Estimates on $I_2+II_2$}\\
Since $II_2$ is obtained by interchanging $\nu_1$ and $\nu_2$, we immediately have
\begin{equation}
\|I_2+II_2\|_{n,\alpha}
\leq C\bigg\{\|\nu_1\|_{n+1,\alpha}\|\nu_2\|_{1,\alpha}+\|\nu_1\|_{1,\alpha}\|\nu_2\|_{n+1,\alpha}
+\|\omega\|_{n+2,\alpha}\|\nu_1\|_{2,\alpha}\|\nu_2\|_{2,\alpha}\bigg\}.
\end{equation}
\stopstep%

\noindent\textbf{Estimates on $III_2$}\\
We have
\begin{eqnarray}
&&\|III_2\|_{n,\alpha}\\
&\leq&C\Bigg\{
(\|\omega\|_{n+1,\alpha}+1)\left\|\nu_1\right\|_{0,\alpha}\left\|\nu_2\right\|_{0,\alpha}\\
&&+\bigg(\left\|\nu_1\right\|_{n,\alpha}+\|\omega\|_{n+1,\alpha}\|\nu_1\|_{1,\alpha}\bigg)\left\|\nu_2\right\|_{0,\alpha}\\
&&+\|\nu_1\|_{0,\alpha}\bigg(\left\|\nu_2\right\|_{n,\alpha}+\|\omega\|_{n+1,\alpha}\|\nu_2\|_{1,\alpha}\bigg)
\Bigg\}\\
&\leq&C\Bigg\{
\|\nu_1\|_{n,\alpha}\|\nu_2\|_{0,\alpha}+\|\nu_1\|_{0,\alpha}\|\nu_2\|_{n,\alpha}+\|\omega\|_{n+1,\alpha}\|\nu\|_{1,\alpha}\|\nu_2\|_{1,\alpha}
\Bigg\}
\end{eqnarray}
\stopstep%

\noindent\textbf{Estimates on $IV_2$}\\
We have
\begin{eqnarray}
&&\|IV\|_{n,\alpha}\\
&\leq&C\Bigg\{
(\|\omega\|_{n+1,\alpha}+1)\left\|\frac{\nu_1\nu_2N}{|\nabla\omega|}\right\|_{1,\alpha}
+\left\|\frac{\nu_1\nu_2N}{|\nabla\omega|}\right\|_{n+1,\alpha}
\Bigg\}\\
&\leq&C\Bigg\{
(\|\omega\|_{n+1,\alpha}+1)\|\nu_1\|_{1,\alpha}\|\nu_2\|_{1,\alpha}\\
&+&\|\nu_1\|_{n+1,\alpha}\|\nu_2\|_{0,\alpha}+\|\nu_1\|_{0,\alpha}\|\nu_2\|_{n+1,\alpha}
+(\|\omega\|_{n+2,\alpha}+1)\|\nu_1\|_{1,\alpha}\|\nu_2\|_{1,\alpha}\Bigg\}\\
&\leq&C\Bigg\{\|\nu_1\|_{n+1,\alpha}\|\nu_2\|_{1,\alpha}+\|\nu_1\|_{1,\alpha}\|\nu_2\|_{n+1,\alpha}
+\|\omega\|_{n+2,\alpha}\|\nu_1\|_{1,\alpha}\|\nu_2\|_{1,\alpha}\Bigg\}.
\end{eqnarray}
\stopstep%

\noindent\textbf{Conclusion}\\
Putting the estimates on $I_2, II_2, III_2, IV_2$ together,
\begin{equation}
\|D^2Q(\omega)(\nu_1,\nu_2)\|_{n,\alpha}
\leq C\Bigg\{\|\nu_1\|_{n+1,\alpha}\|\nu_2\|_{1,\alpha}+\|\nu_1\|_{1,\alpha}\|\nu_2\|_{n+1,\alpha}
+\|\omega\|_{n+2,\alpha}\|\nu_1\|_{2,\alpha}\|\nu_2\|_{2,\alpha}\Bigg\}.
\end{equation}
\stopproof %

\section{Proof of Theorem~\ref{theorem:MainTheorem}}\label{section:ProofOfMainTheorem}
\noindent\textbf{Outline of proof}\\
That $T$ is smooth tame is immediate since $T(F)=(Q\circ \Delta \circ S)(F)$ is a composition of smooth tame maps.
The crucial part in the surjective part of the Nash-Moser theorem is to establish that $DT(F)f$ has a tame family of right-inverses.
We emphasize that the non-degeneracy condition (ND2) is only made 
\textit{at} the reference steady-state, and 
\textit{not in an entire neighborhood} of the reference steady-state.
The problem of finding such right-inverse for $DT(F)f$, given in (\ref{eq:DT(F)fReminiscentOfFredholm}), 
is equivalent to inverting a map
of the form $g+K(F)g=h$ where $K(F)g$ can be thought of as a ``compact perturbation'' of the first term $g$,
see (\ref{eq:K(F)=Ktilde(F)VB(F)}).
This is precisely what was done in Lemma~\ref{proposition:VE(c)kIsSmoothTame} 
with the elliptic operator $\Delta \phi+c\phi=k$ (augmented with suitable boundary conditions):
$c\phi$ is a ``compact perturbation'' of $\Delta\phi$.
There, the estimates on the (bilinear) term $c\phi$ were standard.
Here, the term $K(F)g$ is more complicated and requires considerably more work.

The injective part of Theorem~\ref{theorem:MainTheorem} is proved in Section~\ref{section:InjectivePart}.
At the conceptual level, the proof is an adjustment of the injective part of the Nash-Moser
Inverse Function Theorem as presented in Section~1.3, Part~III \cite{HamiltonIFTNashMoser}.
(One cannot use this theorem directly because of complications created by the lack of injectivity of the map $F\mapsto \psi$.)
\bigskip%

\noindent\textbf{Assumptions}\\
We recall the main assumptions.

The domain $\Omega$ is assumed diffeomorphic to the annulus so that
\begin{equation}\partial\Omega = \Gamma_o\cup\Gamma_i.\end{equation}
We assume that the reference steady-state $\overline\omega=\overline F(\overline\psi)$ is such that $\overline F'\neq 0$,
$\overline \omega$ has no critical points,
and satisfies the non-degeneracy conditions (ND1) and (ND2).
By continuity of $F\in C^1\mapsto \psi\in C^2$, make $\epsilon_S$ in $\mathcal V_S(\overline F)$ smaller if necessary
(see Proposition~\ref{proposition:SolutionMappsi=S(F)}), so that the corresponding $\psi$ has no critical points either
and that $F'\neq 0$.
We will then assume without loss of generality that
\begin{equation}\psi\leq 0.\end{equation}
The interval $I$ introduced in Section~\ref{section:ConstructSolutionOperator} can now be taken of the form
\begin{equation}I=[\underline c,0]\qquad\text{where}\quad \underline c<\min\overline\psi\quad \text{is fixed}.\label{eq:IntervalI}
\end{equation}
For simplicity, the calculations will be performed assuming that
\begin{equation}F'>0.\end{equation}
As observed in the Introduction, the case $F'>0$ is special in that the corresponding solution
automatically satisfies both non-degeneracy conditions (ND1) and (ND2).
But this property will never be used in the following.
\bigskip%

\noindent\textbf{The first derivative $DT(F)f$}\\
The map
\begin{equation}
T\colon\left\{\begin{array}{ccc}
\mathcal V_S^n(\overline F)\\F\end{array}
\longrightarrow\begin{array}{c}C^n_{[0,|\Omega|]}\\A_\omega^{-1}
\end{array}\right\}\qquad (n\geq 2)
\end{equation}
is continuous as a map of Banach spaces
(Lemma~\ref{lemma:QIsContinuous} 
and (\ref{eq:ContinuityOfF->omegaInC^n-Grading}) after Proposition~\ref{proposition:SolutionMappsi=S(F)}).
Write the first derivative as
$DT(F)f=DQ(\omega)\nu$
where  $\omega=\Delta\psi$ and $\psi=S(F)$ solves the steady-state equation
\begin{equation}\Delta\psi=F(\psi),\qquad \psi_{|\Gamma_o}=0,\quad \frac{\partial\psi}{\partial\tau}_{|\Gamma_i}=0,\quad 
\int_{\Gamma_i}\frac{\partial\psi}{\partial N}=\gamma_i,
\end{equation}
and $\nu=\Delta \phi$ where $\phi=DS(F)f$ solves the linearized steady-state equation
\begin{equation}
\Delta\phi=F'(\psi)\phi+f(\psi),\qquad \phi_{|\Gamma_o}=0,\quad \frac{\partial\phi}{\partial\tau}_{|\Gamma_i}=0,\quad \int_{\Gamma_i}\frac{\partial\phi}{\partial N}=0.
\end{equation}
Thus,
\begin{equation}
DT\colon\left\{\begin{array}{ccc}
\mathcal V_S^n(\overline F)&\times&C^{n-1}_I\\F&&f\end{array}
\longrightarrow\begin{array}{c}C^{n-1}_{[0,|\Omega|]}\\DT(F)f
\end{array}\right\}\qquad (n\geq 2)
\end{equation}
is continuous as a map of Banach spaces
(Lemma~\ref{lemma:QIsSmoothTame} and (\ref{eq:ContinuityOf(F,f)->nuInC^n-Grading})
after Proposition~\ref{proposition:SolutionMappsi=S(F)}).
From Proposition~\ref{proposition:AomegaInverseIsSmoothTame},
\begin{equation}
DT(F)f=\frac{J_\omega\frac\nu{|\nabla\omega|}\circ A_\omega^{-1}}{J_\omega\frac1{|\nabla\omega|}\circ A_\omega^{-1}}
=\frac{J_\omega\frac{f(\psi)}{|\nabla\omega|}\circ A_\omega^{-1}}{J_\omega\frac1{|\nabla\omega|}\circ A_\omega^{-1}}
+ \frac{J_\omega\frac{F'(\psi)\phi}{|\nabla\omega|}\circ A_\omega^{-1}}{J_\omega\frac1{|\nabla\omega|}\circ A_\omega^{-1}}
\end{equation}
To simplify this, use the identity (when $F'>0$)
\begin{equation}
T(F)=F\circ A_\psi^{-1}\label{eq:UsefulFactsT(F)=FOfApsiInverse}
\end{equation}
which follows from
$A_\psi(\lambda)=|\{\psi<\lambda\}|=|\{\omega<F(\lambda)\}|=A_\omega(F(\lambda))$.
(A similar expression holds in the case $F'<0$.)
Then, $\omega(x)=A_\omega^{-1}(\lambda)$ if and only if $\psi(x)=A_\psi^{-1}(\lambda)$
and
\begin{equation}
\int_{\omega=A_\omega^{-1}(\lambda)}\frac{f(\psi)}{|\nabla\omega|}dl
=f(A_\psi^{-1}(\lambda))\int_{\omega=A_\omega^{-1}(\lambda)}\frac1{|\nabla\omega|}dl.
\end{equation}
Also, $\nabla\omega=F'(\psi)\nabla\psi$ so that
after some elementary calculations we obtain
\begin{eqnarray}
DT(F)\cdot f
&=&f\circ A_\psi^{-1}+(F'\circ A_\psi^{-1})\left(\frac{J_\psi\frac\phi{|\nabla\psi|}\circ A_\psi^{-1}}{J_\psi\frac1{|\nabla\psi|}\circ A_\psi^{-1}}\right)\\
&=&f\circ A_\psi^{-1} +\left(\frac{dA_\omega^{-1}}{d\mu}\right)\left(J_\psi\frac\phi{|\nabla\psi|}\circ A_\psi^{-1}\right)\\
&=:&B(F)f+\tilde K(F)f.
\label{eq:DT(F)fReminiscentOfFredholm}
\end{eqnarray}
\bigskip%

\noindent\textbf{Construction of a tame right-inverse $f=L(F)h$ to $h=DT(F)\cdot f$}\\
In order to construct a right-inverse $f=L(F)h$ to $h=DT(F)f$ we will construct first the \foremph{inverse}
to a modification $h=M(F)g$ of $h=DT(F)f$.
Set
\begin{equation}
M(F):=DT(F)\cdot VB(F)
=\text{Id}_{C^\infty_{[0,|\Omega|]}}+K(F)\end{equation}
where
\begin{equation}
K(F)\cdot g:=\tilde K(F)\cdot VB(F)\cdot g=\left(\frac{d A_\omega^{-1}}{d\mu}\right)\left(J_\psi\frac\phi{|\nabla\psi|}\circ A_\psi^{-1}\right),
\label{eq:K(F)=Ktilde(F)VB(F)}\end{equation}
and $\Delta \phi=F'(\psi)\phi+f(\psi)$, $f=VB(F)\cdot g$,
and $f=VB(F)g$ is a tame right-inverse to $g=B(F)f$ (intuitively, it is ``$f=g\circ A_\psi$'').
The latter is constructed in Lemma~\ref{lemma:B(F)fHasASmoothTameFamilyOfRightInverses}.
We will show that $h=M(F)g$ has a tame inverse $g=VM(F)h$, and thus 
that $h=DT(F)f$ has a tame right-inverse by setting
\begin{equation}L(F)\cdot h:=VB(F)\cdot VM(F)\cdot h.\end{equation}
This is possible because $K(F)g$ can be viewed as a ``compact perturbation'' of $g$.
This is due to the fact that $\phi$ gains sufficient regularity from $f$.
In fact, the proof is completely analogous to the proof that $\Delta \phi+c\phi=k$ has 
a family of tame inverses:
compare Lemma~\ref{lemma:EstimatesForLinearEllipticEquations}
with Lemma~\ref{lemma:EstimatesOnh=g+K(F)g},
and Proposition~\ref{proposition:VE(c)kIsSmoothTame} with Proposition~\ref{proposition:TameInverseToh=g+K(F)g}.
\bigskip%

\subsection{Right-inverse to $B(F)f=f\circ A_\psi^{-1}$}\label{section:RightInverseToB(F)f}
First, as an auxiliary step, we need to construct a right-inverse $f=VB(F)g$ to $g=B(F)f$.
Since $B(F)$ is surjective for each $F$, we know that a right-inverse exists 
for each $F$.
However, we need the inverse $f=VB(F)g$ to $g=B(F)f$ to be continuous in both $F$ and $h$,
and to satisfy tame estimates.
Naively, the inverse of $g=f\circ A_\psi^{-1}$ should be ``$f=g\circ A_\psi$'',
but $f$ is defined on an interval larger than the domain of $A_\psi$.
\begin{lemma}\label{lemma:B(F)fHasASmoothTameFamilyOfRightInverses}
The map $B(F)\cdot f$ has a smooth tame family of right-inverses  $f=VB(F)\cdot g$
defined on a sufficiently small $\|\cdot\|_{1,\alpha}$-neighborhood $\mathcal V_B(\overline F)\subset \mathcal V_S(\overline F)$ of $\overline F$:
\begin{equation}
VB(F)\colon C^\infty_{[0,|\Omega|]}\rightarrow C^\infty_{[\underline c,0]},\qquad B(F)\cdot VB(F) = \text{Id}_{C^\infty_{[0,|\Omega|]}}.
\end{equation}
For $n\geq 2$,  $F\in\mathcal V_B(\overline F)$ and any $g$ (without restriction),
\begin{equation}\|f\|_{n,\alpha}\leq C\cdot \bigg\{\|g\|_{n,\alpha}+\|F\|_{n-2,\alpha}\|g\|_{1,\alpha}\bigg\}
\end{equation}
while $\|f\|_{0,\alpha}\leq C\|g\|_{0,\alpha}$, $\|f\|_{1,\alpha}\leq C\|g\|_{1,\alpha}$.

\end{lemma}
\startproof %
Fix $D>0$ and let $\overline B\in C^\infty_{[-D,|\Omega|+D]}$ be a monotone increasing extension 
of $A_{\overline\psi}^{-1}\in C^\infty_{[0,|\Omega|]}$.
It can be arranged so that $\text{range}(\overline B)=[\underline c,|\underline c|]$ (see (\ref{eq:IntervalI})).
See proof of Corollary~1.3.7, p.~138, Part~II of \cite{HamiltonIFTNashMoser}.
Let 
\begin{equation}
\mathcal E_0\colon C^\infty_{[0,|\Omega|]}\quad\longrightarrow\quad \{b\in C^\infty_{[-D,|\Omega|+D]}~|~b(-D)=b(|\Omega|+D)=0\}
\end{equation}
be an extension operator taking functions on $[0,|\Omega|]$ to functions on $[-D,|\Omega|+D]$ vanishing at the endpoints.
(The target space is easily seen to be a tame Fr\'echet space).
It can be made tame linear of degree $0$: again from the proof of Corollary~II.1.3.7, p.~138, [Hamilton]),
extend $b=b(\lambda)\in C^\infty_{[0,|\Omega|]}$ to $\lambda\leq 0$, 
then to $\lambda\geq |\Omega|$,
and finally multiply by a smooth cut-off function with support in $(-D,|\Omega|+D)$ and equal to $1$ on $[0,|\Omega|]$.
Then we have the tame estimates for $n\geq 0$
\begin{equation}\|\mathcal E_0b\|_{n,\alpha}\leq C\|b\|_{n,\alpha}.\label{ineq:TameEstimatesForExtensionOperator}
\end{equation}
Define now 
\begin{equation}
\mathcal E\colon
\left\{\begin{array}c C^\infty_{[0,|\Omega|]}\\B\end{array}\quad \longrightarrow\quad 
\begin{array}c C^\infty_{[-D,|\Omega|+D]}\\\overline B+\mathcal E_0(B-A_{\overline\psi}^{-1})
\end{array}\right\}
\end{equation}
which extends maps defined on $[0,|\Omega|]$ to maps defined on $[-D,|\Omega|]$ 
with fixed endpoint values $\underline c$ and $|\underline c|$ at $-D$ and $|\Omega|+D$ respectively.
$\mathcal E$ is smooth tame since it is affine with tame linear part $\mathcal E_0$.
For a sufficiently small $\|\cdot\|_1$-neighborhood $\mathcal V(A_{\overline\psi}^{-1})$ of $A_{\overline\psi}^{-1}$,
it also defines a map
\begin{equation}
\mathcal E\colon (\mathcal V(A_{\overline\psi}^{-1})\subset C^\infty_{[0,|\Omega|]})
\quad\longrightarrow\quad\mathcal D^\infty_{I_1,I_2}
\end{equation}
where $\mathcal D^\infty_{I_1,I}$ is the set of smooth diffeomorphisms from $I_1=[-D,|\Omega|]$ to $I_2=[\underline c,|\underline c|]$.

From $g=f\circ A_\psi^{-1}$ we find $\mathcal E(g)=\mathcal E(f\circ A_\psi^{-1})=f\circ \mathcal E(A_\psi^{-1})$.
Let then
\begin{equation}V\colon
\left\{\begin{array}c \mathcal D^\infty_{I_1,I_2}\\B
\end{array}\quad \longrightarrow\quad
\begin{array}c \mathcal D^\infty_{I_2,I_1}\\B^{-1}\end{array}
\right\}
\end{equation}
denote the operator which takes inverses.
Choose now
\begin{equation}\mathcal V_B(\overline F)\subset \mathcal V_S(\overline F)
\end{equation}
a sufficiently small $\|\cdot\|_{1,\alpha}$-neighborhood of 
$\mathcal V_S(\overline F)$ from Proposition~\ref{proposition:SolutionMappsi=S(F)}
so that the corresponding $A_\psi^{-1}$ remains in $\mathcal V(A_{\overline\psi}^{-1})$
(use also Lemma~\ref{lemma:QIsContinuous}).
We have constructed a right-inverse
\begin{equation}
VB(F)\cdot g:= \left.\left(\mathcal E g\circ V(\mathcal  E(A_\psi^{-1}))\right)\right|_{[\underline c,0]}
\label{eq:VB(F)g}
\end{equation}
defined for any $F\in\mathcal V_B(\overline F)$ and any $g\in C^\infty_{[0,|\Omega|]}$.

Clearly, $f=VB(F)g$ is smooth tame.
Using the above estimates on $\mathcal E_0$, 
on $A_\psi^{-1}$ from Proposition~\ref{proposition:AomegaInverseIsSmoothTame},
and on $\psi=S(F)$ from Proposition~\ref{proposition:SolutionMappsi=S(F)}, we have
for $n\geq 2$
\begin{equation}
\|\mathcal E(A_\psi^{-1})\|_{n,\alpha}
\leq C\left(\|A_\psi^{-1}\|_{n,\alpha}+1\right)
\leq C\left(\|\psi\|_{n,\alpha}+1\right)
\leq C\left(\|F\|_{n-2,\alpha}+1\right).
\end{equation}
From tame estimates on composition of functions from Lemma~\ref{lemma:CompositionIsSmoothTame}
and on the inversion operator $V$ from Lemma~\ref{lemma:InversionIsSmoothTame} of the Appendix,
for $F\in\mathcal V_B(\overline F)$ and any $g\in C^\infty_{[0,|\Omega|}$,
\begin{eqnarray}
\|f\|_{n,\alpha}
&\leq&C\left(\|\mathcal E(g)\|_{n,\alpha}+\|\mathcal E(A_\psi^{-1})\|_{n,\alpha}\|\mathcal E(g)\|_{1,\alpha}\right)\\
&\leq&C\left(\|g\|_{n,\alpha}+(\|F\|_{n-2,\alpha}+1)\|\|g\|_{1,\alpha}\right)\\
&\leq&C\left(\|g\|_{n,\alpha}+\|F\|_{n-2,\alpha}\|\|g\|_{1,\alpha}\right).
\end{eqnarray}
Note that $\|\mathcal E(A_\psi^{-1})\|_{1,\alpha}$ remains bounded for $F\in \mathcal V_B(\overline F)$.
Thus, we deduce easily the desired estimates on $\|f\|_{0,\alpha}$ and $\|f\|_{1,\alpha}$.
\stopproof%
 \setcounter{counterpar}{0}

\subsection{Summary of tame estimates}
Here we collect tame estimates which will be used abundantly in the next Sections.
Some estimates will be given in two equivalent forms, the second being particularly useful for the estimates on 
the difference $K(F)g-K(\overline F)g$.
\bigskip%

\noindent\textbf{From Proposition~\ref{proposition:SolutionMappsi=S(F)}:}
If $\omega=\Delta\psi=F(\psi)$ with $F\in\mathcal V_S(\overline F)$, then for $n\geq 0$
\begin{equation}
\|\omega\|_{n,\alpha}\quad \text{and}\quad \|\psi\|_{n+2,\alpha}
~\leq~ C\Big\{\|F\|_{n,\alpha}+1\Big\}~\leq~ C\Big\{\|F-\overline F\|_{n,\alpha}+1\Big\}
\label{ineq:TEomegaAndpsi}
\end{equation}
from the triangle inequality 
$\|F\|_{n,\alpha}\leq \|\overline F\|_{n,\alpha}+\|F-\overline F\|_{n,\alpha}\leq C\Big(\|F-\overline F\|_{n,\alpha}+1\Big)$.\\
If $\Delta\phi=F'(\psi)\phi+f(\psi)$ with $F\in\mathcal V_S(\overline F)$ and $f\in C^\infty_{\overline\Omega}$,
then for $n\geq 0$
\begin{eqnarray}
\|\phi\|_{n+2,\alpha}&\leq& C\Big\{\|f\|_{n,\alpha}+\|F\|_{n+1,\alpha}\|f\|_{1,\alpha}\Big\},\label{ineq:TEphi=phi(F,f)}\\
\|\phi\|_{n+2,\alpha}&\leq& C\Big\{\|f\|_{n,\alpha}+\|F-\overline F\|_{n+1,\alpha}\|f\|_{1,\alpha}\Big\}\label{ineq:ATEphi=phi(F,f)}
\end{eqnarray}
again by the triangle inequality and using $\|f\|_{1,\alpha}\leq \|f\|_{n,\alpha}$
for $n\geq 1$
(for $n=0$, the last term in (\ref{ineq:TEphi=phi(F,f)}) and
(\ref{ineq:ATEphi=phi(F,f)}) is actually not needed, see Proposition~\ref{proposition:SolutionMappsi=S(F)}).

\noindent\textbf{From Lemma~\ref{lemma:B(F)fHasASmoothTameFamilyOfRightInverses}}
If $f=VB(F)g$ where $F\in\mathcal V_B(\overline F)\subset \mathcal V_S(\overline F)$ and $g\in C^\infty_{[0,|\Omega|]}$,
then for $n\geq 2$,
\begin{eqnarray}
\|f\|_{n,\alpha}&\leq& C\Big\{\|g\|_{n,\alpha}+\|F\|_{n-2,\alpha}\|g\|_{1,\alpha}\Big\},\label{ineq:TEf=VB(F)g}\\
\|f\|_{n,\alpha}&\leq& C\Big\{\|g\|_{n,\alpha}+\|F-\overline F\|_{n-2,\alpha}\|g\|_{1,\alpha}\Big\}\label{ineq:ATEf=VB(F)g}
\end{eqnarray}
while $\|f\|_{0,\alpha}\leq C\|g\|_{0,\alpha}$, $\|f\|_{1,\alpha}\leq C\|g\|_{1,\alpha}$.

Combining the above, we have for $n\geq 0$, 
\begin{eqnarray}
\|\phi\|_{n+2,\alpha}&\leq& C\Big\{\|g\|_{n,\alpha}+\|F\|_{n+1,\alpha}\|g\|_{1,\alpha}\Big\},\label{ineq:TEphi=phi(F,g)}\\
\|\phi\|_{n+2,\alpha}&\leq& C\Big\{\|g\|_{n,\alpha}+\|F-\overline F\|_{n+1,\alpha}\|g\|_{1,\alpha}\Big\}.\label{ineq:ATEphi=phi(F,g)}
\end{eqnarray}

\noindent\textbf{From Lemma~\ref{lemma:QIsSmoothTame}:}
For $\psi\in\mathcal F_+$ in a $\|\cdot\|_{2,\alpha}$-neighborhod, $\phi\in C^\infty_{\overline\Omega}$, and $n\geq 0$,
\begin{eqnarray}
\left\|J_\psi\frac\phi{|\nabla\psi|}\circ A_\psi^{-1}\right\|_{n,\alpha}~\text{and}~ \|DQ(\psi)\phi\|_{n,\alpha}
&\leq& C\left(\|\phi\|_{n,\alpha}+\|\psi\|_{n+1,\alpha}\|\phi\|_{1,\alpha}\right),\label{ineq:TEDQ(psi)phi}\\
\left\|J_\psi\frac\phi{|\nabla\psi|}\circ A_\psi^{-1}\right\|_{n,\alpha}~\text{and}~ \|DQ(\psi)\phi\|_{n,\alpha}
&\leq& C\left(\|\phi\|_{n,\alpha}+\|\psi-\overline\psi\|_{n+1,\alpha}\|\phi\|_{1,\alpha}\right).\label{ineq:ATEDQ(psi)phi}
\end{eqnarray}
In turn, for $n\geq 0$ and $\psi\in\mathcal F_+$ in a $\|\cdot\|_{2,\alpha}$-neighborhood
(see proof below),
\begin{equation}
\|A_\psi^{-1}-A_{\overline\psi}^{-1}\|_{n,\alpha}\leq C\|\psi-\overline\psi\|_{n+1,\alpha}.
\label{ineq:TEDifferenceA(psi)Inverse}
\end{equation}

For $n\geq 0$, $F$ in the $\|\cdot\|_{2,\alpha}$-neighborhood $\mathcal V_S(\overline F)$
(see proof below),
\begin{equation}\|\omega-\overline\omega\|_{n,\alpha}\leq C\|\psi-\overline\psi\|_{n+2,\alpha}\leq C\|F-\overline F\|_{n+1,\alpha}.
\label{ineq:TEDifferenceomegaAndpsi}
\end{equation}
Combining these two estimates we obtain
\begin{equation}
\|A_\psi^{-1}-A_{\overline\psi}^{-1}\|_{n,\alpha}\leq C\|F-\overline F\|_{n,\alpha}.
\end{equation}
One would expect from the above that $\|A_\omega^{-1}-A_{\overline\omega}^{-1}\|_{n,\alpha}\leq  C\|F-\overline F\|_{n+2,\alpha}$.
However, the regularizing effect of $\omega=\Delta\psi=F(\psi)$ gives the  better estimates
(see proof below):
\begin{equation}
\|A_\omega^{-1}-A_{\overline\omega}^{-1}\|_{n,\alpha}\leq C\|F-\overline F\|_{n,\alpha}.
\label{ineq:TEDifferenceA(omega)InverseForSteadyStates}
\end{equation}

For $n\geq 0$, $h,\overline h\in C^\infty_I$, and $\psi,\overline \psi\in C^\infty_{\overline\Omega}$ in a $\|\cdot\|_{1,\alpha}$-neighborhood
(see proof below),
\begin{equation}
\|h(\psi)-\overline h(\overline \psi)\|_{n,\alpha}
\leq C\bigg\{\|h\|_{n+1,\alpha}\|\psi-\overline\psi\|_{0,\alpha}+\|h\|_{2,\alpha}\|\psi-\overline\psi\|_{n,\alpha}+\|h-\overline h\|_{n,\alpha}\bigg\}.
\label{ineq:TEDifferenceh(psi)}
\end{equation}
\bigskip%

Observe then that for $F$ in the $\|\cdot\|_{2,\alpha}$-neighborhood $\mathcal V_B(\overline F)\subset \mathcal V_S(\overline F)$,
\begin{equation}
\|F\|_{2,\alpha},\quad \|\omega\|_{2,\alpha},\quad \|\psi\|_{4,\alpha},\quad \|A_\omega^{-1}\|_{2,\alpha},\quad \|A_\psi^{-1}\|_{4,\alpha}
\quad \leq \quad C.
\end{equation}
This allows to incorporate a number of terms of lower order into a constant.
\bigskip%

We also recall the standard interpolation inequalities
\begin{equation}
\begin{array}{ccl}
\|u\|_{n+r,\alpha}\|v\|_{s+m,\alpha}&\leq&C\left\{\|u\|_{n+r+m,\alpha}\|v\|_{s,\alpha}+\|u\|_{s,\alpha}\|v\|_{n+r+m,\alpha}\right\},\\
\|u\|_{n+r,\alpha}\|u\|_{s+m,\alpha}&\leq&C\|u\|_{n+r+m,\alpha}\|u\|_{s,\alpha}
\end{array}
\label{ineq:SummaryInterpolationInequalities}
\end{equation}
which are consequences of (\ref{ineq:IntroInterpolationInequalities}).
\bigskip%

\noindent\textbf{Proof of (\ref{ineq:TEDifferenceA(psi)Inverse})}\qquad
Set $\omega_t=\overline\omega+t(\omega-\overline\omega)$ and note that
$Q(\omega)-Q(\overline\omega)=\int_0^1DQ(\omega_t)(\omega-\overline\omega)dt$.
Then, $\omega_t$ remains in a $\|\cdot\|_{2,\alpha}$-neighborhood and by (\ref{ineq:ATEDQ(psi)phi}) we have
\begin{eqnarray}
\|A_\omega^{-1}-A_{\omega_0}^{-1}\|_{m,\alpha}
&\leq& \int_0^1\left\|DQ(\omega_t)(\omega-\overline\omega)\right\|_{m,\alpha}dt\\
&\leq& C\Big\{\|\omega-\overline\omega\|_{m,\alpha}+\|\omega-\overline\omega\|_{m+1,\alpha}\|\omega-\overline\omega\|_{1,\alpha}\Big\}\\
&\leq&C\|\omega-\overline\omega\|_{m+1,\alpha}
\end{eqnarray}
since $\|\omega_t-\overline\omega\|_{j,\alpha}\leq \|\omega-\overline\omega\|_{j,\alpha}$ for $0\leq t\leq 1$.
\stopproof%

\noindent\textbf{Proof of (\ref{ineq:TEDifferenceomegaAndpsi})}\qquad
Write $\omega-\overline\omega=\Delta\psi-\Delta\overline\psi$, and
set $F_t=\overline F+t(F-\overline F)$ and $\psi_t=S(F_t)$.
Then, $\psi-\overline\psi=\int_0^1\phi_tdt$ 
where $\Delta\phi_t=F_t'(\psi_t)\phi_t+(F-\overline F)(\psi_t)$.
Then by  (\ref{ineq:TEphi=phi(F,g)}),
\begin{eqnarray}
\|\phi_t\|_{n+2,\alpha}
&\leq& C\left(\|F-\overline F\|_{n,\alpha}+\|F_t\|_{n+1,\alpha}\|F-\overline F\|_{1,\alpha}\right)\\
&\leq&C\left(\|F-\overline F\|_{n,\alpha}+(\|F-\overline F\|_{n+1,\alpha}+1)\|F-\overline F\|_{1,\alpha}\right)\\
&\leq& C\|F-\overline F\|_{n+1,\alpha}.
\end{eqnarray}
Integrating, this gives the desired estimate.
\stopproof%
\noindent\textbf{Proof of (\ref{ineq:TEDifferenceA(omega)InverseForSteadyStates})}
Write
\begin{eqnarray}
A_\omega^{-1}-A_{\overline\omega}^{-1}
&=&(F-\overline F)(A_\psi^{-1})+\overline F(A_\psi^{-1})-\overline F(A_{\overline \psi}^{-1})\\
&=&(F-\overline F)(A_\psi^{-1})+\int_0^1\overline F'(A_{\overline\psi}^{-1}+t(A_\psi^{-1}-A_{\overline \psi}^{-1}))(A_\psi^{-1}-A_{\overline \psi}^{-1})dt.
\end{eqnarray}
The first term is estimated by
\begin{eqnarray}
\|(F-\overline F)(A_\psi^{-1})\|_{n,\alpha}
&\leq& \|F-\overline F\|_{n,\alpha}+\|A_\psi^{-1}\|_{n,\alpha}\|F-\overline F\|_{1,\alpha}\\
&\leq&\|F-\overline F\|_{n,\alpha}+\|A_\psi^{-1}-A_{\overline \psi}^{-1}\|_{n,\alpha}\|F-\overline F\|_{1,\alpha}\\
&\leq&C\|F-\overline F\|_{n,\alpha}
\end{eqnarray}
while we have
\begin{eqnarray}
\|\overline F'(A_{\overline\psi}^{-1}+t(A_\psi^{-1}-A_{\overline \psi}^{-1}))\|_{n,\alpha}
&\leq&C\left(1+\|A_{\overline\psi}^{-1}+t(A_\psi^{-1}-A_{\overline \psi}^{-1})\|_{n,\alpha}\right)\\
&\leq&C\left(1+\|A_\psi^{-1}-A_{\overline\psi}^{-1}\|_{n,\alpha}\right)\\
&\leq&C\left(1+\|F-\overline F\|_{n,\alpha}\right)
\end{eqnarray}
so that the second term is also bounded by $\|F-\overline F\|_{n,\alpha}$.
\stopproof%

\noindent\textbf{Proof of (\ref{ineq:TEDifferenceh(psi)})}\qquad
Write $\psi_t=\overline \psi+t(\psi-\overline\psi)$.
Then
\begin{eqnarray}
h(\psi)-\overline h(\overline\psi)=h(\psi)-h(\overline \psi)+(h-\overline h)(\overline \psi)
=\int_0^1h'(\psi_t)(\psi-\overline\psi)dt+(h-\overline h)(\overline\psi)
\end{eqnarray}
so that by (\ref{ineq:TameEstimatesForProductOfFunctions}) and (\ref{ineq:TameEstimatesForCompositionOperator}),
\begin{eqnarray}
&&\|h(\psi)-\overline h(\overline\psi)\|_{n,\alpha}\\
&\leq&C
\bigg\{\Big(\|h\|_{n+1,\alpha}+\|\psi-\overline\psi\|_{n,\alpha}\|h\|_{2,\alpha}\Big)\|\psi-\overline\psi\|_{0,\alpha}\\
&&+\|h\|_{2,\alpha}\|\psi-\overline\psi\|_{n,\alpha}+\|h-\overline h\|_{n,\alpha}\bigg\}\\
&\leq& C\bigg\{\|h\|_{n+1,\alpha}\|\psi-\overline\psi\|_{0,\alpha}+\|h\|_{2,\alpha}\|\psi-\overline\psi\|_{n,\alpha}+\|h-\overline h\|_{n,\alpha}\bigg\}.
\end{eqnarray}
\stopproof%

\subsection{$\text{Id}+K(F)g$ has a tame family of inverses}\label{section:DT(F)IsSurjectiveWithSmoothTameRightInverse}
Recall that $K(F)g$ is defined in (\ref{eq:K(F)=Ktilde(F)VB(F)}).
The procedure to construct a tame family of inverses to $\text{Id}+K(F)g$ completely
parallels the proof that $\Delta\phi+c\phi=k$ (augmented with suitable boundary conditions)
has a tame family of inverses, see Lemma~\ref{lemma:EstimatesForLinearEllipticEquations}
and Proposition~\ref{proposition:VE(c)kIsSmoothTame}.
In this case though, the estimates on $K(F)g-K(\overline F)g$ require significantly more work
and they are derived in the separate Lemma~\ref{lemma:TameEstimatesForK(F)g-K(Fbar)g}.

\begin{lemma}[Estimates for $h=g+K(F)g$]\label{lemma:EstimatesOnh=g+K(F)g}
Let $h=g+K(F)g$ where  $F\in\mathcal V_B(\overline F)\subset \mathcal V_S(\overline F)$ and $g\in C^\infty_{[0,|\Omega|]}$.
Then, for any $n\geq 0$,
\begin{equation}
\|g\|_{n,\alpha}\leq C\Big\{\|h\|_{n,\alpha}+\|F\|_{n+1,\alpha}\|g\|_{1,\alpha}\Big\}
\label{ineq:EstimatesForh=g+K(F)g}
\end{equation}
where $C$ may depend on $n$ and $\mathcal V_B(\overline F)$.
\end{lemma}
\startproof%
The estimates on $g$ are derived by writing $g=h-K(F)g$.
We first estimate $K(F)g$: by (\ref{ineq:TameEstimatesForProductOfFunctions}), (\ref{ineq:TEDQ(psi)phi}), (\ref{ineq:TEomegaAndpsi}), and (\ref{ineq:TEphi=phi(F,g)}),
\begin{eqnarray}
\|K(F)g\|_{n,\alpha}&\leq&C\left\|A_\omega^{-1}\right\|_{n+1,\alpha}\left\|J_\psi \frac\phi{|\nabla\psi|}\circ A_\psi^{-1}\right\|_{0,\alpha}
 +C\left\|A_\omega^{-1}\right\|_{1,\alpha}\left\|J_\psi\frac\phi{|\nabla\psi|}\circ A_\psi^{-1}\right\|_{n,\alpha}\\
&\leq&C(\|\omega\|_{n+1,\alpha}+1)\|\phi\|_{1,\alpha}
+C\big(\|\phi\|_{n,\alpha}+\|\psi\|_{n+1,\alpha}\left\|\phi\right\|_{1,\alpha}\big)\\
&\leq&C\Big\{\|\phi\|_{n,\alpha}+\|F\|_{n+1,\alpha}\|\phi\|_{1,\alpha}\Big\}\\
&\leq&C\Big\{\|g\|_{n-2,\alpha}+\|F\|_{n+1,\alpha}\|g\|_{1,\alpha}\Big\}.
\end{eqnarray}
for $F\in\mathcal V_B(\overline F)\subset \mathcal V_S(\overline F)$ and any $g$.
In turn,
\begin{eqnarray}
\|g\|_{n,\alpha}&\leq&\|h\|_{n,\alpha}+\|K(F)g\|_{n,\alpha}\\
&\leq&\|h\|_{n,\alpha}+C\left\{\|g\|_{n-2,\alpha}+\|F\|_{n+1,\alpha}\|g\|_{1,\alpha}\right\}
\end{eqnarray}
and the term $\|g\|_{n-2,\alpha}$ can be incorporated into the left-hand side as usual and we find
\begin{eqnarray}
\|g\|_{n,\alpha}&\leq&C\Big\{\|h\|_{n,\alpha}+\|F\|_{n+1,\alpha}\|g\|_{1,\alpha}\Big\}.
\end{eqnarray}
\stopproof%

We will need estimates on the difference $K(F)g-K(\overline F)g$ for $F\in\mathcal V_B(\overline F)$.
Recall that $K(F)g$ is defined in (\ref{eq:K(F)=Ktilde(F)VB(F)}).
\begin{lemma}[Estimates on $K(F)g-K(\overline F)g$]\label{lemma:TameEstimatesForK(F)g-K(Fbar)g}
For $n\geq 0$, $F\in\mathcal V_S(\overline F)$, and $g\in C^\infty_{[0,|\Omega|]}$, 
\begin{equation}
\|K(F)g-K(\overline F)g\|_{n,\alpha}
\leq C\Big\{\|F-\overline F\|_{n+1,\alpha}\|g\|_{2,\alpha}
+\|F-\overline F\|_{2,\alpha}\|g\|_{n-1,\alpha}\Big\}.
\label{ineq:TameEstimatesForK(F)g-K(Fbar)g}
\end{equation}
\end{lemma}
\startproof%
\setcounter{counterpar}{0}
$F$ and $g$ being given, let  $f=VB(F)g$, $\overline f=VB(\overline F)g$
and $\phi, \overline\phi$ such that
\begin{equation}\Delta\phi=F'(\psi)\phi+f(\psi),\qquad \Delta\overline\phi=\overline F'(\overline\psi)\overline\phi+\overline f(\overline\psi).
\end{equation}
Write then
\begin{eqnarray}
K(F)g-K(\overline F)g
&=&\left(A_\omega^{-1}-A_{\overline \omega}^{-1}\right)'\left(J_\psi\circ \frac\phi{|\nabla\psi|}\circ A_\psi^{-1}\right)\\
&&+\left(A_{\overline\omega}^{-1}\right)'\left(J_\psi\frac\phi{|\nabla\psi|}\circ A_\psi^{-1}-J_{\overline\psi}\frac\phi{|\nabla\psi|}\circ A_{\overline\psi}^{-1}\right)\\
&&+\left(A_{\overline\omega}^{-1}\right)'J_{\overline\psi}\left(\frac\phi{|\nabla\psi|}-\frac{\overline\phi}{|\nabla\overline\psi|}\right)\circ A_{\overline\psi}^{-1}\\
&=&I+II+III.
\end{eqnarray}

\nextpar{Estimates on I in terms of $\|\phi\|_{j,\alpha}$}
By (\ref{ineq:TameEstimatesForProductOfFunctions}), 
(\ref{ineq:TEDifferenceA(omega)InverseForSteadyStates}), (\ref{ineq:ATEDQ(psi)phi}), and (\ref{ineq:TEDifferenceomegaAndpsi}), 
\begin{eqnarray}
&&\|I\|_{n,\alpha}\\
&\leq& C\bigg\{\|A_\omega^{-1}-A_{\overline\omega}^{-1}\|_{n+1,\alpha}\left\|J_\psi\frac\phi{|\nabla\psi|}\circ A_\psi^{-1}\right\|_{0,\alpha}
+\|A_\omega^{-1}-A_{\overline\omega}^{-1}\|_{1,\alpha}\left\|J_\psi\frac\phi{|\nabla\psi|}\circ A_\psi^{-1}\right\|_{n,\alpha}\bigg\}\\
&\leq& C\bigg\{\|F-\overline F\|_{n+1,\alpha}\|\phi\|_{1,\alpha}
+\|F-\overline F\|_{2,\alpha}\left[\|\phi\|_{n,\alpha}+\|\psi-\overline\psi\|_{n+1,\alpha}\|\phi\|_{1,\alpha}\right]\bigg\}\\
&\leq&C\bigg\{\|F-\overline F\|_{n+1,\alpha}\|\phi\|_{1,\alpha}
+\|F-\overline F\|_{2,\alpha}\|\phi\|_{n,\alpha}+\|F-\overline F\|_{2,\alpha}\|F-\overline F\|_{n,\alpha}\|\phi\|_{1,\alpha}\bigg\}\\
&\leq&C\Big\{\|F-\overline F\|_{n+1,\alpha}\|\phi\|_{1,\alpha}
+\|F-\overline F\|_{2,\alpha}\|\phi\|_{n,\alpha}\Big\}.
\label{ineq:K(F)g-K(Fbar)gEstimatesOnI}
\end{eqnarray}
\stopstep%

\nextpar{Estimates on II in terms of $\|\phi\|_{j,\alpha}$}
As a preliminary we derive the following estimates for $\psi\in\mathcal F_+$ in a $\|\cdot\|_{1,\alpha}$-neighborhood 
and any $u\in C^\infty_{\overline\Omega}$.
Write $\psi_t=\overline \psi+t(\psi-\overline\psi)$, $N_t=\frac{\nabla\psi_t}{|\nabla\psi_t|}$  and
\begin{eqnarray}
\left\|J_\psi u\circ A_\psi^{-1}-J_{\overline\psi} u \circ A_{\overline\psi}^{-1}\right\|_{n,\alpha}
\leq\int_0^1\left\|D_\psi\left(J_{\psi_t}u\circ A_{\psi_t}^{-1}\right)(\psi-\overline\psi)\right\|_{n,\alpha}dt.
\end{eqnarray}
Using tame estimates (\ref{ineq:TameEstimatesForDerivativeOfJ(omega)uOfA(omega)InverseInomega}) of the
derivative of $J_\psi u\circ A_\psi^{-1}$ in $\psi$, and (\ref{ineq:TEDifferenceomegaAndpsi}),
\begin{eqnarray}
&&\left\|D_\psi\left(J_{\psi_t}u\circ A_{\psi_t}^{-1}\right)(\psi-\overline\psi)\right\|_{n,\alpha}\\
&\leq&C\bigg\{\|u\|_{n+1,\alpha}\|\psi-\overline \psi\|_{1,\alpha}+\|u\|_{2,\alpha}\|\psi-\overline\psi\|_{n,\alpha}
+\|\psi_t\|_{n+2,\alpha}\|u\|_{2,\alpha}\|\psi-\overline\psi\|_{1,\alpha}\bigg\}\\
&\leq&C\bigg\{\|u\|_{n+1,\alpha}\|\psi-\overline \psi\|_{1,\alpha}+\|u\|_{2,\alpha}\|\psi-\overline\psi\|_{n,\alpha}
+\|\psi-\overline \psi\|_{n+2,\alpha}\|u\|_{2,\alpha}\bigg\}\\
&\leq&C\bigg\{\|u\|_{n+1,\alpha}\|\psi-\overline \psi\|_{1,\alpha}+\|\psi-\overline\psi\|_{n+2,\alpha}\|u\|_{2,\alpha}\bigg\}\\
&\leq&C\bigg\{\|u\|_{n+1,\alpha}\|F-\overline F\|_{1,\alpha}+\|F-\overline F\|_{n+1,\alpha}\|u\|_{2,\alpha}\bigg\}.
\end{eqnarray}
using $\|\psi_t\|_{n+2,\alpha}\leq \|\overline\psi\|_{n+2,\alpha}+\|\psi-\overline\psi\|_{n+2,\alpha}\leq C(1+\|\psi-\overline\psi\|_{n+2,\alpha})$.
Now with $u=\frac\phi{|\nabla\psi|}$, we have by (\ref{ineq:TEDifferenceomegaAndpsi}),
\begin{eqnarray}
\left\|u\right\|_{m,\alpha}
&\leq& C\Big\{\|\phi\|_{m,\alpha}+\|\psi-\overline\psi\|_{m+1,\alpha}\|\phi\|_{0,\alpha}\Big\}\\
&\leq& C\Big\{\|\phi\|_{m,\alpha}+\|F-\overline F\|_{m,\alpha}\|\phi\|_{0,\alpha}\Big\}
\end{eqnarray}
so that 
\begin{eqnarray}
&&\|II\|_{n,\alpha}\\
&\leq&C\Big\{\|\phi\|_{n+1,\alpha}\|F-\overline F\|_{1,\alpha}+\|F-\overline F\|_{n+1,\alpha}\|\phi\|_{0,\alpha}
+\|F-\overline F\|_{n+1,\alpha}\|\phi\|_{2,\alpha}\Big\}\\
&\leq&C\Big\{\|\phi\|_{n+1,\alpha}\|F-\overline F\|_{1,\alpha}+\|F-\overline F\|_{n+1,\alpha}\|\phi\|_{2,\alpha}\Big\}.
\label{ineq:K(F)g-K(Fbar)gEstimatesOnII}
\end{eqnarray}
\stopstep%

\nextpar{Estimates on I+II in terms of $\|g\|_{m,\alpha}$}
We may now estimate $\|I\|_{n,\alpha}+\|II\|_{n,\alpha}$ in terms of $g$ instead of $\phi$ 
(the third term $III$ must be dealt with differently).
Putting (\ref{ineq:K(F)g-K(Fbar)gEstimatesOnI}) and (\ref{ineq:K(F)g-K(Fbar)gEstimatesOnII}) together,
and using interpolation inequalities (\ref{ineq:SummaryInterpolationInequalities}),
and (\ref{ineq:ATEphi=phi(F,g)}),
\begin{eqnarray}
&&\|I\|_{n,\alpha}+\|II\|_{n,\alpha}\\
&\leq& C\Big\{\|F-\overline F\|_{n+1,\alpha}\|\phi\|_{2,\alpha}+\|F-\overline F\|_{2,\alpha}\|\phi\|_{n+1,\alpha}\Big\}\\
&\leq& C\Big\{\|F-\overline F\|_{n+1,\alpha}\|g\|_{1,\alpha}
+\|F-\overline F\|_{2,\alpha}\left(\|g\|_{n-1,\alpha}+\|F-\overline F\|_{n,\alpha}\|g\|_{1,\alpha}\right)\Big\}\\
&\leq& C\Big\{\|F-\overline F\|_{n+1,\alpha}\|g\|_{1,\alpha}
+\|F-\overline F\|_{2,\alpha}\|g\|_{n-1,\alpha}\Big\}.
\label{ineq:K(F)g-K(Fbar)gEstimatesOnIAndII}
\end{eqnarray}

\nextpar{Estimates on III}
Since $\overline\omega$ is fixed, we have by (\ref{ineq:TameEstimatesForProductOfFunctions})
\begin{eqnarray}
\|III\|_{n,\alpha}
&\leq&C
\left\|J_{\overline\psi}\left\{\left(\frac\phi{|\nabla\psi|}-\frac{\overline\phi}{|\nabla\overline\psi|}\right)\right\}\circ A_{\overline\psi}^{-1}\right\|_{n,\alpha}\\
&\leq&C\left\|\left(\frac1{|\nabla\psi|}-\frac1{|\nabla\overline\psi|}\right)\phi\right\|_{n,\alpha}
+C\left\|\frac1{|\nabla\overline\psi|}(\phi-\overline\phi)\right\|_{n,\alpha}\\
&=&III_1+III_2.
\end{eqnarray}
Using (\ref{ineq:TameEstimatesForProductOfFunctions}), (\ref{ineq:TEDifferenceomegaAndpsi}) and (\ref{ineq:ATEphi=phi(F,g)}), 
the first term is estimated by
\begin{eqnarray}
III_1&\leq&C\left\{\left\|\frac1{|\nabla\psi|}-\frac1{|\nabla\overline\psi|}\right\|_{0,\alpha}\|\phi\|_{n,\alpha}
 +\left\|\frac1{|\nabla\psi|}-\frac1{|\nabla\overline\psi|}\right\|_{n,\alpha}\|\phi\|_{0,\alpha}\right\}\\
&\leq&C\Big\{\|\psi-\overline \psi\|_{2,\alpha}\|\phi\|_{n,\alpha}+\|\psi-\overline \psi\|_{n+2,\alpha}\|\phi\|_{0,\alpha}\Big\}\\
&\leq&C\Big\{\|F-\overline F\|_{1,\alpha}\|\phi\|_{n,\alpha}+\|F-\overline F\|_{n+1,\alpha}\|\phi\|_{0,\alpha}\Big\}\\
&\leq&C\Big\{\|F-\overline F\|_{1,\alpha}\left[\|g\|_{n-2,\alpha}+\|F-\overline F\|_{n-1,\alpha}\|g\|_{1,\alpha}\right]
+\|F-\overline F\|_{n+1,\alpha}\|g\|_{1,\alpha}\Big\}\\
&\leq&C\Big\{\|F-\overline F\|_{1,\alpha}\|g\|_{n-2,\alpha}+\|F-\overline F\|_{n+1,\alpha}\|g\|_{1,\alpha}\Big\}.
\label{ineq:K(F)g-K(Fbar)gEstimatesOnIII-1}
\end{eqnarray}
For $III_2$ we ``just'' need to estimate $\|\phi-\overline\phi\|_{n,\alpha}$.
Note that
\begin{equation}
\Delta(\phi-\overline\phi)-\overline F'(\overline\psi)(\phi-\overline\phi)
=(F'(\psi)-\overline F'(\overline\psi))\phi+f(\psi)-\overline f(\overline\psi).
\end{equation}
From tame estimates (\ref{ineq:CnalphaTameEstimatesOnVE(c)k}), 
(observe as well  that $\overline F$ and $\overline \psi$ are fixed), we have for $j\geq 0$
\begin{eqnarray}
\|\phi-\overline\phi\|_{j+2,\alpha}
&\leq&C\bigg\{\|(F'(\psi)-\overline F'(\overline \psi))\phi\|_{j,\alpha}+\|f(\psi)-\overline f(\overline\psi)\|_{j,\alpha}\bigg\}.
\end{eqnarray}

\nextpar{Estimates on $\|\phi-\overline \phi\|_{j+2,\alpha}$: $\|(F'(\psi)-\overline F'(\overline \psi))\phi\|_{j,\alpha}$}
With the above, using (\ref{ineq:TEDifferenceh(psi)}),
(\ref{ineq:TameEstimatesForProductOfFunctions}) and (\ref{ineq:TameEstimatesForCompositionOperator}), 
the first term is estimated as follows:
\begin{eqnarray}
&&\|(F'(\psi)-\overline F'(\overline \psi))\phi\|_{j,\alpha}\\
&\leq&C\|F'(\psi)-\overline F'(\overline \psi)\|_{j,\alpha}\|\phi\|_{0,\alpha}
+C\|F'(\psi)-\overline F'(\overline \psi)\|_{0,\alpha}\|\phi\|_{j,\alpha}\\
&\leq&
C\Big(\|F\|_{j+2,\alpha}\|\psi-\overline\psi\|_{0,\alpha}+\|F\|_{3,\alpha}\|\psi-\overline\psi\|_{j,\alpha}+\|F-\overline F\|_{j+1,\alpha}\Big)
\|\phi\|_{0,\alpha}\\
&&+C\Big(\|F\|_{3,\alpha}\|\psi-\overline\psi\|_{0,\alpha}+\|F-\overline F\|_{1,\alpha}\Big)\|\phi\|_{j,\alpha}\\
&\leq&
C\bigg((\|F-\overline F\|_{j+2,\alpha}+1)\|\psi-\overline\psi\|_{0,\alpha}+(\|F-\overline F\|_{3,\alpha}+1)\|\psi-\overline\psi\|_{j,\alpha}\\
&&+\|F-\overline F\|_{j+1,\alpha}\bigg) \|\phi\|_{0,\alpha}\\
&&+\Big((\|F-\overline F\|_{3,\alpha}+1)\|\psi-\overline\psi\|_{0,\alpha}+\|F-\overline F\|_{1,\alpha}\Big)\|\phi\|_{j,\alpha}\\
&\leq&C\Big\{\|F-\overline F\|_{j+2,\alpha}\|\phi\|_{0,\alpha}+\|F-\overline F\|_{3,\alpha}\|\phi\|_{j,\alpha}\Big\}
\end{eqnarray}
where we have used that 
\begin{equation}\|F-\overline F\|_{3,\alpha}\|\psi-\overline\psi\|_{j,\alpha}\leq C\|F-\overline F\|_{3,\alpha}\|F-\overline F\|_{j-1,\alpha}
\leq C\|F-\overline F\|_{0,\alpha}\|F-\overline F\|_{j+2,\alpha}.
\end{equation}
In turn, using (\ref{ineq:ATEphi=phi(F,g)}),
\begin{eqnarray}
&&\|(F'(\psi)-\overline F'(\overline \psi))\phi\|_{j,\alpha}\\
&\leq&C\Big\{\|F-\overline F\|_{j+2,\alpha}\|g\|_{1,\alpha}+\|F-\overline F\|_{3,\alpha}\bigg(\|g\|_{j-2,\alpha}+\|F-\overline F\|_{j-1,\alpha}\|g\|_{1,\alpha}\bigg)\Big\}\\
&\leq&C\Big\{\|F-\overline F\|_{j+2,\alpha}\|g\|_{1,\alpha}+\|F-\overline F\|_{3,\alpha}\|g\|_{j-2,\alpha}\Big\}\\
&\leq&C\Big\{\|F-\overline F\|_{j+2,\alpha}\|g\|_{1,\alpha}+\|F-\overline F\|_{1,\alpha}\|g\|_{j,\alpha}\Big\}
\label{ineq:FirstTermIn||phi-phibar||(j+2,alpha)}
\end{eqnarray}
using again interpolation inequalities (\ref{ineq:SummaryInterpolationInequalities}) to get
$\|F-\overline F\|_{3,\alpha}\|g\|_{j-2,\alpha}\leq C\|F-\overline F\|_{1,\alpha}\|g\|_{j,\alpha}+C\|F-\overline F\|_{j,\alpha}\|g\|_{1,\alpha}$,
$\|F-\overline F\|_{3,\alpha}\|F-\overline F\|_{j-2,\alpha}\leq C\|F-\overline F\|_{2,\alpha}\|F-\overline F\|_{j,\alpha}\leq C\|F-\overline F\|_{j,\alpha}$
and to incorporate these into other terms.

\nextpar{Estimates on $\|\phi-\overline \phi\|_{j+2,\alpha}$: $\|f(\psi)-\overline f(\overline \psi)\|_{j,\alpha}$}
Using again (\ref{ineq:TEDifferenceh(psi)}),
\begin{eqnarray}
\|f(\psi)-\overline f(\overline\psi)\|_{j,\alpha}
&\leq& C\bigg\{\|f\|_{j+1,\alpha}\|\psi-\overline\psi\|_{0,\alpha}+\|f\|_{2,\alpha}\|\psi-\overline\psi\|_{j,\alpha}+\|f-\overline f\|_{j,\alpha}\bigg\}.
\end{eqnarray}
Using (\ref{ineq:ATEf=VB(F)g}), 
the first two terms in the right-hand side are estimated by
\begin{eqnarray}
&&\|f\|_{j+1,\alpha}\|\psi-\overline\psi\|_{0,\alpha}+\|f\|_{2,\alpha}\|\psi-\overline\psi\|_{j,\alpha}\\
&\leq&C\Big((\|g\|_{j+1,\alpha}+\|F-\overline F\|_{j-1,\alpha}\|g\|_{1,\alpha})\|F-\overline F\|_{1,\alpha}
+\|g\|_{2,\alpha}\|F-\overline F\|_{j-1,\alpha}\Big)\\
&\leq&C\Big(\|g\|_{j+1,\alpha}\|F-\overline F\|_{1,\alpha}+\|F-\overline F\|_{j-1,\alpha}\|g\|_{2,\alpha}\Big).
\label{ineq:FirstTwoTermsInf(psi)-fbar(psibar)}
\end{eqnarray}
As for the last term in $\|f(\psi)-\overline f(\overline \psi)\|_{j,\alpha}$,
write it as
\begin{equation}
f-\overline f=\mathcal Eg\circ A-\mathcal Eg\circ \overline A
\qquad\text{with}\qquad A=V(\mathcal E(A_\psi^{-1})),\qquad \overline A=V(\mathcal E(A_{\overline\psi}^{-1})),
\end{equation}
to find again from (\ref{ineq:TEDifferenceh(psi)})
\begin{equation}
\|f-\overline f\|_{j,\alpha}
\leq C\bigg\{\|\mathcal Eg\|_{j+1,\alpha}\|A-\overline A\|_{0,\alpha}+\|\mathcal Eg\|_{2,\alpha}\|A-\overline A\|_{j,\alpha}\bigg\}.
\end{equation}
In order to estimates $\|A-\overline A\|_{j,\alpha}$,
pose $\beta=\mathcal E(A_\psi^{-1})$, $\overline\beta=\mathcal E(A_{\overline \psi}^{-1})$,
and $\beta_t=\overline\beta+t(\beta-\overline\beta)$, 
and use the inversion operator $V$, see Lemma~\ref{lemma:InversionIsSmoothTame} in the Appendix:
\begin{equation}
A-\overline A=V(\beta)-V(\overline\beta)=\int_0^1\dot a_tdt,
\qquad \dot a_t=
DV(\beta_t)\cdot(\beta-\overline \beta)
=-\frac{(\beta-\overline\beta)\circ \beta_t^{-1}}{\beta_t'\circ \beta_t^{-1}}.
\end{equation}
We have by (\ref{ineq:TEDifferenceA(psi)Inverse}) and (\ref{ineq:TEDifferenceomegaAndpsi})
\begin{equation}
\|\beta-\overline\beta\|_{m,\alpha}
\leq C\|A_\psi^{-1}-A_{\overline\psi}^{-1}\|_{m,\alpha}
\leq C\|\psi-\overline\psi\|_{m+1,\alpha}
\leq C\|F-\overline F\|_{m,\alpha}.
\label{ineq:||beta-betabar||(m,alpha)}
\end{equation}
Using the triangle inequality on $\beta_t=\overline \beta+t(\beta-\overline\beta)$
and (\ref{ineq:TEDifferenceomegaAndpsi}), 
\begin{equation}\|\beta_t\|_{m,\alpha}\leq C(1+\|\beta-\overline\beta\|_{m,\alpha})\leq C\left\{\|F-\overline F\|_{m,\alpha}+1\right\}.
\label{ineq:||betat||(m,alpha)}
\end{equation}
This along with Lemma~\ref{lemma:InversionIsSmoothTame} implies
\begin{equation}
\|\beta_t^{-1}\|_{m,\alpha}\leq C\left\{\|F-\overline F\|_{m,\alpha}+1\right\}.
\label{ineq:||betatInverse||(m,alpha)}
\end{equation}
In particular, $\beta_t$ and $\beta_t^{-1}$ remain in $\|\cdot\|_{2,\alpha}$-neighborhoods for $F,\overline F\in\mathcal V_B(\overline F)$.
Now the fraction in $\dot a_t$ is linear in the numerator, so that
by (\ref{ineq:TameEstimatesForProductOfFunctions}) and (\ref{ineq:TameEstimatesForCompositionOperator}),
\begin{eqnarray}
\|\dot a_t\|_{j,\alpha}
&\leq&C\|(\beta-\overline\beta)\circ \beta_t^{-1}\|_{j,\alpha}
+C\|\beta_t'\circ \beta_t^{-1}\|_{j,\alpha}\|(\beta-\overline\beta)\circ \beta_t^{-1}\|_{0,\alpha}\\
&\leq&C\Big\{\|\beta-\overline\beta\|_{j,\alpha}+\|\beta_t^{-1}\|_{j,\alpha}\|\beta-\overline\beta\|_{1,\alpha}\bigg\}\\
&&+C\Big(\|\beta_t\|_{j+1,\alpha}+ \|\beta_t^{-1}\|_{j,\alpha}\|\beta_t\|_{2,\alpha}\Big)
\|\beta-\overline\beta\|_{1,\alpha}\\
&\leq&C\|F-\overline F\|_{j,\alpha}
+C\Big(\|F-\overline F\|_{j+1,\alpha}+1\Big)\|F-\overline F\|_{1,\alpha}\\
&\leq&C\|F-\overline F\|_{j+1,\alpha}
\end{eqnarray}
from (\ref{ineq:||beta-betabar||(m,alpha)}), (\ref{ineq:||betat||(m,alpha)}), and 
(\ref{ineq:||betatInverse||(m,alpha)}).
With the above, we have
\begin{equation}
\|A-\overline A\|_{j,\alpha}\leq \int_0^1\|\dot a_t\|_{j,\alpha}dt
\leq C\|F-\overline F\|_{j+1,\alpha}
\end{equation}
so that by $\|\mathcal Eg\|_{m,\alpha}\leq C\|g\|_{m,\alpha}$ (see (\ref{ineq:TameEstimatesForExtensionOperator}))
\begin{equation}
\|f-\overline f\|_{j,\alpha}\leq C\Big\{\|g\|_{j+1,\alpha}\|F-\overline F\|_{1,\alpha}+\|g\|_{2,\alpha}\|F-\overline F\|_{j+1,\alpha}\Big\}.
\end{equation}
Adding this to (\ref{ineq:FirstTwoTermsInf(psi)-fbar(psibar)}), we find
\begin{equation}
\|f(\psi)-\overline f(\overline \psi)\|_{j,\alpha}
\leq C\Big\{\|F-\overline F\|_{1,\alpha}\|g\|_{j+1,\alpha}+\|F-\overline F\|_{j+1,\alpha}\|g\|_{2,\alpha}\Big\}.
\label{ineq:SecondTermIn||phi-phibar||(j+2,alpha)}
\end{equation}

\nextpar{Conclusion on III}
Putting (\ref{ineq:FirstTermIn||phi-phibar||(j+2,alpha)}) and (\ref{ineq:SecondTermIn||phi-phibar||(j+2,alpha)}) together
(with $n=j+2$)
and using interpolation inequalities (\ref{ineq:SummaryInterpolationInequalities}), 
\begin{equation}
\|III_2\|_{n,\alpha}\leq C\|\phi-\overline\phi\|_{n,\alpha}
\leq C\bigg\{\|F-\overline F\|_{n,\alpha}\|g\|_{2,\alpha}+\|F-\overline F\|_{1,\alpha}\|g\|_{n-1,\alpha}\bigg\}.
\end{equation}

Adding this to the estimate (\ref{ineq:K(F)g-K(Fbar)gEstimatesOnIII-1}) on $\|III_1\|_{n,\alpha}$ gives
\begin{equation}
\|III\|_{n,\alpha}\leq C\bigg\{\|F-\overline F\|_{n+1,\alpha}\|g\|_{2,\alpha}+\|F-\overline F\|_{1,\alpha}\|g\|_{n-1,\alpha}\bigg\}.
\label{ineq:K(F)g-K(Fbar)gEstimatesOnIII}
\end{equation}
\stopstep%

\nextpar{Conclusion on Lemma~\ref{lemma:TameEstimatesForK(F)g-K(Fbar)g}}
Putting (\ref{ineq:K(F)g-K(Fbar)gEstimatesOnIAndII})
and (\ref{ineq:K(F)g-K(Fbar)gEstimatesOnIII}) together, we have
\begin{equation}
\|K(F)g-K(\overline F)g\|_{n,\alpha}
\leq C\bigg\{\|F-\overline F\|_{n+1,\alpha}\|g\|_{2,\alpha}+\|F-\overline F\|_{2,\alpha}\|g\|_{n-1,\alpha}\bigg\}.
\end{equation}
\stopproof%

\begin{proposition}[$\text{Id}+K(F)$ is invertible for $F$ near $\overline F$]\label{proposition:TameInverseToh=g+K(F)g}
Suppose $\overline \omega=\overline F(\overline \psi)$ satisfies (ND2).
Then, there exists a $\|\cdot\|_{3,\alpha}$-neighborhood 
\begin{equation}\mathcal V_I(\overline F)\subset \mathcal V_B(\overline F)\subset \mathcal V_S(\overline F)\end{equation}
of $\overline F$ such that $h=g+K(F)g$ has a tame family of inverses $g=VM(F)h$ satisfying
\begin{equation}
\|g\|_{n,\alpha}\leq C\bigg\{\|h\|_{n,\alpha}+\|F\|_{n+1,\alpha}\|h\|_{2,\alpha}\bigg\},\qquad n\geq 2.
\label{ineq:TameEstimatesForg=VM(F)h}
\end{equation}
\end{proposition}
\noindent\textbf{Remark}\qquad
The Nash-Moser Inverse Function Theorem only requires a continuous and tame inverse. 
However, with a little more work, one can show that $g=VM(F)h$ is continuously differentiable,
hence smooth tame by Theorem~5.3.1, p.~102, Part~I, and Theorem~3.1.1, p.~150, Part~II of \cite{HamiltonIFTNashMoser}.
\medskip%

\startproof
\setcounter{counterpar}{0}

\nextpar{The estimate $\|g\|_{2,\alpha}\leq C\|g+K(\overline F)g\|_{2,\alpha}$ holds}
We first show that (ND2) implies that $\text{Id}_{C^\infty_{[0,|\Omega|]}}+K(\overline F)$ has trivial kernel
as a map $C^\infty_{[0,|\Omega|]}\rightarrow C^\infty_{[0,|\Omega|]}$.
Suppose $g\in C^\infty_{[0,|\Omega|]}$ is in the kernel 
of $\text{Id}_{C^\infty_{[0,|\Omega|]}}+K(\overline F)\colon C^\infty_{[0,|\Omega|]}\rightarrow C^\infty_{[0,|\Omega|]}$.
Set $f=VB(\overline F)\cdot g$ so that $DT(\overline F)\cdot f=0$ and thus $\nu=\{\overline\omega,\alpha\}$
for some $\alpha\in \mathcal U$, see Proposition~\ref{proposition:AomegaCharacterizesO*(omega)}.
Then $f=0$ on $\text{range}(\overline\psi)$, precisely by the non-degeneracy condition (ND2)
and in turn $g=B(\overline F)\cdot VB(\overline F)\cdot g=B(\overline F)\cdot f=0$, 
i.e. $\text{Id}_{C^\infty_{[0,|\Omega|]}}+K(\overline F)\colon C^\infty_{[0,|\Omega|]}\rightarrow C^\infty_{[0,|\Omega|]}$
has trivial kernel.

Next, we show that $\text{Id}_{C^{2,\alpha}_{[0,|\Omega|]}}+K(\overline F)\colon C^{2,\alpha}_{[0,|\Omega|]}\rightarrow C^{2,\alpha}_{[0,|\Omega|]}$,
satisfies
\begin{equation}\|g\|_{2,\alpha}\leq C\|g+K(\overline F)g\|_{2,\alpha}.\label{ineq:||g||2,alpha<C||g+K(Fbar)g||2,alpha}
\end{equation}
First observe from the tame estimates on $K(F)g$ derived in the proof of Lemma~\ref{lemma:EstimatesOnh=g+K(F)g}
that $K(\overline F)$ maps $C^{n,\alpha}_{[0,|\Omega|]}$ into $C^{n+1,\alpha}_{[0,|\Omega|]}$ for each $n\geq 2$.
Let $g\in C^{2,\alpha}_{[0,|\Omega|]}$ be in the kernel of 
$\text{Id}_{C^{2,\alpha}_{[0,|\Omega|]}}+K(\overline F)\colon C^{2,\alpha}_{[0,|\Omega|]}\rightarrow C^{2,\alpha}_{[0,|\Omega|]}$.
Then $g=-K(\overline F)g\in C^{3,\alpha}_{[0,|\Omega|]}$ and, repeating, $g\in C^\infty_{[0,|\Omega|]}$.
Hence $g=0$ which shows that 
$\text{Id}_{C^{2,\alpha}_{[0,|\Omega|]}}+K(\overline F)\colon C^{2,\alpha}_{[0,|\Omega|]}\rightarrow C^{2,\alpha}_{[0,|\Omega|]}$
has trivial kernel.
By the Fredholm alternative (for Banach spaces), 
$\text{Id}_{C^{2,\alpha}_{[0,|\Omega|]}}+K(\overline F)\colon C^{2,\alpha}_{[0,|\Omega|]}\rightarrow C^{2,\alpha}_{[0,|\Omega|]}$
is an isomorphism satisfying (\ref{ineq:||g||2,alpha<C||g+K(Fbar)g||2,alpha}).
\stopstep%

\nextpar{Invertibility for $F$ near $\overline F$}
From Lemma~\ref{lemma:TameEstimatesForK(F)g-K(Fbar)g}
\begin{eqnarray}
\|g\|_{2,\alpha}&\leq&C\|g+K(\overline F)g\|_{2,\alpha}\\
&\leq&C\left(\|g+K(F)g\|_{2,\alpha}+\|K(F)g-K(\overline F)g\|_{2,\alpha}\right)\\
&\leq&C\left(\|g+K(F)g\|_{2,\alpha}+\|F-\overline F\|_{3,\alpha}\|g\|_{2,\alpha}\right).
\end{eqnarray}
Define now $\mathcal V_I(\overline F)\subset \mathcal V_B(\overline F)\subset \mathcal V_S(\overline F)$ 
to be a sufficiently small $\|\cdot\|_{3,\alpha}$-neighborhood of $\overline F$ so that the last term can be incorporated into the left-hand side:
\begin{equation}\|g\|_{2,\alpha}\leq C\|g+K(F)g\|_{2,\alpha}.\end{equation}

Now from (\ref{ineq:EstimatesForh=g+K(F)g}) we have for $n\geq 2$
\begin{eqnarray}
\|g\|_{n,\alpha}&\leq& C\left(\|g+K(F)g\|_{n,\alpha}+\|F-\overline F\|_{n+1,\alpha}\|g\|_{1,\alpha}\right)\\
&\leq&C\left(\|g+K(F)g\|_{n,\alpha}+\|F-\overline F\|_{n+1,\alpha}\|g+K(\overline F)g\|_{2,\alpha}\right).
\label{ineq:AlternativeTameEstimatesForh=g+K(F)g}
\end{eqnarray}
which implies the estimates (\ref{ineq:TameEstimatesForg=VM(F)h}).
That is, $h=M(F)g=g+K(F)g$ is a family of invertible linear maps
of Fr\'echet spaces for $F\in\mathcal V_I(\overline F)$ with tame inverse denoted $g=VM(F)h$.
\stopstep%

\nextpar{Continuity}
Let $F,\tilde F\in\mathcal V_I(\overline F)$, $g,\tilde g\in C^\infty_{[0,|\Omega|]}$, and set
$h=g+K(F)g$, $\tilde h=\tilde g+K(\tilde F)\tilde g$.
Then,
\begin{equation}
(g-\tilde g)+K(\tilde F)(g-\tilde g)=h-\tilde h-(K(F)g-K(\tilde F)g)
\end{equation}
and from (\ref{ineq:AlternativeTameEstimatesForh=g+K(F)g}) 
we deduce that
\begin{eqnarray}
\|g-\tilde g\|_{n,\alpha}
&\leq &C\bigg\{\|h-\tilde h\|_{n,\alpha}+\|K(F)g-K(\tilde F)g\|_{n,\alpha}\\
&&+\|F-\tilde F\|_{n+1,\alpha}\Big(\|h-\tilde h\|_{2,\alpha}+\|K(F)g-K(\overline F)g\|_{2,\alpha}\Big)\bigg\}\\
\end{eqnarray}
Now from the tame estimates (\ref{ineq:TameEstimatesForg=VM(F)h}) it is clear that $\|g\|_{n,\alpha}$
remains bounded as $\|h-\tilde h\|_{n,\alpha}$ and $\|F-\tilde F\|_{n+1,\alpha}$ tend to zero.
In turn the estimates (\ref{ineq:TameEstimatesForK(F)g-K(Fbar)g}) (also valid for $K(F)g-K(\tilde F)g$)
show that $\|K(F)g-K(\tilde F)g\|_{n,\alpha}$ tends to zero as well,
and thus clearly $\|g-\tilde g\|_{n,\alpha}$ tends to zero.
\stopproof%

\subsection{Injective part of Theorem~ \ref{theorem:MainTheorem}}\label{section:InjectivePart}
The injective part of Theorem~\ref{theorem:MainTheorem} 
requires a modification of the injective part of the Nash-Moser theorem
since the map $T$ cannot be injective:
defining $T$ on $C^\infty_I$ where $I\supset\text{range}(\overline\psi)$ (see (\ref{eq:IntervalI})),
changing $F$ outside the range of the corresponding solution $\Delta\psi=F(\psi)$ clearly does not 
affect this solution.
\begin{proposition}\label{proposition:InjectivePartOfMainTheorem}
There exists a $C^2$-neighborhood of $\overline F$ such that, if $F,\tilde F$ are such that $T(F)=T(\tilde F)$,
then the corresponding solutions $\psi=S(F)$ and $\tilde \psi=S(\tilde F)$ are the same.
\end{proposition}
\startproof%
Let $F_n\neq\tilde F_n\rightarrow_n \overline F$ in $C^2$ such that $T(F_n)=T(\tilde F_n)$, and let
\begin{equation}
\tilde F_n-F_n=\epsilon_nG_n,\quad \psi_n=S(F_n),\quad \tilde \psi_n=S(\tilde F_n),\quad \tilde \psi_n-\psi_n=\epsilon_nv_n
\end{equation}
where $\epsilon_n$ is to be chosen appropriately.
Assume without loss of generality that $\text{range}(\psi_n)\subset \text{range}(\tilde\psi_n)$.
Supposing that we can normalize according to
\begin{equation}
\|G_n\|_{C^0_{\text{range}(\tilde\psi_n)}}=1,\label{eq:AdjustedNormalization}
\end{equation}
we will arrive at a contradiction, thus proving our claim.

\noindent\textbf{Remark}\quad One might expect to normalize in the $C^2$-norm since the $F$'s converge in that norm.
However, (strong) compactness fails in infinite dimensions.
On the other hand, one can use the gain of regularity 
provided by $\Delta\psi=F(\psi)$.

Changing $F_n$ and $\tilde F_n$ outside of $\text{range}(\tilde\psi_n)$ does not affect $\psi_n$ nor $\tilde \psi_n$,
and in turn the assumption $T(\tilde F_n)=T(F_n)$ is preserved.
Therefore, we may adjust $F_n$ and $\tilde F_n$ in such a way that, without loss of generality,
for the \foremph{first} derivatives we have the bounds
\begin{equation}\|G_n\|_{C^1_{[\overline c,0]}}\leq 2\|G_n\|_{C^1_{\text{range}(\tilde\psi_n)}}.\label{ineq:BoundedExtensionInC1}
\end{equation}

Since $T(\tilde F_n)=T(F_n)$, we have
$\epsilon_nG_n=F_n\circ A_{\psi_n}^{-1}\circ A_{\tilde\psi_n}-F_n$
hence 
\begin{equation}
\epsilon_nG_n(\lambda)
=\left(\int_0^1F_n'(\lambda+t(A_{\psi_n}^{-1}(A_{\tilde\psi_n}(\lambda))-\lambda))dt\right)(A_{\psi_n}^{-1}(A_{\tilde\psi_n}(\lambda))-\lambda),
\qquad \lambda\in\text{range}(\tilde\psi_n).
\end{equation}
Letting $d_n$ denote the integral factor, we write this as
\begin{equation}
\epsilon_nG_n(\lambda)=d_n(\lambda)(A_{\psi_n}^{-1}-A_{\tilde\psi_n}^{-1})(A_{\tilde\psi_n}(\lambda)),
\qquad \lambda\in\text{range}(\tilde\psi_n).
\end{equation}
Now $F_n$ and $\tilde F_n$ are bounded in $C^2$,
$F\in C^2\mapsto \psi\in C^2$ is continuous by (\ref{eq:ContinuityOfSolutionOperatorInCn-Grading}),
and $\psi\in C^2\mapsto A_\psi^{-1}\in C^1$ is continuous by Lemma~\ref{lemma:QIsContinuous}.
Thus, $d_n$ and $A_{\tilde \psi_n}$ are bounded in $C^1$ and in turn
\begin{eqnarray}
\epsilon_n\|G_n\|_{C^1_{\text{range}(\tilde\psi_n)}}
&\leq& C\|A_{\psi_n}^{-1}-A_{\tilde\psi_n}^{-1}\|_1.\label{ineq:LowerIneqForCompactnessOfGn}
\end{eqnarray}

Write
\begin{equation}
A_{\tilde\psi_n}^{-1}-A_{\psi_n}^{-1}
=\epsilon_n\left(\int_0^1DQ(\psi_n+t(\tilde\psi_n-\psi_n))dt\right)v_n
\end{equation}
where $\psi_n$ and $\tilde\psi_n$ are bounded in $C^2$ since $F_n$ and $\tilde F_n$ converge in $C^2$.
The proof of Lemma~\ref{lemma:QIsSmoothTame} easily shows that (\ref{ineq:TEDQ(psi)phi}) holds in the $\|\cdot\|_n$-grading,
so that
\begin{equation}
\epsilon_n\|G_n\|_{C^1_{\text{range}(\tilde\psi_n)}}\leq C\|A_{\tilde\psi_n}^{-1}-A_{\psi_n}^{-1}\|_1\leq \epsilon_nC\|v_n\|_1.\label{ineq:MiddleIneqForCompactnessOfGn}
\end{equation}

Finally, rewrite $\Delta \tilde \psi_n-\Delta \psi_n=\tilde F_n(\tilde \psi_n)-F_n(\psi_n)$ as
\begin{equation}
\epsilon_n\Delta v_n=\epsilon_n G_n(\tilde \psi_n)+\left(F_n(\tilde\psi_n)-F_n(\psi_n)\right).
\end{equation}
This can be written as $\Delta v_n+c_nv_n=G_n(\tilde\psi_n)$
\begin{equation}-c_n=\left(\int_0^1F_n'(\psi_n+t(\tilde \psi_n-\psi_n))dt\right)v_n.\end{equation}
$\overline F$ satisfies (ND1) hence for large $n$, $c_n$ is sufficiently close to $\overline F'(\overline\psi)$
and in turn $\Delta+c_n$ is invertible with
\begin{eqnarray}
\|v_n\|_{1}&\leq&C\|G_n\|_0.\label{ineq:UpperIneqForCompactnessOfGn}
\end{eqnarray}

Putting (\ref{ineq:BoundedExtensionInC1}), (\ref{ineq:MiddleIneqForCompactnessOfGn}),
and (\ref{ineq:UpperIneqForCompactnessOfGn}) together we obtain
\begin{equation}\|G_n\|_1\leq 2\|G_n\|_{C^1_{\text{range}(\tilde \psi_n)}}\leq C\|v_n\|_1\leq C\|G_n\|_0\leq C.
\end{equation}
Passing to a subsequence, we may assume that 
\begin{equation}G_n\rightarrow_n G\qquad \text{in}\qquad C^0.\end{equation}

Lemma~\ref{lemma:QIsContinuous} says that 
$(\omega,\nu)\mapsto DQ(\omega)\nu$
is continuous as an operator $C^2\times C^0\rightarrow C^0$,
(\ref{eq:ContinuityOfF->omegaInC^n-Grading}) shows that $F\mapsto \omega=F(\psi)$ is continuous $C^2\rightarrow C^2$,
while $(F,f)\mapsto \nu=\Delta\phi=F'(\psi)\phi+f(\psi)$ is continuous $C^2\times C^0\rightarrow C^0$
by (\ref{eq:ContinuityOf(F,f)->nuInC^n-Grading}).
In conclusion, we have the continuous operator
\begin{equation}
\left\{\begin{array}{ccc}
C^2&\times&C^0\\F&&f\end{array}
\longrightarrow\begin{array}{c}C^0\\DT(F)f
\end{array}\right\}.
\label{eq:ContinuityOfDT(F)f}
\end{equation}
Taking limits in 
\begin{equation}
0=\frac1{\epsilon_n}\left(T(F_n+\epsilon_n G_n)-T(F_n)\right)=\int_0^1DT(F_n+\epsilon_n(\tilde F_n-F_n))G_ndt
\label{eq:(T(Ftilden)-T(Fn))/epsilonn}
\end{equation}
one finds $DT(\overline F)G=0$.
This means that $G$ vanishes on $\text{range}(\overline\psi)$,
contradicting the normalization $\|G\|_{C^0_{\text{range}(\overline\psi)}}=1$ guaranteed by
(\ref{eq:AdjustedNormalization}).
\stopproof%

\section{Appendix: the Nash-Moser Inverse Function Theorem}\label{section:NashMoserIFT}
Inverse function theorems express the fact that nonlinear problems are as solvable as their linearizations:
a nonlinear map $T$ is (locally) surjective where its first derivative $DT$ is surjective,
and $T$ is (locally) injective where $DT$ is injective.
In case $T\colon (B\subset X)\rightarrow Y$ is a sufficiently smooth 
(e.g. twice continuously differentiable) map of Banach spcaes,
Newton's scheme constructs successive approximations which converge very rapidly.
This ``accelerated convergence'' is visible through an estimate of the form $x_{n+1}\leq x_n^2$.
For comparison, a proof by the Picard approximation method would involve an iteration of the form $x_{n+1}\leq \lambda x_n$
with some fixed $0<\lambda<1$.

Loss of derivatives occurs when, for instance, a surjective first derivative $DT(F)$
uses a number of derivatives, which are not recovered by its right-inverse $L(F)$.
(This is the case of our map $T(F)=A_\omega^{-1}$.)
The Newton algorithm can no longer be implemented as such.
On the other hand, its accelerated convergence indicates that it should tolerate some adjustments 
made in order to overcome loss of derivatives.

In \cite{NashImbeddingProblem} Nash introduced smoothing operators  
for his solution to the isometric imbedding problem of Riemannian manifolds into Euclidean spaces.
From \cite{NashImbeddingProblem} Moser extracted a simple algorithm
solving an inverse function problem even when loss of derivatives occurs \cite{Moser}.
We will refer to this algorithm as the Moser scheme.
It is a modified Newton scheme
where the smoothing operators of Nash introduce an error term having no effect on the convergence
of the algorithm
provided the maps satisfy certain (``tame'') estimates.
The solution so obtained is ``rough'' \textit{a priori}
(e.g. if one works with spaces of functions,
the solution may have fewer derivatives than the formulation of the problem actually allows).
In a second step, one verifies that this rough solution is in fact smooth.
This again uses the ``tameness'' of certain operators, as well as interpolation 
inequalities available in ``tame Fr\'echet-spaces'' in a crucial way.

Various extensions and improvements have been developed subsequently.
In particular, Hamilton introduced in \cite{HamiltonIFTNashMoser} the \defn{tame Fr\'echet category}
(see Sections~II.1 and II.2),
essentially that introduced by Sergeraert in \cite{Sergeraert},
in which the modified Newton algorithm is applicable and therefore an Inverse Function Theorem holds.
That is, an inverse function exists and lives in the tame Fr\'echet category.

We emphasize that the Moser scheme is used to construct a rough solution to $T(x)=y$
when the map $T$ has surjective first derivative ($DT$ is not required to be injective.)
That this solution is smooth is a consequence of the interpolation inequalities
available on ``tame Fr\'echet-spaces'',
and surely the estimates on the successive approximations play a part in the proof.
This is the surjective part of the Inverse Function Theorem.
In case $T$ has injective first derivative $DT$, then the interpolation inequalities (\ref{ineq:IntroInterpolationInequalities})
again show that $T$ is injective as well.
The Moser scheme plays no r\^ole in this injective part of the Inverse Function Theorem.
Our map $T(F)=A_\omega^{-1}$ cannot be injective (see the discussion in the Introduction).
Nevertheless, the injective part of the Inverse Function Theorem for tame Fr\'echet spaces
gives the idea for the proof of the injective part of Theorem~\ref{theorem:MainTheorem}.
\bigskip%

\subsection{Tame estimates on $T(F)=A_\omega^{-1}$}
In this Section we derive the precise tame estimates on $T(F)$ so as to set 
the parameters for the proof of the (existence part of the) Nash-Moser Inverse Function Theorem,
see Theorem~\ref{theorem:SurjectivePartOfNashMoserIFT}.
\bigskip%

Recall from Proposition~\ref{proposition:TameInverseToh=g+K(F)g} that 
$\mathcal V_I(\overline F)$ is a $\|\cdot\|_{3,\alpha}$-neighborhood of $\overline F$.
\begin{proposition}
For $F\in V_I(\overline F)$, any $f_1,f_2\in C^\infty_I$, and $n\geq 0$,
\begin{eqnarray}
\|T(F)\|_{n,\alpha}&\leq&C(\|F\|_{n,\alpha}+1),\\
\|DT(F)f\|_{n,\alpha}&\leq& C(\|f\|_{n,\alpha}+\|F\|_{n+1,\alpha}\|f\|_{1,\alpha}),\\
\|D^2T(F)(f_1,f_2)\|_{n,\alpha}
&\leq&C(\|f_1\|_{n+1,\alpha}\|f_2\|_{2,\alpha}+\|f_1\|_{2,\alpha}\|f_2\|_{n+1,\alpha}\\
&&+\|F\|_{n+2,\alpha}\|f_1\|_{2,\alpha}\|f_2\|_{2,\alpha}),
\end{eqnarray}
and for $n\geq 2$ and $h\in C^\infty_{[0,|\Omega|]}$,
\begin{equation}
\|L(F)h\|_{n,\alpha}\leq C(\|h\|_{n,\alpha}+\|F\|_{n+1,\alpha}\|h\|_{2,\alpha}).
\end{equation}
\end{proposition}
In the proof of the Nash-Moser theorem below, we will
shift the indices in the norms and use the notation
\begin{equation}|\cdot|_n:=\|\cdot\|_{n+2,\alpha},\qquad n\geq 0.\end{equation}
Note that $\mathcal V_I(\overline F)$ is then a $|\cdot|_1$-neighborhood.

\startproof%
With $T(F)=Q(\omega)$ where $\omega=\Delta\psi$ and $\psi=S(F)$,
from Propositions~\ref{proposition:AomegaInverseIsSmoothTame}
and \ref{proposition:SolutionMappsi=S(F)}, we have for $n\geq 0$
and $F$ in the $\|\cdot\|_{2,\alpha}$-neighborhood $\mathcal V_S(\overline F)$, 
see Proposition~\ref{proposition:SolutionMappsi=S(F)},
\begin{equation}
\|T(F)\|_{n,\alpha}\leq C(\|\omega\|_{n,\alpha}+1)
\leq C(\|\psi\|_{n+2,\alpha}+1)
\leq C(\|F\|_{n,\alpha}+1).
\end{equation}

Write the first derivative
as  $DT(F)f=DQ(\omega)\nu$ where $\nu=\Delta\phi=F'(\psi)\phi+f(\psi)$.
Again from Propositions~\ref{proposition:AomegaInverseIsSmoothTame}
and \ref{proposition:SolutionMappsi=S(F)}, we have for $n\geq 0$,
$F\in \mathcal V_S(\overline F)$, and any $f\in C^\infty_I$,
\begin{eqnarray}
&&\|DT(F)f\|_{n,\alpha}\\
&\leq&C\bigg(\|\nu\|_{n,\alpha}+\|\omega\|_{n+1,\alpha}\|\nu\|_{1,\alpha}\bigg)\\
&\leq&C\bigg(\|\phi\|_{n+2,\alpha}+(\|F\|_{n+1,\alpha}+1)\|\phi\|_{3,\alpha}\bigg)\\
&\leq& C\bigg(\|f\|_{n,\alpha}+\|F\|_{n+1,\alpha}\|f\|_{1,\alpha}+(\|F\|_{n+1,\alpha}+1)(\|f\|_{1,\alpha}+\|F\|_{2,\alpha}\|f\|_{1,\alpha})\\
&\leq& C\bigg(\|f\|_{n,\alpha}+\|F\|_{n+1,\alpha}\|f\|_{1,\alpha}\bigg).
\end{eqnarray}

Write the second derivative as
\begin{equation}
D^2T(F)(f_1,f_2)
=D^2Q(\omega)(\nu_1,\nu_2)+DQ(\omega)\nu_{12}\end{equation}
where
\begin{eqnarray}
&\nu_1=\Delta\phi_1=F'(\psi)\phi_1+f_1(\psi),&\\
&\nu_2=\Delta\phi_2=F'(\psi)\phi_2+f_2(\psi),&\\
&\nu_{12}=\Delta\phi_{12}=F'(\psi)\phi_{12}+F'(\psi)\phi_1\phi_2+f_1'(\psi)\phi_2+f_2'(\psi)\phi_1,&
\end{eqnarray}
see Proposition~\ref{proposition:SolutionMappsi=S(F)}.
Thus, for $F$ in the $\|\cdot\|_{3,\alpha}$-neighborhood $\mathcal V_I(\overline F)$,
see Proposition~\ref{proposition:TameInverseToh=g+K(F)g},
and any $f_1,f_2\in C^\infty_I$, we have for $n\geq 0$
\begin{eqnarray}
&&\|D^2T(F)(f_1,f_2)\|_{n,\alpha}\\
&\leq&C\Bigg\{
\|\nu_1\|_{n+1,\alpha}\|\nu_2\|_{1,\alpha}+\|\nu_1\|_{1,\alpha}\|\nu_2\|_{n+1,\alpha}+\|\omega\|_{n+2,\alpha}\|\nu_1\|_{2,\alpha}\|\nu_2\|_{2,\alpha}\\
&&+\|\nu_{12}\|_{n,\alpha}+\|\omega\|_{n+1,\alpha}\|\nu_{12}\|_{1,\alpha}\Bigg\}\\
&\leq&C\Bigg\{
\|\phi_1\|_{n+3,\alpha}\|\phi_2\|_{3,\alpha}+\|\phi_1\|_{3,\alpha}\|\phi_2\|_{n+3,\alpha}+\|\omega\|_{n+2,\alpha}\|\phi_1\|_{4,\alpha}\|\phi_2\|_{4,\alpha}\\
&&+\|\phi_{12}\|_{n+2,\alpha}+\|\omega\|_{n+1,\alpha}\|\phi_{12}\|_{3,\alpha}\Bigg\}\\
&\leq&C\Bigg\{
\bigg(\|f_1\|_{n+1,\alpha}+\|F\|_{n+2,\alpha}\|f_1\|_{1,\alpha}\bigg)\bigg(\|f_2\|_{1,\alpha}+\|F\|_{2,\alpha}\|f_2\|_{1,\alpha}\bigg)\\
&&+\bigg(\|f_1\|_{1,\alpha}+\|F\|_{2,\alpha}\|f_1\|_{1,\alpha}\bigg)\bigg(\|f_2\|_{n+1,\alpha}+\|F\|_{n+2,\alpha}\|f_2\|_{1,\alpha}\bigg)\\
&&+(\|F\|_{n+2,\alpha}+1)\bigg(\|f_1\|_{2,\alpha}+\|F\|_{3,\alpha}\|f_1\|_{1,\alpha}\bigg)\bigg(\|f_2\|_{2,\alpha}+\|F\|_{3,\alpha}\|f_2\|_{1,\alpha}\bigg)\\
&&+\|f_1\|_{n+1,\alpha}\|f_2\|_{1,\alpha}+\|f_1\|_{1,\alpha}\|f_2\|_{n+1,\alpha}+
\|F\|_{n+2,\alpha}\|f_1\|_{2,\alpha}\|f_2\|_{2,\alpha}\\
&&+(\|F\|_{n+1,\alpha}+1)(\|f_1\|_{2,\alpha}\|f_2\|_{2,\alpha}+\|F\|_{3,\alpha}\|f_1\|_{2,\alpha}\|f_2\|_{2,\alpha})
\Bigg\}\\
&\leq&C\Bigg\{
\|f_1\|_{n+1,\alpha}\|f_2\|_{2,\alpha}+\|f_1\|_{2,\alpha}\|f_2\|_{n+1,\alpha}
+\|F\|_{n+2,\alpha}\|f_1\|_{2,\alpha}\|f_2\|_{2,\alpha}\Bigg\}\\
&\leq&C\Bigg\{
\|f_1\|_{n+1,\alpha}\|f_2\|_{2,\alpha}+\|f_1\|_{2,\alpha}\|f_2\|_{n+1,\alpha}
+\|F\|_{n+2,\alpha}\|f_1\|_{2,\alpha}\|f_2\|_{2,\alpha}\Bigg\}.
\end{eqnarray}

Finally we compute the tame estimates for the right-inverse $L(F)h$ to $DT(F)f$.
Recall that it is given by $f=L(F)h=VB(F)\cdot VM(F)h$.
With the tame estimates on $f=VB(F)g$ from Lemma~\ref{lemma:B(F)fHasASmoothTameFamilyOfRightInverses}
and those on $g=VM(F)h$ from Proposition~\ref{proposition:TameInverseToh=g+K(F)g},
we deduce that, for $n\geq 2$,
\begin{eqnarray}
\|f\|_{n,\alpha}
&\leq&C\Bigg\{\|g\|_{n,\alpha}+\|F\|_{n-2,\alpha}\|g\|_{1,\alpha}\Bigg\}\\
&\leq&C\Bigg\{\|h\|_{n,\alpha}+\|F\|_{n+1,\alpha}\|h\|_{2,\alpha}
+\|F\|_{n-2,\alpha}(\|h\|_{1,\alpha}+\|F\|_{2,\alpha}\|h\|_{2,\alpha})\Bigg\}\\
&\leq&C\Bigg\{\|h\|_{n,\alpha}+\|F\|_{n+1,\alpha}\|h\|_{2,\alpha}\Bigg\}.
\end{eqnarray}
\stopproof%

\subsection{Surjective part of the Nash-Moser Inverse Function Theorem}
Our presentation of the Nash-Moser Inverse Function Theorem is a blend of \cite{HamiltonIFTNashMoser}, \cite{Moser}, and \cite{Sergeraert}.
We will limit ourselves to constructing a right-inverse, since in our application to Theorem~\ref{theorem:MainTheorem}
we do not need further properties of this right-inverse (smoothness and tameness).
We refer to Section~III.1 of \cite{HamiltonIFTNashMoser} for further details.
\bigskip%

Consider $\mathcal X$, $\mathcal Y$ tame Fr\'echet spaces, with smoothing operators
$S(t)$, $t>0$, satisfying the estimates (\ref{ineq:IntroEstimatesSmoothingOperators}) described in the Introduction,
and set 
\begin{equation}
\mathcal B:=\{u\in\mathcal X~|~|u-\overline u|_1<\eta\}.\label{eq:NbhdB}
\end{equation}
Let $T\colon (\mathcal B\subset \mathcal X)\rightarrow \mathcal Y$ such that
for any $u\in \mathcal B$, $v_1,v_2\in \mathcal X$,  $h\in\mathcal Y$, and $n\geq 0$,
\begin{eqnarray}
|T(u)|_n&\leq &C(|u-\overline u|_n+1),\\
|DT(u)v|_n&\leq&(|v|_n+|u-\overline u|_{n+1}|v|_0),\\
|D^2T(u)(v_1,v_2)|_n&\leq& C(|v_1|_{n+1}|v_2|_0+|v_1|_0|v_2|_{n+1}+|u-\overline u|_{n+2}|v_1|_0|v_2|_0),\\
|L(u)h|_n&\leq& C(|h|_n+|u-\overline u|_{n+1}|h|_0).
\end{eqnarray}
(It is clear that our map $T(F)=A_\omega^{-1}$ satisfies these conditions where $|\cdot|_n=\|\cdot\|_{n+2,\alpha}$.)
Suppose given a solution $T(\overline u)=\overline g$.

\begin{theorem}[Existence part of the Nash-Moser Inverse Function Theorem]\label{theorem:SurjectivePartOfNashMoserIFT}
\hfill \\%
There is a neighborhood $\mathcal G\subset \mathcal Y$ of $\overline g$ in which
$T(u)=g$ has a solution $u\in\mathcal B$ for any $g\in\mathcal G$.
\end{theorem}
\noindent\textbf{Remark}\qquad
The neighborhood $\mathcal G$ is defined in (\ref{eq:NbhdG}) in terms of a parameter $j$ 
given in (\ref{eq:OptimalValueForj}).
In particular, $\mathcal G$ is a $|\cdot|_8=\|\cdot\|_{10,\alpha}$-neighborhood.
Further, since $\omega\mapsto A_\omega^{-1}$ is continuous $C^{11}\rightarrow C^{11}\hookrightarrow C^{10,\alpha}$,
see Proposition~\ref{proposition:AomegaInverseIsSmoothTame},
the $C^\infty$-neighborhood of $\overline \omega$ in which each co-adjoint orbit
has a steady-state can then be taken as a $\|\cdot\|_{11}$-neighborhood.

\startproof%

\setcounter{counterpar}{0}
\nextpar{The modified Newton scheme}
For $g\in\mathcal Y$, we pose
\begin{equation}P(g,u)=T(u)-g\end{equation}
and the problem is to solve $P(g,u)=0$.
We will think of $g$ as a parameter. 
Then $P$ satisfies the following estimates
for any $u\in \mathcal B$, $v_1,v_2\in \mathcal X$, and $n\geq 0$,
\begin{eqnarray}
|P(g,u)|_n&\leq &C(|u-\overline u|_n+|g-\overline g|_n+1),\\
|D_uP(g,u)v|_n&\leq&(|v|_n+|u-\overline u|_{n+1}|v|_0),\\
|D^2_{uu}P(g,u)(v_1,v_2)|_n&\leq& C(|v_1|_{n+1}|v_2|_0+|v_1|_0|v_2|_{n+1}+|u-\overline u|_{n+2}|v_1|_0|v_2|_0)
\end{eqnarray}
and $v=L(u)h$ is again a right-inverse to $h=D_uP(g,u)v$:
\begin{equation}D_uP(g,u)L(u)h=h.\label{eq:DP(g,u)L(u)h=h}
\end{equation}

The solution is constructed by the Moser scheme,
which is a modified Newton scheme:
\begin{equation}
u_{n+1}-u_n:=-S(t_n)L(u_n)P(u_n),\qquad n\geq 0, \qquad u_0:=\overline u, \quad t_n:=A^{\kappa^n}
\label{eq:u(n+1)-u(n)=-S(t(n))L(u(n))P(u(n))}
\end{equation}
for some $A>1$ and $0<\kappa<2$ to be determined.
Fix $j\geq 1$ which will be specified later.
Let $M, M_j>1$ be constants such that
for all $u, w$ such that $u, u+w\in {\mathcal B}$, 
any $v\in \mathcal X$, any $h\in\mathcal Y$, and $t>0$, we have:
\begin{eqnarray}
|S(t)v|_{1}&\leq&Mt|v|_0,\\
|S(t)v|_j&\leq&M_jt|v|_{j-1}\\
|v-S(t)v|_0&\leq&M_jt^{1-j}|v|_{j-1},\\
|D_uP(u)v|_0&\leq&M|v|_0,\\
|P(u+w)-P(u)-D_uP(u)w|_0&\leq&M|w|_{1}^2,\\
|L(u)h|_0&\leq& M|h|_0.
\end{eqnarray}
(That these hold is immediate from the tame estimates on $T$, the estimates (\ref{ineq:IntroEstimatesSmoothingOperators})
on the smoothing operators $S(t)$, and Taylor's expansion with remainder (\ref{eq:TaylorSecondOrder})
given in the Introduction).
Next, with the first requirement that $\mathcal G$ be contained in a neighborhood of the form
\begin{equation}
\mathcal G\subset \{|g-\overline g|_{j-1}<C\}\subset \mathcal Y
\end{equation}
and increasing $M_j$ if necessary, we have for any $u\in\mathcal B$ and $g\in \mathcal G$
\begin{eqnarray}
|L(u)P(g,u)|_{j-1}&\leq&M_j(1+|u-\overline u|_j)
\end{eqnarray}
which holds since $|u|_0$ and $|g|_0$ remain bounded and
\begin{eqnarray}
|L(u)P(g,u)|_{j-1}&\leq& C(|P(g,u)|_{j-1}+|u-\overline u|_j|P(g,u)|_0)\\
&\leq& C(|u-\overline u|_{j-1}+|g-\overline g|_{j-1}+1+|u-\overline u|_j). 
\end{eqnarray}

Write now
\begin{eqnarray}
|u_{n+1}-u_n|_1&=&|S(t_n)L(u_n)P(g,u_n)|_1\\
&\leq&Mt_n|L(u_n)P(g,u_n)|_0\\
&\leq&M^2t_n|P(g,u_n)|_0
\end{eqnarray}
and by (\ref{eq:u(n+1)-u(n)=-S(t(n))L(u(n))P(u(n))}), (\ref{eq:DP(g,u)L(u)h=h})
\begin{eqnarray}
&&|P(g,u_n)|_0\\
&\leq&|P(g,u_n)-P(g,u_{n-1})-D_uP(g,u_{n-1})(u_n-u_{n-1})|_0\\
&& + |P(g,u_{n-1})+D_uP(g,u_{n-1})(u_n-u_{n-1})|_0\\
&\leq& M|u_n-u_{n-1}|_1^2+|D_uP(g,u_{n-1})(1-S(t_{n-1}))L(u_{n-1})P(g,u_{n-1})|_0\\
&\leq&M|u_n-u_{n-1}|_1^2+M|(1-S(t_{n-1}))L(u_{n-1})P(g,u_{n-1})|_0\\
&\leq&M|u_n-u_{n-1}|_1^2+MM_jt_{n-1}^{1-j}|L(u_{n-1})P(g,u_{n-1})|_{j-1}\\
&\leq&M|u_n-u_{n-1}|_1^2+MM_j^2t_{n-1}^{1-j}(1+|u_{n-1}-\overline u|_j)\label{ineq:EstimateFor|P(u(n))|0}
\end{eqnarray}
so that
\begin{equation}
|u_{n+1}-u_n|_1\leq t_nM^3|u_n-u_{n-1}|_1^2
 +M^3M_j^2t_nt_{n-1}^{1-j}(1+|u_{n-1}-\overline u|_j).
\end{equation}
For some $\mu>0$ to be determined later, let
\begin{equation}\delta_n:=t_n^\mu M^3|u_n-u_{n-1}|_1.\end{equation}
Then, 
\begin{equation}
\delta_{n+1}\leq A^{\kappa^n(1+\mu(\kappa-2))}\delta_n^2+e_n,\qquad
e_n:=M^6M_j^2A^{\mu\kappa^{n+1}+\kappa^n+(1-j)\kappa^{n-1}}(1+|u_{n-1}-\overline u|_j).
\label{ineq:FirstInequalityFordelta(n)}
\end{equation}
The parameters $A, \kappa,$ etc. will be determined in order to view (\ref{ineq:FirstInequalityFordelta(n)})
as a perturbation of $x_{n+1}\leq x_n^2$.
With 
\begin{equation}1+\mu(\kappa-2)\leq 0\label{ineq:1+mu(kappa-2)<0}\end{equation}
we have $\delta_{n+1}\leq \delta_n^2+e_n$.
By inspection, the graphs of $y=x$ and $y=x^2+\frac18$ intersect at some $x\in[\frac23,1]$.
Thus, if one can impose
\begin{equation}\delta_1<\frac23,\qquad e_n<\frac18,\label{ineq:ConditionsSoThatdelta(n)RemainsBounded}
\end{equation}
the inequality $\delta_{n+1}\leq \delta_n^2+e_n$ guarantees that $\delta_n$ is bounded for all $n$ (say by $1$)
and therefore
\begin{equation}|u_n-u_{n-1}|_1\leq \frac{A^{-\mu\kappa^n}}{M^3}.\label{ineq:EstimateOn|u(n)-u(n-1)|1}\end{equation}
In order to control $e_n$ we need to estimate the growth of $1+|u_n-\overline u|_j$:
\begin{eqnarray}
1+|u_{n+1}-\overline u|_j&\leq& 1+|u_n-\overline u|_j+|u_{n+1}-u_n|_j\\
&=& 1+|u_n-\overline u|_j+|S(t_n)L(u_n)P(g,u_n)|_j\\
&\leq&1+|u_n-\overline u|_j+M_jt_n|L(u_n)P(g,u_n)|_{j-1}\\
&\leq&1+|u_n-\overline u|_j+M_j^2t_n(1+|u_n-\overline u|_j)\\
&\leq&2M_j^2A^{\kappa^n}(1+|u_n-\overline u|_j).\label{ineq:InductionInequalityOn(1+|u(n)-ubar|j)}
\end{eqnarray}
Let $\beta\geq 0$ to be determined later and write
\begin{eqnarray}
A^{-\beta\kappa^{n+1}}(1+|u_{n+1}-\overline u|_j)
&\leq& 2M_j^2 A^{(-\beta(\kappa-1)+1)\kappa^n}\left(A^{-\beta\kappa^n} (1+|u_n-\overline u|_j)\right)\\
&\leq& 2M_j^2 A^{-\beta(\kappa-1)+1}\left(A^{-\beta\kappa^n}(1+|u_n-\overline u|_j)\right)
\end{eqnarray}
provided $-\beta(\kappa-1)+1<0$. Since we also want the multiplicative factor to be $\leq 1$,
we impose the more stringent condition that
\begin{equation}(-\beta(\kappa-1)+1)\ln A+\ln(2M_j^2)\leq 0.\label{ineq:ConditionOnbeta}
\end{equation}
In turn, as long as the terms $u_n\in {\mathcal B}$ exist, we have (recall $u_0=\overline u$, $A>1$, and $\beta\geq 0$)
\begin{equation}1+|u_n-\overline u|_j\leq A^{\beta\kappa^n}
\end{equation}
hence
$e_n \leq M^6M_j^2A^{(\mu\kappa^2+\kappa+1-j+\beta)\kappa^{n-1}}$.
In order to satisfy (\ref{ineq:ConditionsSoThatdelta(n)RemainsBounded}),
we will therefore impose that
\begin{equation}(\mu\kappa^2+\kappa+1-j+\beta)\ln A+\ln(M^6M_j^2)<\ln\frac18.
\label{ineq:ConditionForControlOn(1+|u(n)-ubar|j)}
\end{equation}

Finally, we show that the $u_n\in {\mathcal B}$ exists for all $n\geq 0$ provided $|P(g,\overline u)|_0<\epsilon$ is sufficiently small,
and that the sequence is Cauchy in the $|\cdot|_1$-norm.
Estimate (\ref{ineq:EstimateOn|u(n)-u(n-1)|1}) holds 
provided $\delta_1=A^{\mu\kappa}M^3|u_1-u_0|_1<2/3$.
But 
$|u_1-u_0|_{1}\leq Mt_0^2|L(\overline u)P(\overline u)|_{0}\leq M^2A^2|P(\overline u)|_{0}$
so we take
\begin{equation}\epsilon<\frac23\frac1{M^5A^{\mu\kappa+2}}.\label{ineq:BoundOnepsilon}\end{equation}

This in turn determines the neighborhood $\mathcal G$:
since $P(g,\overline u)=T(\overline u)-g=\overline g-g$, 
\begin{equation}
\mathcal G=\{g\in\mathcal Y~|~|g-\overline g|_0<\epsilon, |g|_{j-1}<C\}.
\label{eq:NbhdG}
\end{equation}

Now we verify that $u_n\in {\mathcal B}$ is defined for all $n$:
\begin{equation}
|u_{n+1}|_{1}
\leq \sum_{m=0}^n |u_{m+1}-u_m|_{1}
\leq \frac1{M^3}\sum_{m=0}^\infty A^{-\mu\kappa^{m+1}}
\leq \frac1{M^3}\sum_{m=0}^\infty A^{-\mu\kappa(1+m\ln\kappa)} 
=\frac1{M^3}\frac{A^{-\mu\kappa}}{1-A^{-\mu\kappa\ln\kappa}}
\end{equation}
(we have used 
$\kappa^i\geq \kappa^l+(i-l)\ln \kappa$).
We thus impose
\begin{equation}
\frac1{M^3}\frac{A^{-\mu\kappa}}{1-A^{-\mu\kappa\ln\kappa}}<\eta.
\label{ineq:ConditionForu(n)InB}
\end{equation}

The sequence is Cauchy since
\begin{eqnarray}
|u_m-u_n|_1
&\leq& \sum_{l=n}^{m-1}|u_{l+1}-u_l|_1\\
&\leq& \frac1{M^3}\sum_{l=n}^\infty A^{-\mu\kappa^{l+1}}\\
&\leq& \frac1{M^3}\sum_{l=n}^\infty A^{-\mu\kappa(\kappa^n+\ln\kappa (l-n))} \\
&=&\frac1{M^3}\frac{A^{-\mu\kappa^{n+1}}}{1-A^{-\mu\kappa\ln\kappa}}
\qquad\underset{n,m\rightarrow \infty}\rightarrow\qquad 0.
\end{eqnarray}
We denote $u_\infty$ its limit in the $|\cdot|_{1}$-norm.
Observe from (\ref{ineq:EstimateFor|P(u(n))|0}) that
\begin{eqnarray}
|P(g,u_n)|_0
\leq M|u_n-u_{n-1}|_1^2+MM_j^2A^{(1-j+\beta)\kappa^{n-1}}\rightarrow 0
\end{eqnarray}
since $1-j+\beta<0$ by (\ref{ineq:ConditionForControlOn(1+|u(n)-ubar|j)}).
Thus, once $u_\infty$ is proven to be the limit in each $|\cdot|_k$-norm, the above shows that it is in fact
a solution to $P(g,u)=0$ in $\mathcal B$.
\\

\nextpar{Setting the parameters}
The constants $M$, $M_j$ are imposed by the problem. 
The conditions on the parameters $\kappa$, $\mu$, $\beta$, $j$, and $A$ are
(\ref{ineq:1+mu(kappa-2)<0}), (\ref{ineq:ConditionOnbeta}), (\ref{ineq:ConditionForControlOn(1+|u(n)-ubar|j)}),
(\ref{ineq:ConditionForu(n)InB}):
\begin{eqnarray}
&1+\mu(\kappa-2)\leq 0,\qquad 0<\kappa<2&\\
&(-\beta(\kappa-1)+1)\ln A+\ln(2M_j^2)\leq 0&\\
&(\mu\kappa^2+\kappa+1-j+\beta)\ln A+\ln(M^6M_j^2)<\ln\frac18&\\
&\frac1{M^3}\frac{A^{-\mu\kappa}}{1-A^{-\mu\kappa\ln\kappa}}<\eta.&
\end{eqnarray}
These conditions are satisfied if the parameters $\kappa$, $\mu$, $\beta$, $j$, and $A$ 
are chosen \textit{in this order} so as to satisfy the following:
\begin{eqnarray}
&1<\kappa<2&\\
&\mu\geq \frac1{2-\kappa}&\\
&-\beta(\kappa-1)+1<0&\\
&\mu\kappa^2+\kappa+1-j+\beta<0&
\end{eqnarray}
and $A$ sufficiently large so that the three inequalities where it is involved are satisfied.
It is not difficult to see that, in order to minimize $j$, it should be chosen the smallest integer 
strictly larger than 
$\frac{\kappa^2}{2-\kappa}+\kappa+1+\frac1{\kappa-1}$ over $\kappa\in (1,2)$.
Using a computer one finds that 
\begin{equation}j=9\qquad \text{is attained with any}\qquad 1.2<\kappa<1.5.\label{eq:OptimalValueForj}
\end{equation}
\stopstep

\nextpar{The rough solution is smooth}
We will use that for each $m\geq 1$ there exists a constant $C_m$ such that
\begin{equation}1+|u_n-\overline u|_m \leq C_m A^{\beta \kappa^n}\label{ineq:1+|u(n)-ubar|m}\end{equation}
which is proven below.
The important point is that the inequality holds with the same $\beta$ regardless of $m$.
\\

Fix then $i\geq 1$, and let $m\geq 1$ which will be determined later.
Denote $M_m$ (or make $M_j$ larger if $j-1=m$) a constant such that the following estimates 
hold for any $v\in {\mathcal Y}$, $u\in  {\mathcal B}$, and $t>0$:
\begin{eqnarray}
|S(t)v|_m&\leq&M_m |v|_m,\\
|L(u)P(g,u)|_m&\leq&M_m(1+|u-\overline u|_{m+{1}}).
\end{eqnarray}
By interpolation inequalities, 
($C$ denotes constants depending on $i, m$, but not on $n$)
\begin{eqnarray}
|u_{n+1}-u_n|_i &\leq& |u_{n+1}-u_n|_1^\frac{m-i}{m-1}|u_{n+1}-u_n|_m^\frac{i-1}{m-1}\\
&\leq&CA^{-\mu\kappa^{n+1}\frac{m-i}{m-{1}}}|S(t_n)L(u_n)P(g,u_n)|_m^\frac{i-{1}}{m-{1}}\\
&\leq&CA^{-\mu\kappa^{n+1}\frac{m-i}{m-{1}}} |L(u_n)P(g,u_n)|_m^\frac{i-{1}}{m-{1}}\\
&\leq&CA^{-\mu\kappa^{n+1}\frac{m-i}{m-{1}}} (1+|u_n|_{m+{1}})^\frac{i-{1}}{m-{1}}\\
&\leq&CA^{-\mu\kappa^{n+1}\frac{m-i}{m-{1}}} A^{\beta\kappa^n\frac{i-{1}}{m-{1}}}\\
&\leq&CA^{(-\mu\kappa(m-i)+\beta(i-{1}))\frac{\kappa^n}{m-{1}}}.
\end{eqnarray}
Now choosing $m$ sufficiently large that
\begin{equation}-\mu\kappa(m-i)+\beta(i-{1})<0\label{ineq:ConditionOnm}
\end{equation}
makes the exponent negative and the increment $|u_{n+1}-u_n|_i$ 
decays with a double exponential rate.  It is then easy to see that $u_n$ is Cauchy in 
the $|\cdot|_i$-norm.
\stopstep %

\nextpar{Proof of estimates (\ref{ineq:1+|u(n)-ubar|m})}
Fix $m\geq {1}$.
As for (\ref{ineq:InductionInequalityOn(1+|u(n)-ubar|j)}) we find
\begin{eqnarray}
1+|u_{n+1}-\overline u|_m   
&\leq&2M_m^2A^{\kappa^n}(1+|u_n-\overline u|_m)
\end{eqnarray}
so that
\begin{eqnarray}
A^{-\beta\kappa^{n+1}}(1+|u_{n+1}-\overline u|_m)
\leq 2M_m^2 A^{(-\beta(\kappa-1)+1)\kappa^n}\left(A^{-\beta\kappa^n}(1+|u_n-\overline u|_m)\right).
\end{eqnarray}
Given $m, \beta, \kappa$, 
let $n^*(m)$ such that for $n\geq n^*(m)$, $2M_m^2A^{(-\beta(\kappa-1)+1)\kappa^n}\leq 1$
(recall that $-\beta(\kappa-1)+1<0$). 
Then choose $C_m$ so that $A^{-\beta\kappa^n}(1+|u_n-\overline u|_m)\leq C_m$  for $n<n^*(m)$.
\stopproof 

\subsection{The injective part of the Nash-Moser Inverse Function Theorem}
Even though the injective part of Nash-Moser Inverse Function Theorem cannot be used 
as such for the injective part of Theorem~\ref{theorem:MainTheorem},
it is instructive to show its proof as our result is an adaptation of it.
We follow \cite{HamiltonIFTNashMoser}.

The assumptions here are different from those of Theorem~\ref{theorem:SurjectivePartOfNashMoserIFT}.
We are not concerned with existence of solutions, only uniqueness.
Thus, $D_uP(g,u)$ is assumed \textit{injective}, with a \textit{left}-inverse again denoted $L(u)$.
Note that the Moser scheme plays no r\^ole here.
\begin{theorem}[Nash-Moser IFT - injective part]
Consider $g\in\mathcal Y$ in a $|\cdot|_0$-neighborhood of $\overline g$.
Suppose that $D_uP(g,u)v$ has a left-inverse  $L(u)h$, which is a tame of degree $0$ in $g$ and $h$, and $1$ in $u$.
Then, there exists a $|\cdot|_1$-neighborhood $\mathcal B':=\{u\in\mathcal X~|~|u-\overline u|_1<\eta'\}$
of $\overline u$ such that, \textbf{if} $P(g,u_1)=P(g,u_2)$ where $u_1,u_2\in\mathcal B'$,
then $u_1=u_2$.
\end{theorem}
\startproof
Use Taylor's formula,
\begin{equation}
P(g,u_2)=P(g,u_1)+D_uP(g,u_1)(u_2-u_1)+\int_0^1(1-t)D^2_{uu}P(g,u_1+t(u_2-u_1))(u_2-u_1,u_2-u_1)dt
\end{equation}
so that
\begin{equation}
u_2-u_1= -L(u_1)\int_0^1(1-t)D^2_{uu}P(g,u_1+t(u_2-u_1))(u_2-u_1,u_2-u_1)dt.
\end{equation}
Tame estimates on $D^2_{uu}P(u)(v_1,v_1)$ give
\begin{equation}
|u_2-u_1|_0\leq c|u_2-u_1|_0|u_2-u_1|_1
\end{equation}
where the constant is independent of $u_1, u_2\in\mathcal B'$ and $g$ in the restricted neighborhood.
Making $\eta'$ sufficiently small, we can make $c|u_2-u_1|_1<1$ for any $u_1,u_2\in\mathcal B'$.
This forces $|u_2-u_1|=0$.
\stopproof%

\subsection{Examples of smooth tame maps}\label{section:ExamplesAndUsefulRemarks}
We list in this Section some smooth tame maps which are used throughout the present work.
In this Section, by smoothness we do mean that derivatives of all orders exist.
$K$ denotes a compact subset of Euclidean space with smooth boundary.
$\mathcal V$ denotes an open subset of some Fr\'echet space.

\setcounter{counterpar}{0}
\nextpar{Linear differential operators with constant coefficients}
\begin{lemma}\label{lemma:LinearDifferentialOperatorsAreSmoothTame}
A linear differential operator of order $r$  with
constant coefficients $L\colon C^\infty_K\rightarrow C^\infty_K$ is a smooth tame  map of Fr\'echet spaces:
$L\colon C^{n+r}_K\rightarrow C^{n}_K$ is continuous for each $n\geq 0$.
$Lu$ has degree $r$ and base $0$: for all $u\in C^\infty_K$.
\begin{equation}\|Lu\|_{n,\alpha}\leq C\cdot \|u\|_{n+r,\alpha},\qquad n\geq 0.\end{equation}
\end{lemma}

\nextpar{Product of functions}
\begin{lemma}\label{lemma:ProductOfFunctionsIsSmoothTame}
The bilinear map  
\begin{equation}
B\colon\left\{\begin{array}{ccc}C^\infty_K&\times& C^\infty_K\\F&&G\end{array}
\longrightarrow\begin{array}{c}C^\infty_K\\FG\end{array}\right\}
\end{equation}
is a smooth tame map of Fr\'echet spaces:
for each $n\geq 0$, $B\colon C^n_K\times C^n_K\rightarrow C^n_K$ is continuous 
as well as $B\colon C^{n,\alpha}_K\times C^{n,\alpha}_K\rightarrow C^{n,\alpha}_K$.
$B(F,G)$ has degree $0$ in $F$ and $G$, and base $0$:
\begin{equation}\|B(F,G)\|_{n,\alpha}\leq C\cdot (\|G\|_{0,\alpha}\|F\|_{n,\alpha}+\|F\|_{0,\alpha}\|G\|_{n,\alpha}),\qquad n\geq 0
\label{ineq:TameEstimatesForProductOfFunctions}
\end{equation}
for all $F$ and $G$.
The first derivative is given by
\begin{equation}B(F,G)\cdot (f,g)=fG+Fg.\end{equation}
\end{lemma}

There are obvious generalizations of the above for the product of an arbitrary 
number of functions $(F_1,\dots, F_l)\mapsto F_1\cdots F_l$.\\

\nextpar{The Nemitskii operator}
If  $p(x,z)=p:K\times \mathbb R\rightarrow \mathbb R$ is a smooth function, define
$P(F)(x):=p(x,F(x))$, $x\in K$, $F\in \mathcal V\subset C^\infty_K$.
\begin{lemma}\label{lemma:NemitskiiIsSmoothTame}
$P:\mathcal V\rightarrow C^\infty_K$ is a smooth tame map of Fr\'echet spaces:
for each $n\geq 0$, $P\colon C^n_K\rightarrow C^n_K$ is continuous as well
as $P\colon C^{n,\alpha}_K\rightarrow C^{n,\alpha}_K$.
$P(F)$ has degree $0$ in $F$ and base $1$:
\begin{equation}
  \|P(F)\|_{n,\alpha}\leq C\cdot (\|F\|_{n,\alpha}+1),\qquad n\geq 1
\label{ineq:TameEstimatesForNemitskiiOperator}
\end{equation}
for $F$ in a neighborhood where $\|F\|_{1,\alpha}$ is bounded.
The first derivative $DP(F)\cdot f\in C^\infty_K$ is given by
\begin{equation}(DP(F)\cdot f)(x)=D_zp(x,F(x))f(x),\qquad x\in K.
\end{equation}
\end{lemma}

\nextpar{Composition of functions}
Let $K\subset \mathbb R^d, K'\subset \mathbb R^{d'}, K''\subset \mathbb R^{d''}$
be compact subsets of Euclidean spaces,
and let $G_0\in C^\infty_{(K,\mathbb R^{d'})}$ such that $G_0(K)\subset V'\subset K'$
for some open set $V'$.
If $G$ is in a suitable $\|\cdot\|_0$-neighborhood $\mathcal V$ of $G_0$, then
$G(K)\subset K'$.
Thus, we may define the \defn{composition operator}
\begin{equation}
C\colon\left\{\begin{array}{ccc}C^\infty_{(K',\mathbb R^{d''})}&\times&\mathcal V\\F&&G\end{array}
\longrightarrow\begin{array}{c}C^\infty_{(K,\mathbb R^{d''})}\\F\circ G\end{array}\right\}.
\end{equation}
\begin{lemma}\label{lemma:CompositionIsSmoothTame}
$C$ is a smooth tame map of Fr\'echet spaces.
In the $C^n$-grading we have that
\begin{equation}C\colon C^{n}_{(K',\mathbb R^{d''})}\times C^n_{K,\mathbb R^{d'}}\rightarrow C^n_{(K,\mathbb R^{d''})},
\qquad n\geq 0
\end{equation}
is continuous, while in the $C^{n,\alpha}$-grading we only have that 
\begin{equation}
C\colon C^{n+1,\alpha}_{(K',\mathbb R^{d''})}\times C^{n,\alpha}_{K,\mathbb R^{d'}}\rightarrow C^{n,\alpha}_{(K,\mathbb R^{d''})},
\qquad n\geq 0\end{equation}
is continuous.
For $G$ in a $\|\cdot\|_{1,\alpha}$-neighborhood, and all $F$ (without restriction),
\begin{equation}
\|F\circ G\|_{n,\alpha}\leq C\cdot(\|F\|_{n,\alpha}+\|G\|_{n,\alpha}\|F\|_{1,\alpha}),\qquad n\geq 1
\label{ineq:TameEstimatesForCompositionOperator}
\end{equation}
The derivative is given by
\begin{equation}DC(F,G)\cdot (f,g)= F'(G)g+f(G).
\end{equation}
\end{lemma}
\startproof
See \cite{delaLlaveObaya}. 
Note that composition $C^{n,\alpha}\times C^{n,\alpha}\rightarrow C^{n,\alpha}$ is well-defined,
even though it is not continuous.
\stopproof%

\nextpar{The inversion operator}
Let $I$ be a compact interval denote $\mathcal D^\infty_I$ the group of increasing $C^\infty$-diffeomorphisms of $I$.
Denote $V(F)=F^{-1}$ the inverse of $F\in\mathcal D^\infty_I$.
\begin{lemma}\label{lemma:InversionIsSmoothTame}
The inversion operator
\begin{equation}
V\colon\left\{\begin{array}{c}\mathcal D^\infty_I\\F\end{array}
\longrightarrow\begin{array}{c}\mathcal D^\infty_I\\F^{-1}\end{array}\right\}
\end{equation}
is a smooth tame map.
It is continuous $\mathcal D^n_I\rightarrow \mathcal D^n_I$ for each $n\geq 1$.
$V(F)$ has degree $0$ in $F$ and base $1$ in the $\|\cdot\|_{n,\alpha}$-grading:
\begin{equation}\|V(F)\|_{n,\alpha}\leq C\cdot (\|F\|_{n,\alpha}+1),\qquad n\geq 1
\label{ineq:TameEstimatesForInversionOperator}
\end{equation}
for $F$ in a neighborhood where $\|F\|_{1,\alpha}$ is bounded.
The first derivative is given by
\begin{equation}DV(F)\cdot f=-\frac{f(F^{-1})}{F'(F^{-1})}.
\end{equation}
\end{lemma}
\startproof%
That $V\colon \mathcal D^n_I\rightarrow \mathcal D^n_I$ is continuous for each $n\geq 1$ is standard, 
see Example~4.4.6, p.~92, Part~I, and Theorem~2.3.5, p.~148, Part~II of \cite{HamiltonIFTNashMoser}.
We need tame estimates in the $C^{n,\alpha}$, while those of \cite{HamiltonIFTNashMoser} are given in the $C^n$-grading.

Let $y=g(x)$ be a diffeomorphism in $\mathcal D^\infty_I$ with inverse $x=f(y)$.
Then, $f'(y)=\frac1{g'(x)}$ so that 
\begin{eqnarray}
|f'(y_1)-f'(y_2)|&\leq& \left|\frac1{g'(x_1)}-\frac1{g'(x_2)}\right|\\
&\leq&\frac{|g'(x_2)-g'(x_1)|}{|g'(x_1)g'(x_2)|}\\
&\leq& C[g']_\alpha|x_2-x_1|^\alpha
\end{eqnarray}
for $g$ in a $\|\cdot\|_1$-neighborhood.
Since $|x_2-x_1|\leq \|f'\|_0|y_2-y_1|$, we are done with the tame estimates for $n=1$.

Suppose the tame estimates verified for $1\leq m<n$.
From (\ref{eq:nThDerivativeOff}) and tame estimates on product of functions, we have for $g$ in a $\|\cdot\|_{1,\alpha}$-neighborhood
\begin{eqnarray}
\|f^{(n)}\|_{0,\alpha}
&\leq& C\|f'\|_{0,\alpha}^n\sum_{k=1}^{n-1}\sum_{\begin{subarray}{c}j_1+\dots+j_k=n\\j_1,\dots,j_k\geq 1\end{subarray}}\|f^{(k)}\|_{0,\alpha}
\|g^{(j_1)}\circ f\|_{0,\alpha}\cdots \|g^{(j_k)}\circ f\|_{0,\alpha}\\
&\leq& C\|f'\|_{0,\alpha}^n\sum_{k=1}^{n-1}\sum_{\begin{subarray}{c}j_1+\dots+j_k=n\\j_1,\dots,j_k\geq 1\end{subarray}}\|f^{(k)}\|_{0,\alpha}
\|g^{(j_1)}\|_{0,\alpha}\|f\|_{0,\alpha}\cdots \|g^{(j_k)}\|_{0,\alpha}\|f\|_{0,\alpha}\\
&\leq& C\sum_{k=1}^{n-1}\sum_{\begin{subarray}{c}j_1+\dots+j_k=n\\j_1,\dots,j_k\geq 1\end{subarray}}(1+\|g^{(k)}\|_{0,\alpha})
\|g\|_{j_1,\alpha}\cdots \|g\|_{j_k,\alpha}.
\end{eqnarray}
We interpolate each factor between their $\|\cdot\|_{1,\alpha}$- and $\|\cdot\|_{n,\alpha}$-norms:
\begin{eqnarray}
&&\|g\|_{k,\alpha}\|g\|_{j_1,\alpha}\cdots \|g\|_{j_k,\alpha}\\
&\leq& C\cdot \|g\|_{1,\alpha}^\frac{n-k}{n-1}\|g\|_{n,\alpha}^\frac{k-1}{n-1}\cdot
\|g\|_{1,\alpha}^\frac{n-j_1}{n-1}\|g\|_{n,\alpha}^\frac{j_1-1}{n-1}\cdots
\|g\|_{1,\alpha}^\frac{n-j_k}{n-1}\|g\|_{n,\alpha}^\frac{j_k-1}{n-1}\\
&\leq&C\cdot \|g\|_{n,\alpha}
\end{eqnarray}
since $\|g\|_{1,\alpha}$ remains bounded and $j_1+\cdots+j_k=n$.
The other terms in the sum are handled similarly.
\stopproof%

\bibliography{the-big-paper}


\end{document}